\documentclass[11pt]{amsart}	
\usepackage[utf8]{inputenc}
\usepackage{amsmath,amsfonts,amssymb,dsfont}
\usepackage{lmodern}
\usepackage{graphicx}
\usepackage[english]{babel}
\usepackage{listings}
\usepackage[a4paper]{geometry}
\usepackage[T1]{fontenc}
\usepackage{xspace}
\usepackage[colorinlistoftodos]{todonotes}
\usepackage{pgf,pgfkeys,pgfmath}
\usepackage{upgreek}
\usepackage{fullpage}

\def\one#1{\mathds{1}_{#1}}

\usepackage[backref=page]{hyperref}
\hypersetup{
	plainpages=false,
	unicode=false,          
	pdftoolbar=true,        
	pdfmenubar=true,        
	pdffitwindow=true,      
	pdftitle={My title},    
	pdfauthor={Author},     
	pdfsubject={Subject},   
	pdfcreator={Creator},   
	pdfproducer={Producer}, 
	pdfkeywords={keywords}, 
	pdfnewwindow=true,      
	colorlinks=true,       
	linkcolor=blue,          
	citecolor=blue,        
	filecolor=blue,      
	urlcolor=magenta,          
	urlbordercolor=0 1 1
}
\newcommand*{\Scale}[2][4]{\scalebox{#1}{\ensuremath{#2}}}%
\newcommand{\bigmu}{\Scale[1.3]{\upmu}}
\def\div{\operatorname{div}}
\def\dist{\operatorname{dist}}

\newtheorem{thm}{Theorem}[section]

\newtheorem{prop}[thm]{Proposition}
\newtheorem{lem}[thm]{Lemma}

\usepackage{todonotes}

\newcommand\numeq[1]%
  {\stackrel{\scriptscriptstyle(\mkern-1.5mu#1\mkern-1.5mu)}{=}}

\newcommand{\myqed}{\thinspace\null\nobreak\hfill\hbox{\vbox{\kern-.2pt\hrule
height.2pt depth.2pt\kern-.2pt\kern-.2pt \hbox to2.5mm{\kern-.2pt\vrule
width.4pt \kern-.2pt\raise2.5mm\vbox to.2pt{}\lower0pt\vtop to.2pt{}\hfil
\kern-.2pt \vrule width.4pt\kern-.2pt}\kern-.2pt\kern-.2pt\hrule
height.2pt depth.2pt \kern-.2pt}}\par\medbreak}

\definecolor{lbcolor}{rgb}{0.95,0.95,0.95}
\definecolor{cblue}{rgb}{0.,0.0,0.6}

\newcommand{\ordres}[2]{\ensuremath{(\mathcal{O}(#1), \mathcal{O}(#2) )}}

\definecolor{Lightgray}{rgb}{0.85, 0.85, 0.85}

\title{Approximation of surface diffusion flow: a second order variational Cahn--Hilliard model with degenerate mobilities}
\author{Elie Bretin}
\address{Univ Lyon, INSA de Lyon, CNRS UMR 5208, Institut Camille Jordan\\ 20 avenue Albert Einstein, F-69621 Villeurbanne, France\\ elie.bretin@insa-lyon.fr}
\author{Simon Masnou}
\address{Univ Lyon, Universit\'e Claude Bernard Lyon 1, CNRS UMR 5208, Institut Camille Jordan \\43 boulevard du 11 novembre
1918, F-69622 Villeurbanne, France\\masnou@math.univ-lyon1.fr}
\author{Arnaud Sengers}
\address{Univ Lyon, Universit\'e Claude Bernard Lyon 1, CNRS UMR 5208, Institut Camille Jordan\\43 boulevard du 11 novembre
1918, F-69622 Villeurbanne, France\\sengers@math.univ-lyon1.fr}
\author{Garry Terii}
\address{Univ Lyon, Universit\'e Claude Bernard Lyon 1, CNRS UMR 5208, Institut Camille Jordan\\43 boulevard du 11 novembre
1918, F-69622 Villeurbanne, France\\ terii@math.univ-lyon1.fr}
\makeatletter
\@namedef{subjclassname@2020}{\textup{2020} Mathematics Subject Classification}
\makeatother
 
\begin{document} 
\maketitle
\begin{abstract}

This paper tackles the approximation of surface diffusion flow using a Cahn--Hilliard-type model. We introduce and analyze a new second order variational phase field model which associates the classical Cahn--Hilliard energy with two degenerate mobilities. This association allows to gain an order of approximation of the sharp limit.
 In a second part, we propose some 
simple and efficient numerical schemes to approximate the solutions,
and we provide numerical 2D and 3D experiments that illustrate the interest
of our model in comparison with other Cahn--Hilliard models. 

\end{abstract}
\section{Introduction}

This paper addresses the approximation of surface diffusion flow, which is the evolution of a time-dependent surface $\Gamma:\,t\mapsto\Gamma(t)$ moving with normal velocity at every time $t$:
$$ V(t) =  \Delta_{\Gamma(t)} H(t),$$
where $H(t)$ is the mean curvature vector on $\Gamma(t)$, and $\Delta_{\Gamma(t)}$ the Laplace-Beltrami operator defined on the surface. For simplicity, we shall frequently omit the time dependence.\\
The starting point of our approximation model is the classical Cahn--Hilliard equation
$$\varepsilon^2\partial_t u  = \Delta \left(  W'(u) - \varepsilon^2 \Delta u \right),$$  
where  $u:\,(t,x)\mapsto u(t,x)$ is a smooth function whose level surface $\{u(t,\cdot)=\frac 1 2\}$ approximates $\Gamma(t)$, $\varepsilon>0$ is a small parameter, and $W$ is a reaction potential, typically $W(s) = \frac{1}{2}s^2 (1-s)^2$.

The Cahn--Hilliard equation has been introduced as a mathematical model for phase separation and phase coarsening 
in binary alloys \cite{cahnspin,CAHN:1958}, but it has also been used for applications as diverse as the modeling of two evolving components of intergalactic material or the description of a bacterial film, see the references in~\cite{novick2008cahn}, or the modeling of multiphase fluid flows~\cite{MuliflowCH1,MuliflowCH2}. More recently it was proposed as an inpainting model in image processing, see~\cite{inpaintCH,inpaintCH2,inpaintCH3}. We refer to \cite{novick2008cahn} for an inspiring general introduction to the Cahn--Hilliard equation, see also the recent book \cite{bookMiranville} where state-of-art results and many applications of the Cahn--Hilliard equation are presented.
   
\medskip
\paragraph*{\bf \bf Sharp interface limit and mobilities}~   \\
Pego determined with formal arguments in \cite{pego1989front}, and Alikakos et al proved rigorously in \cite{alikakos1994convergence}, that the sharp limit flow of the Cahn--Hilliard equation (for suitable time regimes as $\varepsilon\to 0$) is the Mullins-Sekerka interface motion. 

Observe now that the Cahn--Hilliard equation can be equivalently written as 
\begin{equation}\label{CH-mob}
\varepsilon^2\partial_t u  = \div\left(M(u)\nabla (W'(u) - \varepsilon^2 \Delta u)\right)
\end{equation}
with the particular choice $M(u)\equiv 1$.  If $M$ is now chosen to be non constant, it plays the role of a concentration-dependent mobility. Cahn et al. showed formally in~\cite{cahn1996cahn} that if one uses a degenerate mobility $M(u) = u(1-u)$ (degenerate in the sense that there is no motion where $u=0$ or $1$) and a logarithmic potential
$$W(s) = \frac{1}{2}\theta \left[ s\ln(s) + (1-s)\ln(1-s) \right] + \frac{1}{2}s(1-s),$$
the sharp limit motion is the surface diffusion flow. However, the singularity of such a logarithmic potential makes the model not well suited for numerical simulations. We shall see in this paper that a different model can be proposed which leads to the surface diffusion flow as well, but involves rather the smooth potential 
$$W(s) = \frac{s^2(1-s)^2}{2}.$$

The choice of appropriate degenerate mobility and potential is important. It was observed in the review paper~\cite{gugenberger2008comparison} that some choices lead to inconsistencies, in the sense that depending on how terms are identified in the matched asymptotic analysis expansion, one can either show the convergence to surface diffusion flow, or to a stationary flow with null velocity. The authors of \cite{lee2015degenerate,lee2016sharp} suggested
that such inconsistencies come from the presence of an additional bulk diffusion term is the limit motion, i.e. the limit velocity is:
\begin{equation*}
V = \frac{2}{3}\Delta_\Gamma H + \alpha H \nabla_n H
\end{equation*}
This has been corroborated numerically in \cite{dai2012motion,dai2014coarsening} where undesired coarsening
effects are observed. The additional term in the velocity depends on the derivative of the mobility $M'(u_0)$, where $u_0$ is the outer solution in the matched asymptotics which equal $0$ or $1$. To obtain a pure motion by surface 
diffusion, one needs to take a higher order mobility, for example $M(s) = s^2(1-s)^2$~\cite{lee2016sharp}. With such a choice, the bulk diffusion appears in 
higher order terms and the correct velocity is recovered (with a different multiplicative constant):
\begin{equation*}
V = \alpha\Delta_\Gamma H
\end{equation*}
These conclusions have been extended to the anisotropic case in~\cite{dziwnik2017anisotropic}. 

\medskip
\paragraph*{\bf \bf Positivity property  and order of phase field model}~\\
We now turn to the following question: starting from an initial $u(0)$ with values in $[0,1]$ and using the above mobility, does the solution $u$ remain valued in $[0,1]$ ? 
This is often referred to as the \emph{positivity condition} as it implies that all phase functions remain positive in a multiphase
context. This condition is important also because it means the function remains within the pure state phase boundaries. 

The theoretical results of \cite{lee2016sharp,salvalaglio2019doubly} and the numerical evidences
of \cite{dai2012motion,dai2014coarsening} establish that it is not 
the case with the Cahn--Hilliard model with mobility~\eqref{CH-mob}.
More precisely the profile of the solution shows some oscillations when reaching the pure states. 
This comes from the influence in the asymptotic expansion of the solution
of the first order error term which does not vanish for this kind of phase field models. 

To circumvent this problem a non variational model has been introduced in \cite{ratz2006surface}:
\begin{equation}\label{ratz}
\left\lbrace
\begin{aligned}
& \varepsilon^2 \partial_t u = \div\left( M(u) \nabla \mu \right) \\ 
& g(u)\mu = W'(u) - \varepsilon^2 \Delta u
\end{aligned}
\right.
\end{equation}
with $g(s) = \gamma |s|^p|1-s|^p$, $p\geq 0$.

The idea is to add another degenerate term $g$ that acts as a diffusion preventing term and forces the aforementioned error 
term to be smaller and to converge to zero far from the interface. 
This model is known to achieve better numerical accuracy than the classical model~\eqref{CH-mob}, and it has been successfully adapted in various applications, see for example \cite{albani2016dynamics,naffouti2017complex,salvalaglio2017morphological,salvalaglio2015faceting}.

Several choices have been made for $p$, the most acclaimed ones being $p=1$ and $p=2$ but they were motivated by better numerical 
results rather than from a theoretical standpoint. In this paper we explain why $p=1$ is the correct choice as 
it imposes the leading error term of the solution to be zero. This result is new to the extent of our knowledge.
 
While it has excellent numerical properties, the above model~\eqref{ratz} does not derive from an energy and it is thus more difficult to prove
rigorously theoretical properties and to extend the model to complex multiphase applications. Therefore, a variational adaptation is proposed in \cite{salvalaglio2019doubly}:
\begin{equation*}
\left\lbrace
\begin{aligned}
& \varepsilon^2 \partial_t u = \div\left( M(u) \nabla \mu \right) \\ 
& \mu = g(u) W'(u) - \varepsilon^2 \div\left(g(u)\nabla u\right) + g'(u) \left( W'(u) + \frac{\varepsilon^2}{2} |\nabla u |^2 \right)
\end{aligned}
\right.
\end{equation*}
The idea is to inject the second degeneracy $g$ in the energy. 
The model conserves the same advantages as the non variational version, in particular the fact that the choice $p=1$ remains
the correct one and nullifies the leading error term of the solution. However it relies on changing the energy, 
thus making it harder to extend to complex multiphase application or to add an anisotropy. Also, it seems more appropriate to incorporate 
the mobility in the metric rather than in the geometry of the evolution problem. This is what we propose in this paper.

\medskip
\paragraph*{\bf A new variational Cahn--Hilliard model of order two}~ \\
Like in \cite{salvalaglio2019doubly}, we want to approximate the
surface diffusion flow using a second order variational phase field model, but we want it closer to the original Cahn--Hilliard model. \\

The new Cahn--Hilliard model we propose reads as:
$$\begin{cases}
   \varepsilon^2 \partial_t u &= N(u) \div\left(M(u) \nabla (N(u)\mu) \right) \\
 \mu &= W'(u) - \varepsilon^2 \Delta u,\\
  \end{cases}
$$ 
and at least in the case where
$$ W(s) = \frac{1}{2} s^2(1-s)^2, \quad  M(s) = s^2(1-s)^2 \text{ and }  N(s) = \frac{1}{s(1-s)},$$
we will show that this model is of order two and converges to the surface diffusion flow. 
Consequently, this model has all desired properties while conserving the correct energy to dissipate.

Moreover, as the conservation of volume is one key feature of the Cahn--Hilliard equation, 
we review how well each model manages to preserve this property.
Furthermore, because of the higher consistency of the solution profile, 
we will show that we achieve very good numerical approximation in this area as in
\cite{Bretin_brassel,bretin_largephases}.

\medskip
\paragraph*{\bf Outline of the paper:}
The paper is organized as follows. First, we review in Section \ref{Sec:Review} 
the properties of the Cahn--Hilliard model with mobility and the drawbacks that need to be improved. In Section \ref{Sec:ModelNMN} we present our new variational Cahn--Hilliard model and review its properties.
We prove these properties in Section \ref{Sec:DeuxMob} using the formal method of matched asymptotic expansion. 
The necessary tools are presented at the beginning of the proof. In the numerical section~\ref{Sec:Numerical}, 
we first explain how to derive a simple and efficient scheme using a convex splitting 
of Cahn--Hilliard energy and exploiting the variational mobility structure.  
Finally, we propose some numerical experiments that compare the various Cahn--Hilliard
models and highlight the advantages of our new model.

\section{Review of the properties of the Cahn--Hilliard equation with mobility}
\label{Sec:Review}
In this section we summarize the properties of the existing models and explain why we will
introduce a new one in the next section. A motion by surface diffusion can be obtained as the sharp limit of the 
Cahn--Hilliard equation with mobility if we choose a mobility $M$ that is of sufficiently high order. 

\subsection{The Cahn--Hilliard model with mobility to approximate surface diffusion flow}
We recall that the normal velocity associated with a surface diffusion flow is:
\begin{equation*}
V = \Delta_\Gamma H
\end{equation*}
We also recall that, if $\Omega$ denotes the inner domain enclosed by $\Gamma$, the phase field method consists in approximating the characteristic 
function $\one{\Omega}$ by a smooth function of the form $u_\varepsilon = q(\dist(\cdot,\Omega)/\varepsilon)$ where $q$ is the so-called {\it optimal profile} associated with the potential $W$, $\varepsilon$ represents the thickness of the smooth transition from $0$ to $1$, and $\dist$ denotes the signed distance function. The one associated with our choice $W(s) = \frac{1}{2} s^2(1-s)^2$ verifies the following properties:
\begin{equation}
\label{ProfilQ}
\left\lbrace
\begin{aligned}
& q(z) = \frac{1-\tanh(z)}{2} \\
& q'(z) = - \sqrt{2W(q)} \\
& q''(z) = W'(q)\\
\end{aligned}
\right.
\end{equation}
We denote:
\begin{equation*}
\left\lbrace
\begin{aligned}
& c_W = \int_\mathbb{R} (q'(z))^2 dz = \frac{1}{6}\\
& c_M = \int_\mathbb{R} \frac{M(z)}{q(z)(1-q(z))} dz\\
\end{aligned}
\right.
\end{equation*}
With the choice $M(s) = s(1-s)$, we have $c_M=1$ and with the choice $M(s) = s^2 (1-s)^2$, we have $c_M = \frac{1}{6}$. 
A higher order mobility $M$ will inevitably lower the constant in front of the velocity but will 
prove to be necessary to find the right motion as stated in Result \ref{Result1}.

We start off with the classical Cahn--Hilliard model with non negative mobility $M$ and potential $W(s) = \frac{s^2(1-s)^2}{2}$, that we refer to as {\bf M-CH} from now on:
\begin{equation}
\label{Eq:ClassicalCahnHilliard}
\left\lbrace
\begin{aligned}
& \varepsilon^2 \partial_t u = \div\left( M(u) \nabla \mu \right) \\ 
& \mu = W'(u) - \varepsilon^2 \Delta u
\end{aligned}
\right.
\end{equation}

When the mobility $M$ is a scalar positive weight independent of $u$, we recall that the equation
\begin{equation*}
\left\lbrace
\begin{aligned}
& \varepsilon^2 \partial_t u = \div\left( M \nabla \mu \right) \\ 
& \mu = W'(u) - \varepsilon^2 \Delta u
\end{aligned}
\right.
\end{equation*}
is the $H^{-1}$ gradient 
flow of the Cahn--Hilliard energy
\begin{equation}
\label{Eq:EnergieCahnHilliard}
E(u) = \int_Q \frac{\varepsilon}{2} |\nabla u|^2 + \frac{1}{\varepsilon} W(u) \ \ dx
\end{equation}
when considering the following scalar product in $H_0^1$ weighted by the mobility $M$:
\begin{equation}
\label{Eq:H01MobilityScalar}
\langle f, g\rangle_{H_0^1} = \int_Q M \nabla f \cdot \nabla g  \ \ dx.
\end{equation}
It is important to note that the mobility is incorporated in the metric with respect to which the gradient flow
is computed, and not as a geometric parameter in the energy. 

Equation~\eqref{Eq:ClassicalCahnHilliard} is an extension of the above equation to the case where $M$ depends on $u$.

\medskip
The {\bf M-CH} model has been extensively studied and it is well understood that
the mobility needs to be a polynomial of order at least 2. Indeed, a mobility of order 1 would give a quicker motion, but as already mentioned 
the authors of \cite{lee2016sharp} showed that it yields an additional undesired bulk diffusion term in the limit velocity.
This term can be removed the pure surface diffusion motion can be recovered by choosing a higher order mobility, which is what we will do.

\subsection{Properties of the classical Cahn--Hilliard model with mobility}
The properties of {\bf M-CH} are summarized in the following result, see~\cite{lee2016sharp}:
\begin{prop}
\label{Result1}
With the choice $M(s) = s^2(1-s)^2$, the solution $u_\varepsilon$ to \eqref{Eq:ClassicalCahnHilliard} expands formally near the interface $\Gamma_\varepsilon(t)= \left\lbrace u_\varepsilon(t,\cdot) =\frac{1}{2} \right\rbrace$ as:
\begin{equation}
\label{Eq:ProfileCH}
u_\varepsilon(t) = q\left( \frac{\dist(x,\Omega_\varepsilon(t))}{\varepsilon} \right) + \mathcal{O}(\varepsilon) \\
\end{equation}
with $\Omega_\varepsilon(t)= \left\lbrace u_\varepsilon(t,\cdot) \leq \frac{1}{2} \right\rbrace$. The associated normal velocity satisfies:
\begin{equation}
V = c_W c_{M} \Delta_{\Gamma} H + \mathcal{O}(\varepsilon)
\end{equation}
If we do not require $M'(0)=M'(1)=0$, i.e. $M$ to be a double well polynomial mobility with roots $0$ and $1$ of multiplicity at least two,
then the velocity contains an additional undesired bulk diffusion term. For example, if we set $M(s)=s(1-s)$, then
\begin{equation}
V = c_W c_M \Delta_{\Gamma} H \pm c_W^2 c_M H \nabla_n H 
\end{equation}
Moreover, the volume is preserved only up to an order $\mathcal{O}(\varepsilon)$:
\begin{equation}
|\Omega_\varepsilon(t)| = |\Omega_\varepsilon(0)| + \mathcal{O}(\varepsilon)
\end{equation}
\end{prop}

\bigskip
\noindent From now on, we fix the mobility $M$ to be 
$$M(s) = s^2(1-s)^2.$$
{\bf M-CH} has a well identified drawback. The leading error term  in \eqref{Eq:ProfileCH}
is of order $\varepsilon$ and has a dependence in the curvature. 
This means it becomes relevant in high curvature regions. This is especially problematic when reaching the 
pure states because an overshoot due to oscillations occurs, and the solution does not
stay within its physical range $[0,1]$. This problem proves to be even more problematic in the multiphase
case because the solutions may not be positive anymore and phantom phases may appear.

The volume conservation is a standard property of the Cahn--Hilliard model on a domain $Q$ with periodic or Neumann homogeneous boundary conditions on $\partial Q$:
\begin{equation*}
\frac{d}{dt} \int_Q u \ dx  = \int_Q \partial_t u \ dx  = \int_Q\div\left(M(u)\nabla \mu \right)dx = 0.
\end{equation*}

However, numerically, the quality of the conservation is constrained by the quality of the approximation
of the solution $u_\varepsilon$. As we will see with the later models, a more accurate solution $u_\varepsilon$
will also lead to a more accurate conservation of the volume.

\section{A new variational model with two mobilities}
\label{Sec:ModelNMN}
In this section, we propose a new variational Cahn--Hilliard model with two mobilities.
In contrast with \cite{salvalaglio2019doubly} where the energy is modified, we propose to  
incorporate the additional degeneracy in the metric used  for defining the gradient flow. First, we derive our model and explain the right choice for its parameters. Then we review its theoretical properties, that appear to be similar as those of the previous model when we choose $p=1$. We compare the numerical behavior of each method in the next section devoted to  numerics.
\subsection{Derivation of the model}
Our model derives from the classical Cahn--Hilliard energy:
\begin{equation*}
E(u) = \int_Q \frac{\varepsilon}{2}|\nabla u|^2 + \frac{1}{\varepsilon} W(u) dx
\end{equation*}
Let us consider a $H_0^1$ scalar product with two scalar positive weights $M$ and $N$:
\begin{equation*}
\langle f,g\rangle_{H_0^1} = \int_Q M \nabla (Nf)\cdot \nabla(Ng) dx
\end{equation*}

Taking the $H^{-1}$ gradient flow of the energy $E$ with respect to this scalar product, we obtain the following equation: 
\begin{equation*}
\left\lbrace
\begin{aligned}
& \varepsilon^2 \partial_t u = N \div(M\nabla(N\mu)) \\
& \mu = -\varepsilon^2 \Delta u + W'(u)
\end{aligned}
\right.
\end{equation*}

Considering now a dependence on $u$ of $M$ and $N$ gives the following equation, that we refer to as the {\bf NMN-CH} model:

\begin{equation}
\label{Eq:modelNMN}
\left\lbrace
\begin{aligned}
& \varepsilon^2 \partial_t u = N(u) \div(M(u)\nabla(N(u)\mu)) \\
& \mu = -\varepsilon^2 \Delta u + W'(u)
\end{aligned}
\right.
\end{equation}

This model has the advantage of conserving the minimized energy and
being variational in the following sense:
\begin{equation*}
\left.
\begin{aligned}
\frac{d}{dt}E(u) & = \int_Q \left(-\varepsilon \Delta u + W'(u) \right) \partial_tu dx \\
& =  \int_Q N(u)\mu \div\left(M(u)\nabla(N(u)\mu) \right) dx\\
& =  - \int_Q M(u) |\nabla(N(u)\mu)|^2 dx + \int_{\partial\Omega} M(u)N(u)\mu \nabla(N(u)\mu)\cdot n d\sigma \\
& = -\int_Q M(u) |\nabla(N(u)\mu)|^2 dx \leq 0
\end{aligned}
\right.
\end{equation*}
in the case of periodic or Neumann homogeneous boundary conditions on $\partial Q$.

\subsection{Choosing $N$} 
As previously stated, $M$ is set to be:
\begin{equation*}
M(s) = s^2(1-s)^2.
\end{equation*}
We want to choose $N$ so that it has an antagonist effect to $M$ and forces the leading error term $U_1$ to be zero, see below. We will show that the correct choice for $N$ is:
\begin{equation}
\label{Eq:CorrectN}
N(s) = \frac{1}{\sqrt{M(s)}} =  \frac{1}{s(1-s)}
\end{equation}
Indeed, the following equation is obtained for $U_1$ (see details below):
\begin{equation*}
\partial_{zz} U_1 - W''(U_0) U_1 = H \partial_z U_0 - \bigmu_1
\end{equation*}
and with the above choice for $N$ we have:
\begin{equation*}
H \partial_z U_0 = \bigmu_1
\end{equation*}
Thus $U_1=0$ while other choices for $N$ only impose $U_1\rightarrow 0$ far from the interface.

We define the following integral:
\begin{equation*}
c_N = \int_{-\infty}^{+\infty} \frac{q'(z)}{N(q(z))} dz
\end{equation*}
In conclusion, the correct choice for $N$ is~\eqref{Eq:CorrectN}.

\subsection{Properties of the {\bf NMN-CH} model }

They are summarized in the following result, to be compared with Proposition~\ref{Result1}.
\begin{prop} \label{Result4} If we choose $M(s) = s^2(1-s)^2$ and $N(s) = \frac{1}{\sqrt{M(s)}} = \frac{1}{s(1-s)}$, the solution $u_\varepsilon$ to \eqref{Eq:modelNMN} expands formally near the interface $\Gamma_\varepsilon(t)$ as:
\begin{equation}
\label{Eq:SolutionNMN}
u_\varepsilon = q\left( \frac{\dist(x,\Omega_\varepsilon(t))}{\varepsilon} \right) + \mathcal{O}(\varepsilon^2) \\
\end{equation}
with $\Omega_\varepsilon(t)= \left\lbrace u_\varepsilon \leq \frac{1}{2} \right\rbrace$. The associated normal velocity satisfies:
\begin{equation}
\label{Eq:VelocityNMN}
V = \frac{c_{W} c_{M}}{c_N^2} \Delta_{\Gamma} H + \mathcal{O}(\varepsilon)
\end{equation}
Moreover, the volume is preserved up to an order $\mathcal{O}(\varepsilon^2)$:
\begin{equation}
\label{Eq:VolumeNMN}
|\Omega_\varepsilon(t)| = |\Omega_\varepsilon(0)| + \mathcal{O}(\varepsilon^2)
\end{equation}
\end{prop}
\section{Proof of Proposition~\ref{Result4}}
\label{Sec:DeuxMob}
In this section we prove Proposition~\ref{Result4} which summarizes the properties of {\bf NMN-CH}.
We start with the volume conservation~\eqref{Eq:VolumeNMN}, assuming the other properties as in \cite{Bretin_brassel}. 
Then, we introduce the tools and 
notations to derive the formal asymptotics and demonstrate \eqref{Eq:SolutionNMN} and \eqref{Eq:VelocityNMN}.
\subsection{Proof of the volume conservation}
In this part, we demonstrate \eqref{Eq:VolumeNMN} when assuming the profile \eqref{Eq:SolutionNMN}, which we will prove in the next part. 
We recall the following relations linking $W$ with $M$ and $N$:
\begin{equation*}
W(s) = \frac{1}{2}s^2(1-s)^2, \quad M(s) = 2W(s),\quad \text{ and }\; N(s) = \frac{1}{\sqrt{2W(s)}}.
\end{equation*}
The proof is done in two steps. First, we give the expression of the volume in terms of an integral of the function $G(s) = \int_0^s \sqrt{2W(s)ds}$ and use it to show the volume conservation of the {\bf NMN-CH}  model using the form of the profile $u$ given by \eqref{Eq:SolutionNMN}. Then we show the validity of this expression to conclude the proof.

The formula linking the volume of $\Omega_\varepsilon(t)$ with $G$ is the following:
\begin{equation}
\label{Eq:ProofVolumeExpression}
|\Omega_\varepsilon(t)| = \int_Q 6 (G\circ q) \left(\frac{d(x,\Omega_\varepsilon(t))}{\varepsilon} \right) dx + \mathcal{O}(\varepsilon^2)
\end{equation}
where $d$ is the signed distance function to the interface of $\Omega$ and $d(x,\Omega)<0$ for $x\in \Omega$.

Under the assumption that the profile $u_\varepsilon$ is given by:
\begin{equation*}
u_\varepsilon(x,t) = q\left( \frac{d(x,\Omega_\varepsilon(t)}{\varepsilon} \right) + \mathcal{O}(\varepsilon^2),
\end{equation*}
we have by composition by $G$ and integration:
\begin{equation*}
\int_{\mathbb{R}^d} G(u_\varepsilon(x,t)) = \int_{\mathbb{R}^d}(G\circ q)\left( \frac{d(x,\Omega_\varepsilon(t)}{\varepsilon} \right) + \mathcal{O}(\varepsilon^2)
\end{equation*}
Using \eqref{Eq:ProofVolumeExpression}, we conclude:
\begin{equation*}
\forall t\geq0, \left\vert \Omega_\varepsilon(t) \right\vert = \int_{\mathbb{R}^d} G\left(u_\varepsilon(x,t) \right)dx + \mathcal{O}(\varepsilon^2)
\end{equation*}

Considering periodic or Neumman boundary condition on $Q$ leads to a conservation of the integral of $G$ along the time:
\begin{equation*}
\left. \begin{aligned}
\frac{d}{dt} \int_Q G(u_\varepsilon) & = \int_Q G'(u_\varepsilon) \partial_t u_\varepsilon \\
			& = \frac{1}{\varepsilon^2} \int_Q \left(\sqrt{2W(u_\varepsilon)}N(u_\varepsilon) \right) \div\left( M(u_\varepsilon) \nabla (N(u_\varepsilon)\mu) \right) \\
			& = \frac{1}{\varepsilon^2} \int_Q \div\left( M(u_\varepsilon) \nabla (N(u_\varepsilon)\mu) \right) \\
			& = \frac{1}{\varepsilon^2} \int_{\partial Q} M(u_\varepsilon) \nabla \left(N(u_\varepsilon)\mu\right)\cdot n  = 0
\end{aligned} \right.
\end{equation*}
This means that the volume is conserved over time and \eqref{Eq:VolumeNMN} is verified if \eqref{Eq:ProofVolumeExpression} is satisfied.

We now turn to the proof of \eqref{Eq:ProofVolumeExpression}.
For the simplicity of the notations of the bounds of the integrals, we work in $\mathbb{R}^d$, 
but the result remains true for any regular bounded domain $Q$. Using the coarea formula, we have:
\begin{equation*}
\int_{\mathbb{R}^d} 6 (G\circ q) \left(\frac{d(x,\Omega_\varepsilon(t)}{\varepsilon} \right) = 6 \int_\mathbb{R} h(s) G\left(q\left(\frac{s}{\varepsilon}\right) \right) ds
\end{equation*}
where $h(s) = \left\vert D\chi_{\lbrace d(x,\Omega_\varepsilon(t)\leq s) \rbrace}\right\vert$ is the perimeter of the 
signed distance function to $\Omega_\varepsilon(t)$. Using the fact that:
\begin{equation*}
6G(q(-s)) = 6G(1-q(s)) = 1- 6 G(q(s))
\end{equation*}
We deduce:
\begin{equation*}
\left. \begin{aligned}
\int_{\mathbb{R}^d} 6 (G\circ q) \left(\frac{d(x,\Omega_\varepsilon(t)}{\varepsilon} \right) & = \int_{-\infty}^0 h(s) + \int h(s) \left( 6 G\left(q\left(\frac{s}{\varepsilon}\right)\right) -1\right) + 6 \int_0^{+\infty} h(s) G\left(q\left(\frac{s}{\varepsilon}\right) \right) dx \\
		& =  \left\vert \Omega_\varepsilon(t) \right\vert - \int_{-\infty}^0 h(s) G\left( q\left( \frac{-s}{\varepsilon}\right) \right) +  6 \int_0^{+\infty} h(s) G\left(q\left(\frac{s}{\varepsilon}\right) \right) dx\\
		& = \left\vert \Omega_\varepsilon(t) \right\vert + 6\varepsilon \int_0^{+\infty} \left[ h(\varepsilon s) - h(-\varepsilon s) \right] G(q(s)) ds
\end{aligned} \right.
\end{equation*}
Equation \eqref{Eq:ProofVolumeExpression} is verified if we manage to show that the second term of the right hand
side is $\mathcal{O}(\varepsilon^2)$.
Using the regularity of $\Omega_\varepsilon(t)$, we have:
\begin{equation*}
\forall s\in \left]0, |\log(\varepsilon)|\right[ h(\varepsilon s) - h (-\varepsilon s) = 2\varepsilon h'(0) + \mathcal{O}(s^2\varepsilon^2)
\end{equation*}
As $q(s)= \frac{1-\tanh(s)}{2}$ and $G$ is an increasing polynomial function,
the moments $\int_0^{+\infty} s^n q(s) ds$ are finite. Then,
\begin{equation*}
\left. \begin{aligned}
\left\vert \int_0^{|\log(\varepsilon)|} \left( h(\varepsilon s) - h(-\varepsilon s) \right) G(q(s)) ds \right\vert & \leq \left\vert \int_0^{|\log(\varepsilon)|} \left( 2s \varepsilon h'(0) + C s^2 \varepsilon^2\right) G(q(s)) ds \right\vert \\
& = \mathcal{O}(\varepsilon) \\
\end{aligned} \right.
\end{equation*}
On the other hand, we know that $h(s) \sim_{s\rightarrow +\infty} s^{d-1}$:
\begin{equation*}
\int_{|\log(\varepsilon)|}^{+\infty} h(\varepsilon s) G(q(s)) ds \leq C \varepsilon^{d-1} \int_{|\log(\varepsilon)|}^{+\infty} s^{d-1} G(q(s)) ds = \mathcal{O}(\varepsilon^{d-1})
\end{equation*}
and that $h$ is bounded in $\mathbb{R}_-^*$:
\begin{equation*}
\int_{|\log(\varepsilon)|}^{+\infty} h(-\varepsilon s) G(q(s)) ds \leq C \int_{|\log(\varepsilon)|}^{+\infty} G(q(s)) ds = \mathcal{O}(\varepsilon)
\end{equation*}
Globally, we conclude that:
\begin{equation*}
 6\varepsilon \int_0^{+\infty} \left[ h(\varepsilon s) - h(-\varepsilon s) \right] G(q(s)) ds = \mathcal{O}(\varepsilon^2)
\end{equation*}
and that \eqref{Eq:ProofVolumeExpression} is true and the property \eqref{Eq:VolumeNMN} is established under the 
condition that \eqref{Eq:SolutionNMN} is verified. This is the object of the next part of this section.
\subsection{Formal asymptotics toolbox}
Before the actual computations, we first recall the tools necessary to derive our formal asymptotic derivation,
following the notations of \cite{alfaro2013convergence,chen2011mass,MR3383330} and the results in differential geometry of~\cite{AmbrosioGeom}. The principle is to study separately 
the behavior of the solution near the interface and far from it. We will do the derivations in dimension 2 for the
sake of simplicity of the notations and readability, but the principle is identical in higher dimension.
 
To derive the method we require that the interface $\Gamma(t,\varepsilon)$ remains smooth enough and that there exists a
neighbourhood $\mathcal{N}=\mathcal{N}_{\delta}(\Gamma(t,\varepsilon))= \lbrace x\in\Omega / |(x,t) < 3\delta \rbrace$ in which
the signed distance function $d$ is well-defined. $\mathcal{N}$ is called the \emph{inner region} near the interface and its 
complementary the \emph{outer region}.

\medskip
\paragraph*{\bf Outer variables:}
Far from the interface, we consider the \emph{outer functions} $(u,\mu)$ depending on the standard \emph{outer variable} $x$. The system 
remains the same:
\begin{equation}
\label{Eq:OuterSystemNMN}
\left\lbrace
\begin{aligned}
& \varepsilon^2 \partial_t u = N(u) \div(M(u)\nabla(N(u)\mu)) \\
& \mu = -\varepsilon^2 \Delta u + W'(u)
\end{aligned}
\right.
\end{equation}

\medskip
\paragraph*{\bf Inner variables:}
Inside $\mathcal{N}$, we define the \emph{inner functions} $(U,\bigmu)$ depending on the \emph{inner variables} $(z,s)$, 
where $z$ is the variable along the normal and $s$ is the variable in the direction of the arc-length parametrization $S$ of
the interface $\Gamma$:
\begin{equation*}
\left\lbrace
\begin{aligned}
& U(z,s,t):= U\left( \frac{d(x,t)}{\varepsilon}, S(x,t), t \right) = u(x,t) \\
& \bigmu(z,s,t): = \bigmu\left( \frac{d(x,t)}{\varepsilon}, S(x,t), t \right) = \mu(x,t) \\
\end{aligned}
\right.
\end{equation*}
In order to express the derivatives of $U$, we first need to calculate the gradient and the laplacian of $d$ and $S$. 
The properties of $d$ are common knowledge in differential geometry, see for instance~\cite{AmbrosioGeom}:
\begin{equation*}
\left\lbrace
\begin{aligned}
& \nabla d(x,t) = n(x,t) \\
& \Delta d(x,t) = \sum_{k=1}^{d-1} \frac{\kappa_k(\pi(x))}{1+\kappa_k(\pi(x))d(x,t)} \\
&\qquad \quad \ \hspace*{2.5pt}  =\frac{H}{1+\varepsilon z H} \text{ in dimension 2}  
\end{aligned}
\right.
\end{equation*}
Let $X_0(s,t)$ be a given point of the interface, then deriving the equation connecting the variable $s$ and
the function $S$ gives:
\begin{equation*}
s = S(X_0(s,t) + \varepsilon z n(s,t),t)
\end{equation*}
with respect to $z$:
\begin{equation*}
\left. \begin{aligned}
0 & = \varepsilon n \cdot \nabla S \\
  & = \varepsilon \nabla d \cdot \nabla S \\
\end{aligned} \right.
\end{equation*}
This means that there are no cross derivative terms. We now derive the same equation with 
respect to $s$:
\begin{equation*}
\left. \begin{aligned}
1 & = \left( \partial_s X_0 + \varepsilon z H \partial_s n \right) \cdot \nabla S \\
  & = (1+\varepsilon z H) \tau \cdot \nabla S \\
\end{aligned} \right.
\end{equation*}
We know that $\nabla S$ is orthogonal to $n$, meaning it is colinear to the tangent $\tau$, then:
\begin{equation*}
\nabla S = \frac{1}{1+\varepsilon z H} \tau
\end{equation*}
Taking the divergence, we find $\Delta S$:
\begin{equation*}
\left. \begin{aligned}
\Delta S & = \div\left( \frac{1}{1+\varepsilon z H} \tau \right) \\
& = \nabla \left( \frac{1}{1+\varepsilon z H} \right) \cdot \tau + \frac{1}{1+\varepsilon z H} \div(\tau) \\
& = \frac{1}{1+\varepsilon z H} \partial_s \left( \frac{1}{1+\varepsilon z H}\right) + \frac{1}{1+\varepsilon z H} \tau \cdot  \partial_s \tau \\
& = - \frac{\varepsilon z \partial_s H}{(1+\varepsilon z H)^3}
\end{aligned} \right.
\end{equation*}
To express the connection between the derivatives of $U,\bigmu$ and $u,\mu$, we come back to the definition of the inner functions:
\begin{equation}
\label{Eq:RefEqualities}
u(x,t) = U\left(\frac{d(x,t)}{\varepsilon}, S(x,t),t) \right)
\end{equation}
Successive derivations with respect to $x$ give the following equations
\begin{equation}
\label{Eq:Formulaire}
\left\lbrace \begin{aligned}
& \nabla u = \nabla d \frac{1}{\varepsilon}\partial_z U + \nabla S \partial_s U \\
& \Delta u = \Delta d \frac{1}{\varepsilon}\partial_z U + \frac{1}{\varepsilon^2} \partial_{zz} U + \Delta S \partial_s U + |\nabla S|^2 \partial_{ss} U \\
& \div\left( M(u) \nabla (N(u)\mu)\right) = \frac{1}{\varepsilon^2}\left( \partial_{z} M \partial_z(N\bigmu)\right) + \frac{M}{\varepsilon} \Delta d \partial_z (N\bigmu)  \\
& \qquad \qquad \qquad \qquad \qquad \quad + |\nabla S|^2 \partial_s \left( M \partial_s (N\bigmu) \right) + \Delta S M \partial_s (N\bigmu)
\end{aligned}\right.
\end{equation}
The \emph{inner system} of the {\bf NMN-CH} model finally reads:
\begin{equation}
\label{Eq:InnerSystemNMN}
\left\lbrace
\begin{aligned}
& \varepsilon^2 \partial_t U +\varepsilon^2 \partial_tS \partial_s U - \varepsilon V \partial_zU = \frac{N}{\varepsilon^2} \partial_{z}\left(M \partial_z(N\bigmu)\right) + \frac{NM}{\varepsilon} \Delta d \partial_z(N\bigmu) + T_1(s) \\
& \bigmu = W'(U) - \partial_{zz} U - \varepsilon \Delta d \partial_z U - \varepsilon^2 T_2(s) \\
& T_1(s) = -\frac{\varepsilon z \partial_s H NM }{(1+\varepsilon z H)^3} \partial_s(N\bigmu) + \frac{N}{(1+\varepsilon z H)^2} \partial_{s}(M\partial_s(N\mu)) \\
& T_2(s) = -\frac{\varepsilon z \partial_s H}{(1+\varepsilon z H)^3} \partial_s U + \frac{1}{(1+\varepsilon z H)^2} \partial_{ss}U \\
& \Delta d = \frac{H}{1+\varepsilon z H} \\
\end{aligned}
\right.
\end{equation}

\medskip
\paragraph*{\bf Independence in $z$ of the normal velocity V:}
The normal velocity of the interface $V(s,t)$ is defined by:
\begin{equation*}
\left. \begin{aligned}
V(s,t) & = \partial_t X_0(s,t) \cdot n(s,t) \\
\end{aligned}\right.
\end{equation*}
In the neighbourhood $\mathcal{N}$, we have the following property (which is a direct consequence of the definition of the signed distance function):
\begin{equation*}
d(X_0(s,t) + \varepsilon z n(s,t),t) = \varepsilon z
\end{equation*}
Deriving this with respect to $t$ yields:
\begin{equation*}
V(s,t) = \partial_t X_0(s,t) \cdot \nabla d(X_0(s,t)+ \varepsilon z n(s,t),t) = - \partial_t d \left( X(z,s,t),t\right)
\end{equation*}
Thus, the function $\partial_t d(x,t)$ is independent of $z$ and we can extend the function everywhere in the neighbourhood by chosing:
\begin{equation*}
V(X_0(s,t) + \varepsilon z n,t) := -\partial_t d(X_0(s,t) +\varepsilon z n,t) = V(s,t)
\end{equation*}
This property of independence is crucial to be able to extract the velocity from integrals in $z$ in the following derivations. 

\medskip
\paragraph*{\bf Taylor expansions:}
We assume the following Taylor expansions for our functions:
\begin{equation*}
\left.
\begin{aligned}
& u(x,t) = u_0(x,t)+\varepsilon u_1(x,t)+ \varepsilon^2 u_2(x,t) + \cdots \\
& U(z,s,t) = U_0(z,s,t) + \varepsilon U_1(z,s,t) + \varepsilon^2 U_2(z,s,t) + \cdots \\
& \mu(x,t) = \mu_0(x,t) + \varepsilon \mu_1(x,t) + \varepsilon^2 \mu_2(x,t) + \cdots \\
& \bigmu(z,s,t) = \bigmu_0(z,s,t) + \varepsilon\bigmu_1(z,s,t) + \varepsilon^2\bigmu_2(z,s,t) + \cdots \\
\end{aligned}
\right.
\end{equation*}
We can then compose these expansions with a regular function $F$:
\begin{equation*}
\left.
\begin{aligned}
F(U) = & \ F(U_0) + \\
 	   & \ + \varepsilon F'(U_0)U_1 + \\
	   & \ + \varepsilon^2 \left[ F'(U_0)U_2+\frac{F''(U_0)}{2} U_1^2\right] + \\
	   & \ + \varepsilon^3\left[F'(U_0)U_3 + F''(U_0) U_1U_2 + \frac{F'''(U_0)}{6}U_1^3 \right] + \cdots \\
\end{aligned}
\right.
\end{equation*}
We can now investigate order by order the behavior of the system. We have to study four orders as the velocity appears
in the fourth order of the first equation of the Cahn--Hilliard system.

To simplify the notation within the asymptotics,
we adopt the following notations for $M(u)$:
\begin{equation*}
M(u) = m_0 + \varepsilon m_1 + \varepsilon^2 m_2 + \cdots,
\end{equation*}
where each term corresponds to:
\begin{equation*}
\left\lbrace \begin{aligned}
& m_0 = M(u_0)\\
& m_1 = M'(u_0)u_1 \\
& m_2 = M'(u_0)u_2 + \frac{M''(u_0)}{2}(u_1)^2\\
\end{aligned} \right.
\end{equation*}
We adopt the same convention for any generic \emph{outer function} $F(u)$ or \emph{inner function} $F(U)$:
\begin{equation*}
\left.
\begin{aligned}
& F(u) = f_0 + \varepsilon f_1 + \varepsilon^2 f_2 + \varepsilon^3 f_3 +\cdots \\
& F(U) = F_0 + \varepsilon F_1 + \varepsilon^2 F_2 + \varepsilon^3 F_3 +\cdots \end{aligned}
\right.
\end{equation*}
We can now investigate order by order the behavior of system \eqref{Eq:modelNMN}. We have to study up to
the fourth order where the leading order of the velocity will appear in the first equation of \eqref{Eq:OuterSystemNMN}.

\medskip
\paragraph*{\bf Flux matching condition between inner and outer equations:}

Instead of using the matching conditions directly between the first equations of the inner and outer systems, 
it is more convenient to perform the matching on the flux $j = M(u) \nabla (\sigma\mu+\lambda)$. It has the following Taylor expansion:
\begin{equation}
\label{Eq:DeuxMobOuterFlux}
\left.
\begin{aligned}
j & = \left[ m_0\nabla(n\mu)_0 \right] \\
& \quad + \varepsilon\left[m_1 \nabla (n\mu)_0 +m_0\nabla (n\mu)_1) \right] \\
		 & \quad + \varepsilon^2 \left[ m_2\nabla(n\mu)_0 + m_1\nabla(n\mu)_1 + m_0 \nabla (n\mu)_2 \right] \\
		 & \quad  + \mathcal{O}(\varepsilon^3)
\end{aligned}
\right.
\end{equation}
In inner coordinates, we only need to express the normal part $J_n:= J \cdot n = \frac{M(U)}{\varepsilon} \partial_z \bigmu$ because
the tangential part terms are of higher order. It expands as:
\begin{equation}
\label{Eq:DeuxMobInnerFlux}
\left.
\begin{aligned}
J_n & = \ \ \frac{1}{\varepsilon} \left[ M_0 \partial_z(N\bigmu)_0 \right]\\
 	& \quad \quad + \ \ \left[M_1 \partial_z (N\bigmu)_0 + M_0 \partial_z (N\bigmu)_1 \right]\\
	&\quad \quad + \varepsilon \left[M_2\partial_z(N\bigmu)_0 + M_1\partial_z (N\bigmu)_1+ M_0 \partial_z (N\bigmu)_2\right] \\
	& \quad \quad + \varepsilon^2 \left[M_3\partial_z (N\bigmu)_0 + M_2\partial_z (N\bigmu)_1 \right. \\
	& \qquad \qquad \qquad \left.  + M_1 \partial_z (N\bigmu)_2 +M_0\partial_z (N\bigmu)_3 \right] \\
	& \quad \quad + \mathcal{O}(\varepsilon^3)
\end{aligned}
\right.
\end{equation}

The flux matching conditions allow to match the limit as $z\rightarrow\pm\infty$ of terms of \eqref{Eq:DeuxMobOuterFlux} 
with the correspond order terms of \eqref{Eq:DeuxMobInnerFlux}. 

\subsection{Formal matched asymptotic analysis for the new {\bf NMN-CH} model}
Now that all the tools necessary are defined, we start the derivation of the proof of Proposition~\ref{Result4}. At first order, we determine the profile of the solution. At second order, we link the curvature with the leading term of $\bigmu$ and prove that the leading error term is zero. The third order is used to establish certain relations between different terms and finally we recover the velocity in the fourth order.

\medskip
\paragraph*{\bf First order:}
At order \ordres{1}{1} the outer system \eqref{Eq:OuterSystemNMN} reads: 
\begin{equation}
\label{Eq:DeuxMob1Outer}
\left\lbrace
\begin{aligned}
& 0 = N_0 \div(M_0 \nabla (N_0\mu_0)) \\
& \mu_0 = W'(u_0)
\end{aligned}
\right.
\end{equation} 
At order \ordres{\varepsilon^{-2}}{1} the inner system \eqref{Eq:InnerSystemNMN} reads:
\begin{equation*}
\left\lbrace
\begin{aligned}
& 0 = N_0 \partial_z \left( M_0 \partial_z(N_0 \bigmu_0) \right) \\
& \bigmu_0 = W'(U_0) - \partial_{zz} U_0 \\
\end{aligned}
\right.
\end{equation*}
The first equation gives that $M_0 \partial_z (N_0\bigmu_0)=B_0$ is constant in $z$. The matching conditions on the outer flux \eqref{Eq:DeuxMobOuterFlux} and the inner flux \eqref{Eq:DeuxMobInnerFlux} at order $\varepsilon^{-1}$ impose this constant to be zero. Then $N_0\bigmu_0$ is constant. The matching conditions with the outer system \eqref{Eq:DeuxMob1Outer} give that:
\begin{equation*}
\bigmu_0 = 0
\end{equation*}
Then $U_0$ satisfies the differential equation:
\begin{equation*}
\partial_{zz} U_0 - W'(U_0) = 0
\end{equation*}
The solution to this equation is the profile $q$ given by \eqref{ProfilQ}. Thus the first order results in:
\begin{equation*}
\left\lbrace
\begin{aligned}
& \bigmu_0 = 0 \\
& U_0 = q(z) := \frac{1-\tanh(\frac{z}{2})}{2} \\
\end{aligned}
\right.
\end{equation*}

\medskip
\paragraph*{\bf Second order:}
At order \ordres{\varepsilon}{\varepsilon} the outer system \eqref{Eq:OuterSystemNMN} reads:
\begin{equation}
\label{Eq:DeuxMob2Outer}
\left\lbrace
\begin{aligned}
& 0 = N_0 \div(M_0 \nabla (N_0\mu_1)) \\
& \mu_1 = W''(u_0)u_1 = u_1
\end{aligned}
\right.
\end{equation}
At order \ordres{\varepsilon^{-1}}{\varepsilon} the inner system \eqref{Eq:InnerSystemNMN} reads:
\begin{equation}
\label{Eq:DeuxMob2Inner}
\left\lbrace
\begin{aligned}
& 0 = N_0 \partial_z \left( M_0 \partial_z (N_0\bigmu_1) \right) \\
& \bigmu_1 = W''(U_0) U_1 - \partial_{zz} U_1 - H \partial_z U_0 \\
\end{aligned}
\right.
\end{equation}
The first equation of \eqref{Eq:DeuxMob2Inner} shows that $N_0\bigmu_1$ is a certain constant $B_1$. The matching conditions between the inner flux \eqref{Eq:DeuxMobInnerFlux} and the outer flux \eqref{Eq:DeuxMobOuterFlux} at order $0$ require that (by removing all null terms):
\begin{equation*}
B_1 = \lim_{z\rightarrow +\infty} M_0\partial_z (N_0\bigmu_1) = 0 
\end{equation*}
Then there exists a function $A_1$ constant in $z$ such that $N_0\bigmu_1=A_1$. The matching from inner to outer for $\mu$ yields:
\begin{equation*}
\mu_1 = \lim_{z\rightarrow \pm \infty} \bigmu_1 = \lim_{z\rightarrow \pm \infty} \frac{A_1}{N_0} = 0
\end{equation*}
From the matching conditions with the second equation of \eqref{Eq:DeuxMob2Outer} we have:
\begin{equation*}
u_1 = \mu_1 = 0
\end{equation*}
We now determine the value of $A_1$ using the second equation of \eqref{Eq:DeuxMob2Inner}. We multiply it by $\partial_z U_0$ and integrate it. We divide the equation in three terms. The left hand side term gives:
\begin{equation*}
\int \bigmu_1 \partial_z U_0 = \int N_0 \bigmu_1 \frac{\partial_z U_0}{N_0} = A_1 \int_{-\infty}^{+\infty}\frac{\partial_z U_0(z)}{N(U_0(z))}dz = A_1 c_N
\end{equation*}
The first two terms in the right hand side vanish:
\begin{equation*}
\left.
\begin{aligned}
\int W''(U_0) U_1 \partial_z U_0 - \partial_{zz} U_1 \partial_z U_0 & =\int \partial_z \left( W'(U_0)\right) U_1 - \partial_{zz} U_1 \partial_z U_0 \\
 & = -\int \underbrace{ (W'(U_0) - \partial_{zz} U_0)}_{=0} \partial_z U_1 \\
 & \qquad \qquad + \left[ W'(U_0) U_1 - \partial_z U_0 \partial_z U_1 \right]_{-\infty}^{+\infty} \\
 & = 0 \\
\end{aligned}
\right.
\end{equation*}
The fact that the functions in the bracket term vanishes at the limit $z\rightarrow \pm \infty$ comes from the matching conditions.
The second right hand side term results in the curvature:
\begin{equation*}
\int - H (\partial_z U_0)^2 dz= - H \int_{-\infty}^{+\infty} q'(z)^2 dz = - c_W H
\end{equation*}
Then:
\begin{equation}
\label{Eq:DeuxMobBigmu1}
N_0\bigmu_1 = A_1 = -\frac{c_W}{c_N} H
\end{equation}
In conclusion, we have the following properties:
\begin{equation*}
\left\lbrace
\begin{aligned}
& \bigmu_1 = - \frac{c_W}{c_N} \frac{H}{N(q)} \\
& \partial_{zz} U_1 - W''(U_0) U_1 = H q'(z) - \bigmu_1\\
& \mu_1 = u_1 = 0 \\
\end{aligned}
\right.
\end{equation*}
Reminding that $N(z) = -\frac{1}{q'(z)}$, the equation verified by $U_1$ is:
\begin{equation*}
\partial_{zz} U_1 - W''(U_0) U_1 = 0
\end{equation*}
To solve this equation, we use the following Lemma, which is now rather standard, see for example \cite{alikakos1994convergence,alfaro2013convergence}:
\begin{lem}\label{Lemma1}
Let A(z) be a bounded function on $-\infty<z<\infty$. Then the problem:
\begin{equation*}
\left\lbrace \begin{aligned}
& \partial_{zz} \psi - W''(q(z)) \psi = A(z) \\
& \psi\left(0\right) = 0, \quad \psi\in\mathrm{L}^{\infty}(\mathbb{R}) \\
\end{aligned} \right.
\end{equation*}
has a solution if and only if:
\begin{equation}
\label{Eq:LemmaCondition}
\int_{-\infty}^{+\infty} A(z) q'(z) dz = 0
\end{equation}
Moreover the solution, if it exists, is unique, satisfies:
\begin{equation}
\forall z \in\mathbb{R},  \ |\psi(z)| \leq C \left\lVert A \right\rVert_{\mathrm{L}^\infty}
\end{equation}
and is given by the formula:
\begin{equation}
\label{Eq:ExplicitLemma}
\psi(z) = q'(z) \int_{0}^z \left( \frac{1}{(q'(s))^2} \int_{-\infty}^s A(\sigma) q'(\sigma) d\sigma \right) ds
\end{equation}
\end{lem}

\medskip
\paragraph*{\bf Sketch of the proof:} Multiplying the equation by $q'$ and integrating by parts, we see that condition \eqref{Eq:LemmaCondition} is necessary. Reciprocally, if the condition \eqref{Eq:LemmaCondition} is verified, we can perform the method of variation of constants to find the solution explicitely \eqref{Eq:ExplicitLemma}.\myqed

Using Lemma~\ref{Lemma1} with $A=0$, we have that $U_1=0$. Therefore the leading error term in $U$ is of magnitude $\varepsilon^2$ 
and \eqref{Eq:SolutionNMN} of Result \eqref{Result4} is verified.

\medskip
\paragraph*{\bf Third order:}
At order $\ordres{\varepsilon^2}{\varepsilon^2}$ the outer system \eqref{Eq:OuterSystemNMN} reads:
\begin{equation}
\label{Eq:OuterOrdre2A}
\left\lbrace
\begin{aligned}
& 0 = n_0 \div\left(m_0 \nabla(n_0\mu_2+n_1\mu_1)\right) \\
& \mu_2 = W'''(u_0) \frac{(u_1)^2}{2} + W''(u_0)u_2 - \Delta u_0 = u_2 \\
\end{aligned}
\right.
\end{equation}
At order $\ordres{1}{\varepsilon^2}$ the inner system \eqref{Eq:InnerSystemNMN} reads:
\begin{equation}
\label{Eq:InnerOrdre2A}
\left\lbrace
\begin{aligned}
& 0 = N_0 \partial_z \left( M_0 \partial_z (N_0\bigmu_2 +N_1 \bigmu_1) \right) \\
& \bigmu_2 = W'''(U_0) \frac{(U_1)^2}{2} + W'(U_0) U_2 - \partial_{zz} U_2 - H \partial_z U_1 + z H^2 \partial_z U_0
\end{aligned}
\right.
\end{equation}
Similarly to previous orders, there exists a constant $B_2$ in $z$ so that:
\begin{equation*}
M_0\partial_z(N_0\bigmu_2 + N_1\bigmu_1) = B_2
\end{equation*}
The matching of the flux terms from \eqref{Eq:DeuxMobOuterFlux} and \eqref{Eq:DeuxMobInnerFlux} of order $\varepsilon$ (removing all the null terms) yields:
\begin{equation*}
B_2 = \lim_{z\rightarrow \pm \infty}M_0\partial_z(N_0\bigmu_2 + N_1\bigmu_1) = 0
\end{equation*}
Thus:
\begin{equation}
\label{Eq:DeuxMob3OrderCOnstant}
N_0 \bigmu_2 +N_1 \bigmu_1 = A_2
\end{equation}
The derivative in $z$ of this term would have appeared at the next order. Now that we know it is constant, 
we can omit it in the next paragraph.

\medskip
\paragraph*{\bf Fourth order:}
At order \ordres{\varepsilon^3}{\varepsilon^3} the outer system \eqref{Eq:OuterSystemNMN} reads:
\begin{equation}
\label{Eq:DeuxMob4outer}
\left\lbrace
\begin{aligned}
& 0 = n_0 \div\left( n_0 \nabla(n_2\mu_1 + n_1\mu_2 + n_0\mu_3) \right)\\
& \mu_3 = W''(u_0) u_3 + W'''(u_0) u_1u_2 + \frac{W''''(u_0)}{6}(u_1)^3 - \Delta u_1 
\end{aligned}
\right.
\end{equation}
At order $\ordres{\varepsilon}{\varepsilon^3}$ the inner system \eqref{Eq:InnerSystemNMN} reads:
\begin{equation}
\label{Eq:DeuxMob4inner}
\left\lbrace
\begin{aligned}
& - V_0 \partial_z U_0 = N_0 \partial_z\left( M_0 \partial_z (N_2 \bigmu_1 + N_1 \bigmu_2 + N_0 \bigmu_3) \right)+ N_0\partial_s\left( M_0 \partial_s(N_0 \bigmu_1 ) \right) \\
& \bigmu_3 = W''(U_0) U_3 + W'''(U_0) U_1U_2 + \frac{W''''(U_0)}{6}(U_1)^3 - \partial_{zz} U_1 \\
& \qquad \qquad - H^3 z^2 \partial_z U_0+zH^2 \partial_z U_1 - H\partial_z U_2 - \partial_{ss} \bigmu_1
\end{aligned}
\right.
\end{equation}
We determine the velocity $V_0$ by multiplying by $\frac{1}{N_0}$ and integrating the first equation of \eqref{Eq:DeuxMob4inner}. We divide the equality in three terms. The left hand side term isolates the velocity:
\begin{equation*}
-V_0 \int \partial_z \frac{U_0}{N(U_0)} = -c_N V_0
\end{equation*}
The first term of the right hand side is a pure derivative:
\begin{equation*}
\int \partial_z \left(M_0 \partial_z(N_2\bigmu_1 + N_1 \bigmu_2 + N_0 \bigmu_3)\right) 
\end{equation*}
Then, by using the matching conditions between the fluxes \eqref{Eq:DeuxMobInnerFlux} and \eqref{Eq:DeuxMobOuterFlux}
at order $\varepsilon^2$ (the equations \eqref{Eq:DeuxMob3OrderCOnstant} and \eqref{Eq:DeuxMobBigmu1} ensuring that 
the other inner terms are zero) and the fact that $M$ goes to $0$ faster than the terms in $N$ goes to infinite, we have:
\begin{equation*}
 \left[M_0 \partial_z(N_2\bigmu_1 + N_1 \bigmu_2 + N_0 \bigmu_3) \right]_{-\infty}^{+\infty} =0
\end{equation*}
Finally, using \eqref{Eq:DeuxMobBigmu1}, the second term of the right hand side gives the surface diffusion part:
\begin{equation*}
\int M_0 \partial_{ss}(N_0\bigmu_1) = - \frac{c_M c_W}{c_N}\partial_{ss} H 
\end{equation*}
In conclusion, we obtain the desired motion \eqref{Eq:VelocityNMN}:
\begin{equation*}
V_0 = \frac{c_M c_W}{(c_N)^2} \partial_{ss} H
\end{equation*}
This concludes the proof of Proposition~\ref{Result4}.\myqed

\section{Numerics: discretization and experiments}
\label{Sec:Numerical}
In this section, we propose a generic numerical scheme  to solve the three different Cahn--Hilliard models:
\begin{itemize}
 \item The classical Cahn--Hilliard equation ({\bf C-CH})
 $$
\begin{cases}
 \partial_t u &=  \Delta \mu \\
 \mu &= \frac{1}{\epsilon^2} W'(u) - \Delta u,
\end{cases}
$$
where $W(s) = \frac{1}{2} s^2 (1-s)^2$.
\item The Cahn--Hilliard  model  with classical mobility ({\bf M-CH})
$$ 
\begin{cases}
 \partial_t u &=  \div( M(u) \nabla \mu) \\
 \mu &= \frac{1}{\epsilon^2} W'(u) - \Delta u.
\end{cases}
$$
where the mobility is defined as $M(u) = c_N^2 2 W(u)$. Here, the constant $c_N = 6$ is added
to get the same limit law as using our new Cahn--Hilliard model.
\item New second order variational Cahn--Hilliard equation: ({\bf NMN-CH})

$$
\begin{cases}
 \partial_t u &=  N(u) \div( M(u) \nabla (N(u) \mu)) \\
 \mu &= \frac{1}{\epsilon^2} W'(u) - \Delta u
\end{cases}
$$
where the mobility is defined as $M(u) = W(u) + \gamma \epsilon^2$ and $N(u) = \sqrt{M(u)}$. 
Here $\gamma>0$ is a smoothing parameter and we take $\gamma = 1$ for all numerical  experiments presented below,
\end{itemize}
 ~\\
 
Our numerical algorithm is constructed as a semi-implicit Fourier spectral method in the spirit of
\cite{Chen_fourier,Bretin_brassel,bretin_droplet,bretin_largephases,BRETIN2018324}, see \cite{DU2020425} for a recent review of numerical methods for the phase field approximation of various geometric flows.

All  schemes proposed here are based on a convex splitting of the Cahn--Hilliard energy, which was first proposed by 
Eyre \cite{MR1676409} and became popular as a simple, efficient, and stable  scheme to approximate
various evolution problems with a gradient flow structure~\cite{MR2418360,MR2519603,MR2799512,MR3100769,MR3564350,MR3682074}.
More recently,  a first- and second-order splitting scheme was proposed in \cite{MR3874087,salvalaglio2019doubly,Doubly_anisotropic}
to address the case of the Cahn--Hilliard equation with mobility. However, these approaches are based 
on the finite element method and are not compatible with a Fourier spectral discretization.

In this paper, we therefore propose to generalize the idea of convex splitting using
an additionally convex splitting of the variational metric associated to the mobility. 
The advantage is to make it a very simple and efficient scheme, even in the case of highly contrasted  and degenerate mobilities.
As an illustration, we present above a numerical implementation of our scheme in {\bf Matlab} that requires less than 40 lines. \\
In this section, we then give some details about these schemes and propose a numerical comparison of phase field models 
in space dimensions $2$ and $3$.

\subsection{Spatial discretization: a Fourier-spectral approach}

All equations are solved on a square-box $Q = [0,L_1]\times \cdots \times [0,L_d]$ with periodic boundary conditions.
We recall that the Fourier $\boldsymbol K$-approximation of a function $u$ defined in a box 
$Q = [0,L_1]\times \cdots \times [0,L_d]$ is given by
$$u^{\boldsymbol K}(x) = \sum_{{\boldsymbol k}\in K_N  } c_{\boldsymbol k} e^{2i\pi{\boldsymbol \xi}_k\cdot x},$$
where  $K_N =  [ -\frac{N_1}{2},\frac{N_1}{2}-1 ]\times [ -\frac{N_2}{2},\frac{N_2}{2}-1] \cdots \times   [ -\frac{N_d}{2},\frac{N_d}{2}-1] $,   ${\boldsymbol k} = (k_1,\dots,k_d)$ and ${\boldsymbol \xi_k} = (k_1/L_1,\dots,k_d/L_d)$. In this formula, the $c_{\boldsymbol k}$'s denote the $K^d$ first discrete Fourier coefficients of $u$. 
The inverse discrete Fourier transform leads to 
$u^{K}_{\boldsymbol k} =   \textrm{IFFT}[c_{\boldsymbol k}]$ 
where $u^{K}_{\boldsymbol k}$ denotes the value of $u$ at the points 
$x_{\boldsymbol k} = (k_1 h_1, \cdots, k_d h_d)$ and where $h_{\alpha} = L_{\alpha}/N_{\alpha}$ for $\alpha\in\{1,\cdots,d\}$. Conversely,
$c_{\boldsymbol k}$ can be computed as the discrete Fourier transform of $u^K_{\boldsymbol k},$ {\em i.e.}, $c_{\boldsymbol k} =
\textrm{FFT}[u^K_{\boldsymbol k}].$

\subsection{Time discretization}

Given a time discretization parameter $\delta_t > 0$, we construct a sequence $(u^n)_{n \geq 0}$ 
of approximations of ${u}$ at times $n \delta_t$.

\subsubsection{An IMEX scheme for the {\bf C-CH} model}

We propose now to use a simple scheme to discretize the classical Cahn--Hilliard equation
$$
\begin{cases}
 \partial_t u &=  \Delta \mu \\
 \mu &= \nabla_u E(u) = \frac{1}{\epsilon^2} W'(u) - \Delta u,
\end{cases}
$$
where the Cahn--Hilliard energy reads as 
$$ E(u) =  \int_Q \varepsilon \frac{|\nabla u|^2}{2} +  \frac{1}{\varepsilon} W(u) dx.$$

\medskip
\paragraph*{\bf A semi-implicit scheme based on a convex-concave splitting of $E$:}

Following the idea of \cite{MR1676409}, we propose to split the energy $E$ as the sum of a convex energy and a concave energy 
$$E(u) = E_c(u) + E_e(u),$$
with, respectively, an implicit and an explicit integration of the convex and concave parts:
$$
\begin{cases}
 (u^{n+1} - u^{n})/\delta_t &=  \Delta \mu^{n+1} \\
 \mu^{n+1} &=  \nabla_u  E_c(u^{n+1}) + \nabla_u  E_e(u^{n})   
\end{cases}
$$
Notice that this scheme can also be interpreted as an implicit discretization of the semi linearized  PDE 
$$
\begin{cases}
 \partial_t u &=  \Delta \mu \\
 \mu &= \nabla_u \overline{E}_{u^{n}}(u) = \nabla_u  E_c(u) + \nabla_u  E_e(u^{n})   
\end{cases},
$$
 where the new associated energy $\overline{E}_{u^{n}}$ reads as
$$\overline{E}_{u^{n}}(u) = E_c(u) + E_e(u^n)+ \langle \nabla_u E_e(u^{n}), (u - u^n) \rangle.$$ 
This continuous  point of view shows that $\overline{E}_{u^{n}}(u)$ is clearly decreasing along the flow  
$$\frac{d}{dt}\left(\overline{E}_{u^{n}}(u) \right)=   \langle  \overline{E}_{u^{n}}(u) , u_t \rangle = - \|  \nabla \overline{E}_{u^{n}}(u) \|^2 < 0.$$
and then
$$ \overline{E}_{u^{n}}(u^{n+1}) \leq   \overline{E}_{u^{n}}(u^{n}) = E(u^{n}),$$
Finally, the assumption on the concavity of $E_e$ implies that  $ E(u) \leq \overline{E}_{u^{n}}(u)$ and gives
the decreasing of $E$, 
$$ E(u^{n+1}) \leq    E(u^{n}).$$
without requiring any assumption on the time step $\delta_t$. 

\medskip
\paragraph*{\bf Application in the case of the Cahn--Hilliard energy:}

In the case of the Cahn--Hilliard equation using the smooth double well potential 
$W(s) = \frac{1}{2}s^2 (1-s)^2$, a standard splitting choice is
$$ E_c(u) =  \frac{1}{2}\int_Q \varepsilon |\nabla u|^2  + \frac{\alpha}{\varepsilon^2} u^2  dx \text{ and }  E_e(u) =  \int_{Q} \frac{1}{\varepsilon} (W(u) - \alpha \frac{u^2}{2})dx.$$
Notice that $E_e$ is clearly concave as soon as  $\alpha \geq \max_{s \in [0,1]} \left| W''(s) \right|$.
In particular, this approach leads to the semi-implicit scheme 
$$
\begin{cases}
 (u^{n+1} - u^{n})/\delta_t &=  \Delta \mu^{n+1} \\
 \mu^{n+1} &=  \left( - \Delta u^{n+1} + \frac{\alpha}{\epsilon^2} u^{n+1} \right) + \left( \frac{1}{\epsilon^2} (W'(u^{n}) - \alpha u^{n}) \right),   
\end{cases}
$$
which also reads as
$$ \begin{pmatrix}
    I_d &  - \delta_t \Delta \\
    \Delta - \alpha/\epsilon^2 & I_d
   \end{pmatrix}  
   \begin{pmatrix}
    u^{n+1} \\
    \mu^{n+1}
   \end{pmatrix}
   =
    \begin{pmatrix}
    u^{n} \\
     \frac{1}{\epsilon^2} (W'(u^{n}) - \alpha u^{n}).
   \end{pmatrix}
  $$
 Finally, the couple $(u^{n+1},\mu^{n+1})$  can be expressed as
  $$u^{n+1} = L \left[ u^{n} + \frac{\delta_t}{\epsilon^2} \Delta \left(  W'(u^{n}) - \alpha u^{n} \right) \right] \text{ and } \mu^{n+1}= L\left[ \frac{1}{\epsilon^2} (W'(u^{n} )  -  \Delta u^{n}  \right].$$
 Here, the operator $L = \left( I_d +  \delta_t \Delta ( \Delta - \alpha/\epsilon^2 I_d ) \right)^{-1}$ can be easily computed  in Fourier space
 like a symbol operator associated to 
 $$\hat{L}(\xi) = 1/(1 + \delta_t 4 \pi^2 |\xi|^2 ( 4 \pi^2 |\xi|^2 + \alpha/\epsilon^2)).$$

\subsubsection{A numerical scheme for the M-CH model }

We now consider the case of the {\bf M-CH} model, which reads
$$
\begin{cases}
 \partial_t u &=  \div \left( M(u) \nabla \mu \right) \\
 \mu &= \frac{1}{\epsilon^2} W'(u) - \Delta u.
\end{cases}
$$
As previously, it should be interesting to consider the following scheme    
$$
\begin{cases}
 (u^{n+1} - u^{n})/\delta_t &=  \div \left( M(u^n) \nabla \mu^{n+1} \right) \\
 \mu^{n+1} &=  \nabla_u  E_c(u^{n+1}) + \nabla_u  E_e(u^{n}).   
\end{cases}
$$
It can also be interpreted  as an implicit discretization of the modified Cahn--Hilliard system
$$
\begin{cases}
 \partial_t u &=  \div \left( M(u^{n}) \nabla \mu \right) \\
 \mu &= \nabla_u \overline{E}_{u^n},
\end{cases}
$$
which shows that  $E(u^{n+1}) \leq E(u^{n})$ as 
$$\frac{d}{dt}\left(\overline{E}_{u^{n}}(u) \right)=   \langle  \overline{E}_{u^{n}}(u) , u_t \rangle = - \| \sqrt{ M(u^{n})} \nabla \overline{E}_{u^{n}}(u) \|^2 < 0.$$

However, such an approach requires the computation of the new operator $L_{M,u^n}$ defined by  
$L_{M,u^n} = \left( I_d +  \delta_t \div ( M(u^n) \nabla  ( \Delta + \alpha/\epsilon^2 ) \right)^{-1},$
which cannot be made in Fourier space. Notice also that this approach has been recently proposed 
in \cite{MR3874087,salvalaglio2019doubly,Doubly_anisotropic} where  the resolution of $(u^{n+1}, \mu^{n+1})$
has been made using finite elements.

\medskip
\paragraph*{\bf Imex approach on the variational mobility term:}

We then propose another approach in this paper keeping in mind  the variational property of mobility:
$$ \begin{cases}
    \partial_t u &= - \nabla_\mu J_{u}(\mu) \\
     \mu &= \nabla_u \overline{E}_{u^n}
   \end{cases}
 $$
where  
$$J_{u}(\mu) = \frac{1}{2}  \int_Q M(u) |\nabla \mu|^2 dx.$$
As for the energy  $E$, we then  propose to split also $J$  as the sum of a convex and a concave term  $J_{u}  = J_{u,c} + J_{u,e}$
with respectively an implicit and explicit treatment  of the convex and concave part:
$$
\begin{cases}
 (u^{n+1} - u^{n})/\delta_t &=  -\nabla_{\mu} J_{u^{n},c}(\mu^{n+1})   -\nabla_{\mu} J_{u^{n},e}(\mu^{n}),\\
 \mu^{n+1} &=  \nabla_u  E_c(u^{n+1}) + \nabla_u  E_e(u^{n}).   
\end{cases}
$$
As previously, this scheme can be interpreted as an Euler implicit discretization of  
$$
\begin{cases}
 \partial_t u &=  - \nabla_\mu \overline{ J}_{u^{n},\mu^{n}}(\mu) \\
 \mu &= \nabla_u \overline{E}_{u^n},
\end{cases}
$$
where the new mobility energy  $\overline{J}_{u^n,\mu^n}$ is given by
$$\overline{J}_{u^{n},\mu^{n}}(\mu) =  J_{u^{n},c}(\mu) + J_{u^{n},e}(\mu^{n}) + \langle \nabla_{\mu} J_{u^{n},e}(\mu^{n}), \mu - \mu^n  \rangle .$$
Then, to ensure the decrease of $t \mapsto \overline{E}_{u^n}(u(\cdot,t))$ along the flow, we require at least  the semi-implicit
metric $\overline{J}_{u^n,\mu^n}$ to be non negative. This corresponds to the concavity condition on $J_{u,e}$, meaning that we have
$$ 0 \leq J_{u^{n}}(\mu) \leq \overline{J}_{u^{n},\mu^{n}}(\mu).$$
Moreover, from the identity
 $$ \frac{d}{dt} \overline{E}_{u^n}(u) =  \langle \nabla_u \overline{E}_{u^n}, u_t \rangle = - \langle \mu ,   \nabla_\mu \overline{ J}_{u^{n},\mu^{n}}(\mu) \rangle,$$
we conclude that it is sufficient to show that
 $$  \langle \mu ,   \nabla_\mu \overline{ J}_{u^{n},\mu^{n}}(\mu) \rangle \geq 0.$$
to ensure the decrease of the energy.

\medskip
\paragraph*{\bf Application to the {\bf M-CH} model:}

Motivated by the previous section, we propose the following splitting of $J$:
$$ J_{u^n,c}(\mu) = \frac{1}{2} \int m |\nabla \mu|^2 dx \quad  \text{ and }  \quad  J_{u^n,e}(\mu) =  \frac{1}{2} \int (M(u^{n})-m) |\nabla \mu|^2 dx    $$
with $m>0$. We take $m  = \max_{s \in [0,1]} \left\{ M(s) \right\}$ in order to obtain the concavity of $J_{u^n,e}(\mu)$, and the scheme reads
$$
\begin{cases}
 (u^{n+1} - u^{n})/\delta_t &=  m \Delta \mu^{n+1} +  \div( (M(u^n) - m) \nabla \mu^n ) \\
 \mu^{n+1} &=  \left( - \Delta u^{n+1} + \frac{\alpha}{\epsilon^2} u^{n+1} \right) + \left( \frac{1}{\epsilon^2} (W'(u^{n}) - \alpha u^{n}) \right),   
\end{cases}
$$
or in a matrix form 
$$ \begin{pmatrix}
    I_d &  - \delta_t m \Delta \\
    \Delta - \alpha/\epsilon^2 & I_d
   \end{pmatrix}  
   \begin{pmatrix}
    u^{n+1} \\
    \mu^{n+1}
   \end{pmatrix}
   =
    \begin{pmatrix}
    u^{n} +  \delta_t \div( (M(u^n) - m) \nabla \mu^n ) \ \\
     \frac{1}{\epsilon^2} (W'(u^{n}) - \alpha u^{n}).
   \end{pmatrix}
   = 
    \begin{pmatrix}
    B^1_{u^{n},\mu^n}\\
    B^2_{u^{n},\mu^n}
   \end{pmatrix}
  $$
  Finally, the couple  $(u^{n+1},\mu^{n+1})$ can be expressed  as 
  
$$u^{n+1} = L_M \left[   B^1_{u^{n},\mu^n} +  \delta_t m \Delta B^2_{u^{n},\mu^n} \right]$$
and
$$\mu^{n+1}= L_M\left[ ( - \Delta  B^1_{u^{n},\mu^n} + \alpha/\epsilon^2 B^1_{u^{n},\mu^n} ) +   B^2_{u^{n},\mu^n} \right],$$
where the operator $L_M$ is now given by $L_M = \left( I_d +  \delta_t m \Delta ( \Delta - \alpha/\epsilon^2 I_d) \right)^{-1}$,
which can be computed efficiently in Fourier space.\\

\subsubsection{Case of the NMN-CH model}

We now turn to the {\bf NMN-CH} model:
$$
\begin{cases}
 \partial_t u &=  N(u) \div( M(u) \nabla (N(u) \mu)) \\
 \mu &= \frac{1}{\epsilon^2} W'(u) - \Delta u,
\end{cases}
$$
where $N(u) = \frac{1}{\sqrt{M(u)}}$ and $M(u) = W(u) + \gamma \epsilon^2$.

In a similar manner to the other models, we study the model rewritten in a variational form
$$ \begin{cases}
    \partial_t u &= - \nabla_\mu J_{u}(\mu) \\
     \mu &= \nabla_u \overline{E}_{u^n}
   \end{cases}
 $$
 with
 $$ J_{u}(\mu) =  \frac{1}{2} \int_Q  M(u) \left| \nabla ( N(u) \mu) \right|^2 dx.$$
$J_u$ can be split in three parts:
 $$ J_{u}(\mu) = \frac{1}{2} \int_Q |\nabla \mu|^2 dx + \int_Q G(u) \cdot \nabla \mu \mu dx + \frac{1}{2}  \int_Q |G(u)|^2 \mu^2  dx,$$
with $$G(u) = -\frac{1}{2} \nabla(log(M(u)))$$
 as $N(u) = \frac{1}{\sqrt{M(u)}}$ and $\sqrt{M(u)} \nabla (N(u)) = - \frac{1}{2} \frac{\nabla M(u)}{ M(u) } = - \nabla ( log(M(u)))$. \\
This suggests that we could use the following splitting of $J_{u}(\mu) = J_{u,c}(\mu) +  J_{u,e}(\mu)$  with
 $$ J_{u,c}(\mu) =   \frac{1}{2} \int_Q m |\nabla \mu|^2 dx + \frac{1}{2}  \int_Q \beta \mu^2  dx$$
 and
 $$ J_{u,e}(\mu) =  \int_Q G(u) \cdot \nabla \mu \mu dx + \frac{1}{2}  \int_Q (|G(u)|^2 - \beta) \mu^2  dx + \frac{1}{2} \int_Q (1 -m) |\nabla \mu|^2 dx ,$$
 with $\beta>0$ and $m>0$. Moreover, as soon as $G(u)$ is bounded is $H_1(Q)$,  a sufficiently large choice fo
 $m$ and $\beta$ should  ensure the concavity of $J_{u,e}(\mu)$. In practice, we take $m=1$ and 
 $\beta = 1/\epsilon^2$ for our numerical experiments and these values
did not show any sign of instability regardless of the choice of the time step $\delta_t$. 
In particular, this leads to the following system
 
 $$
\begin{cases}
 (u^{n+1} - u^{n})/\delta_t &=  m \Delta \mu^{n+1} - \beta \mu^{n+1}  +  H(u^n,\mu^n) \\
 \mu^{n+1} &=  \left( - \Delta u^{n+1} + \frac{\alpha}{\epsilon^2} u^{n+1} \right) + \left( \frac{1}{\epsilon^2} (W'(u^{n}) - \alpha u^{n}) \right),   
\end{cases}
$$
where 
$$   H(u^n,\mu^n) =  N(u^n)\div( (M(u^n) \nabla (N(u^{n}) \mu^n) ) - m \Delta \mu^{n} + \beta \mu^{n} $$
The couple $(u^{n+1},\mu^{n+1})$ is then solution of the system
$$ \begin{pmatrix}
    I_d &  - \delta_t (m \Delta - \beta I_d)  \\
    \Delta - \alpha/\epsilon^2 & I_d
   \end{pmatrix}  
   \begin{pmatrix}
    u^{n+1} \\
    \mu^{n+1}
   \end{pmatrix}
   =
    \begin{pmatrix}
    u^{n} +  \delta_t H(u^n,\mu^n)\ \\
     \frac{1}{\epsilon^2} (W'(u^{n}) - \alpha u^{n}).
   \end{pmatrix}
   = 
    \begin{pmatrix}
    B^1_{u^{n},\mu^n}\\
    B^2_{u^{n},\mu^n}
   \end{pmatrix}
  $$
satisfying
$$u^{n+1} = L_{NMN} \left[   B^1_{u^{n},\mu^n} +  \delta_t ( m \Delta B^2_{u^{n},\mu^n} - \beta B^2_{u^{n},\mu^n} ) \right]$$
and
$$\mu^{n+1}= L_{NMN}\left[ ( - \Delta  B^1_{u^{n},\mu^n} + \alpha/\epsilon^2 B^1_{u^{n},\mu^n} ) +   B^2_{u^{n},\mu^n} \right].$$
Here the operator $L_{NMN}$ is given by
$L_{NMN} = \left( I_d +  \delta_t (m \Delta - \beta I_d)( \Delta - \alpha/\epsilon^2 I_d) \right)^{-1}$,
which can be still computed efficiently in Fourier space.

\subsection{{\bf Matlab} code }
We present in Figure \eqref{fig:matlab_code} an example of {\bf Matlab} script with less than $40$ lines 
which implements the scheme approximating the solutions of the {\bf NMN-CH} model.
In particular : 
\begin{itemize}
 \item We consider here a computation box $Q = [-1/2,1/2]^2$ discretized  with $N = 2^9$ nodes in each direction. 
 The initial condition of $u$ is a uniform noise and the numerical parameters are given 
 by $\epsilon = 2/N$, $\delta_t = 4 \epsilon^2$,  
 $\alpha = 2$, $\beta = 2/\epsilon^2$ and $m=1$.   
 \item Line $14$ corresponds to the definition of the Fourier-symbol associated with operator $L_{NMN}$.
   The application of $L_{NMN}$ can then be performed using 
  a simple multiplication  in Fourier space with the array $M_{LNMN}$.
 \item The computation of $N(u) \div(M(u) \nabla (N(u) \mu ) ) $ is made on  line $28$ and is based on the following equality  
\begin{eqnarray*}
 N(u) \div(M(u) \nabla (N(u)\mu)) &=&  \sqrt{M(u)} \Delta N(u)\mu + N(u) \nabla (M(u)) \cdot \nabla (N(u)\mu) \\
                                  &=&  \sqrt{M(u)} \Delta N(u)\mu +   2 \nabla \left[\sqrt{M(u)}\right] \cdot \nabla (N(u)\mu), 
\end{eqnarray*}
 as $N = 1/\sqrt{M(u)}$.
 \item Each computation of gradient and divergence operator are made in Fourier space.
 For instance the gradient of $\sqrt{M(u)}$ is computed on line $23$.
 \item Figure \eqref{fig_Rand_NMN} shows the phase field function $u^{n}$ computed  
 at different times $t^{n}$ by using this script.   
 \end{itemize}
We believe that this implementation shows the simplicity, efficiency and stability of our numerical scheme. 

\begin{figure}[htbp]
\centering
\begin{lstlisting}
clear all;
%%%%%%%%%%%%%%%%% Numerical parameters %%%%%%%%%%%%%%%%%%%%%%%%
N = 2^9; epsilon =1/N; dt =epsilon^4; T =1;
%%%%%%%%%%%%%%%%% Double well potential, mobilities %%%%%%%%%%%
W = @(U) 1/2*(U.*(U-1)).^2;
W_prim = @(U) (U.*(U-1).*(2*U-1));
MobM = @(U) 1/2*((((U).*(1-U)).^2+epsilon^2) );
MobN = @(U) 1./sqrt(MobM(U) );

%%%%%%%%%%%%%%% Fourier operators %%%%%%%%%
k =  [0:N/2,-N/2+1:-1]; [K1,K2] = meshgrid(k,k);
Delta = -4*pi^2*((K1.^2 + (K2).^2));
alpha = 2; beta = 1/epsilon^2; m = 1;
M_LNMN = 1./(1 + dt*(m*Delta - beta) .*(Delta - alpha/epsilon^2));

%%%%%%%%%%%%%%%%%% Initial condition %%%%%%%%%%%
U = rand(N,N); U_fourier = fft2(U);
Mu = zeros(N,N); Mu_fourier = zeros(N,N);
%%%%%%%%%%%%%%%%%% Scheme loop %%%%%%%%%%%
for i=1:T/dt,
mobMU = MobM(U); mobNU = MobN(U);
sqrtM = sqrt(mobMU); sqrtM_fourier = fft2(sqrtM); 
nabla1_sqrtM=  real(ifft2(2*pi*1i*K1.*sqrtM_fourier )); nabla2_sqrtM=  real(ifft2(2*pi*1i*K2.*sqrtM_fourier ));

muN_fourier = fft2(Mu.*mobNU); muN = real(ifft2(muN_fourier));
nabla1_muN =  real(ifft2(2*pi*1i*K1.*muN_fourier )); nabla2_muN =  real(ifft2(2*pi*1i*K2.*muN_fourier ));
laplacien_muN =  real(ifft2(Delta.*muN_fourier ));
NdivMgradNMu = sqrtM.*laplacien_muN + 2*(nabla1_sqrtM.*nabla1_muN +nabla2_sqrtM.*nabla2_muN);

B1 = U_fourier + dt*(fft2(NdivMgradNMu) - (m*Delta-beta).*Mu_fourier);
B2 = fft2(W_prim(U)/epsilon^2 - alpha/epsilon^2*U);

U_fourier = M_LNMN.*(B1 + dt*(m*Delta-beta).*B2);
U = real(ifft2(U_fourier));
Mu_fourier = M_LNMN.*((alpha/epsilon^2 - Delta).*B1 + B2);  
Mu = real(ifft2(Mu_fourier));

end



 
\end{lstlisting}
	 \caption{Example of {\bf Matlab} implementation of the previous scheme in dimension $2$ to approximate the solutions to the {\bf NMN-CH} model.}
\label{fig:matlab_code}
\end{figure}

\begin{figure}[htbp]
\centering
	\includegraphics[width=3.5cm]{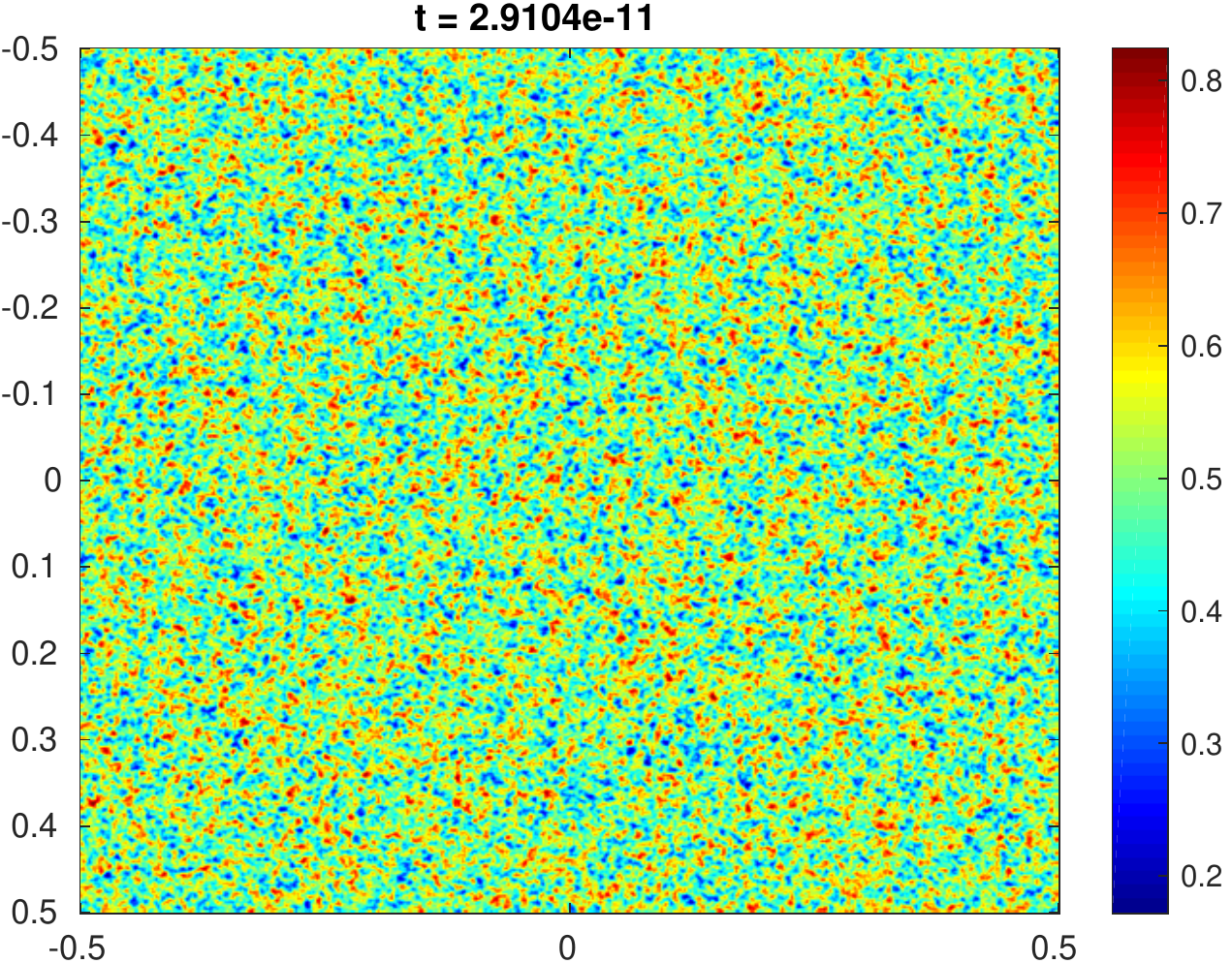}
	\includegraphics[width=3.5cm]{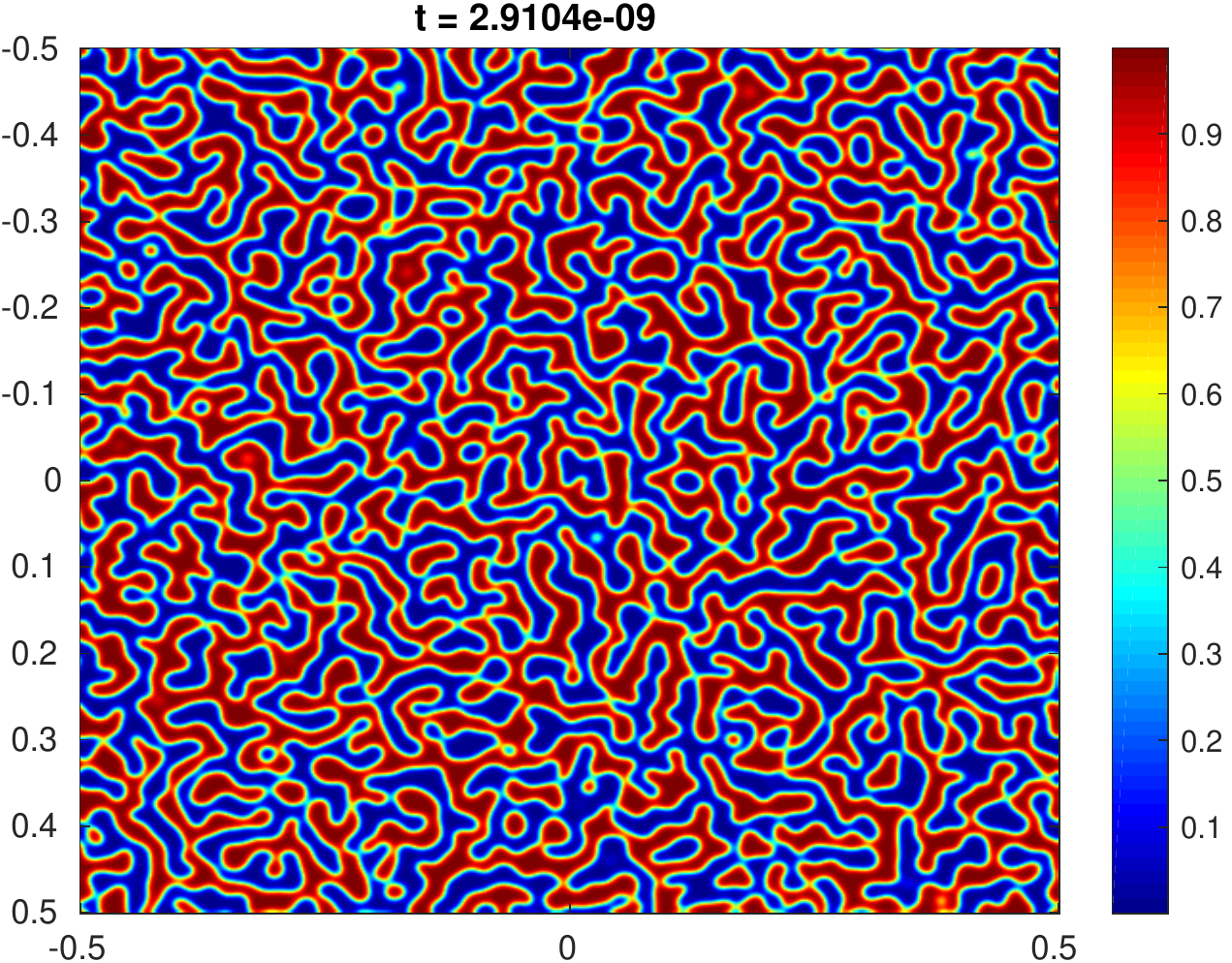}
	\includegraphics[width=3.5cm]{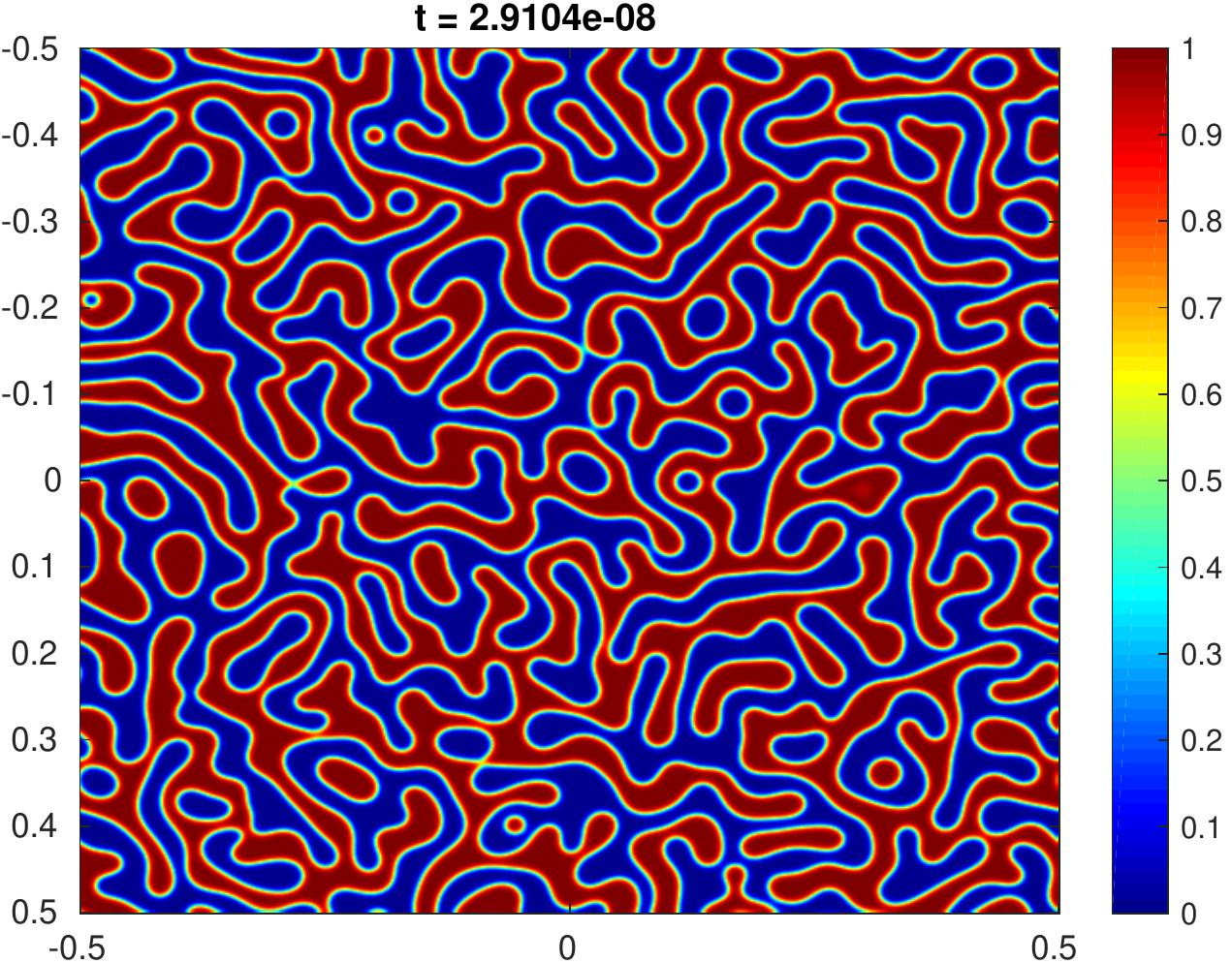}
	\includegraphics[width=3.5cm]{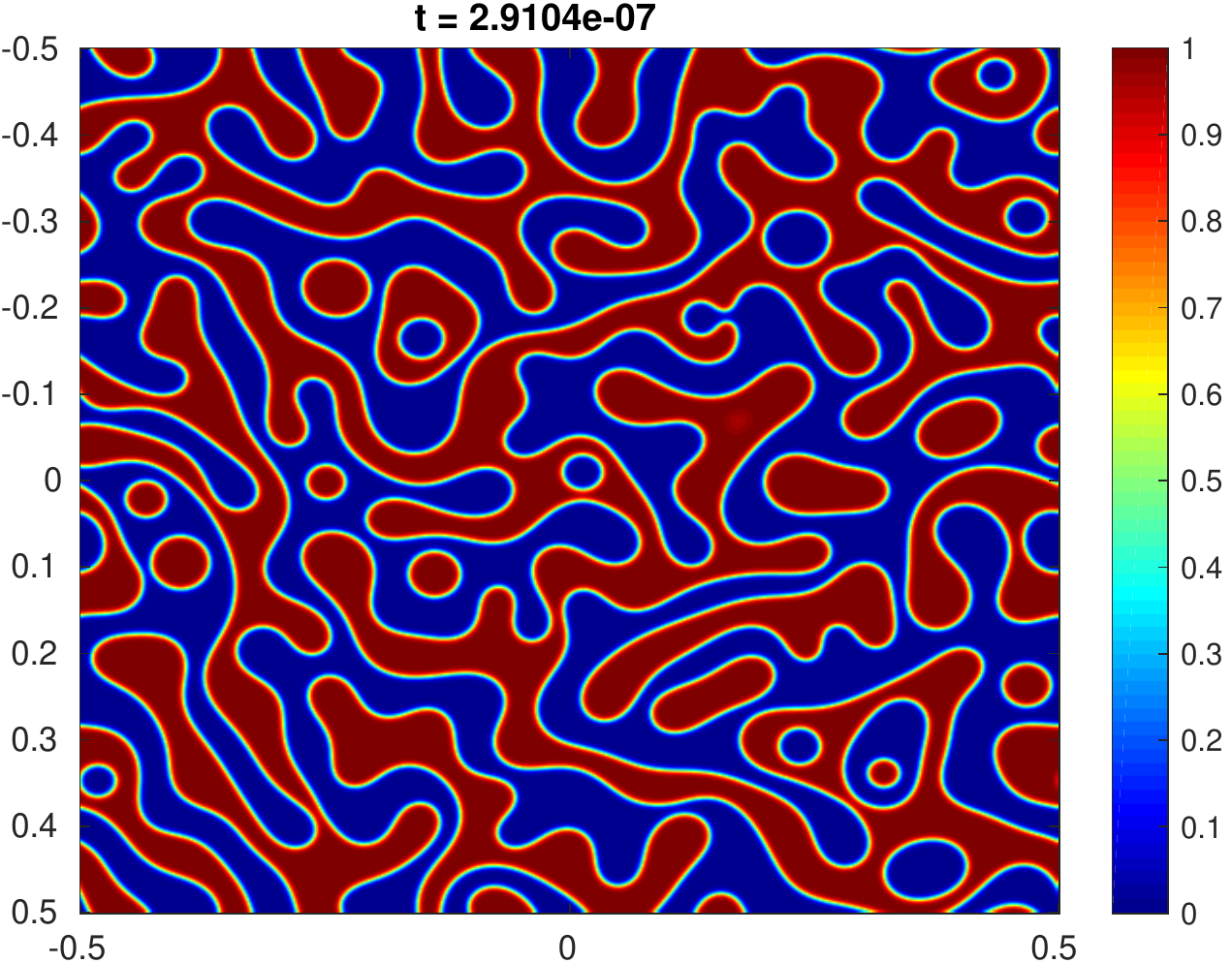}
\caption{First numerical experiment using the {\bf NMN-CH} model; the solutions $u$ are computed with the {\bf Matlab}
script of Figure~\ref{fig:matlab_code}.}
\label{fig_Rand_NMN}
\end{figure}

\subsubsection{Asymptotic expansion and flow: numerical comparison of the different models}

\begin{figure}[htbp]
\centering
	\includegraphics[width=3.5cm]{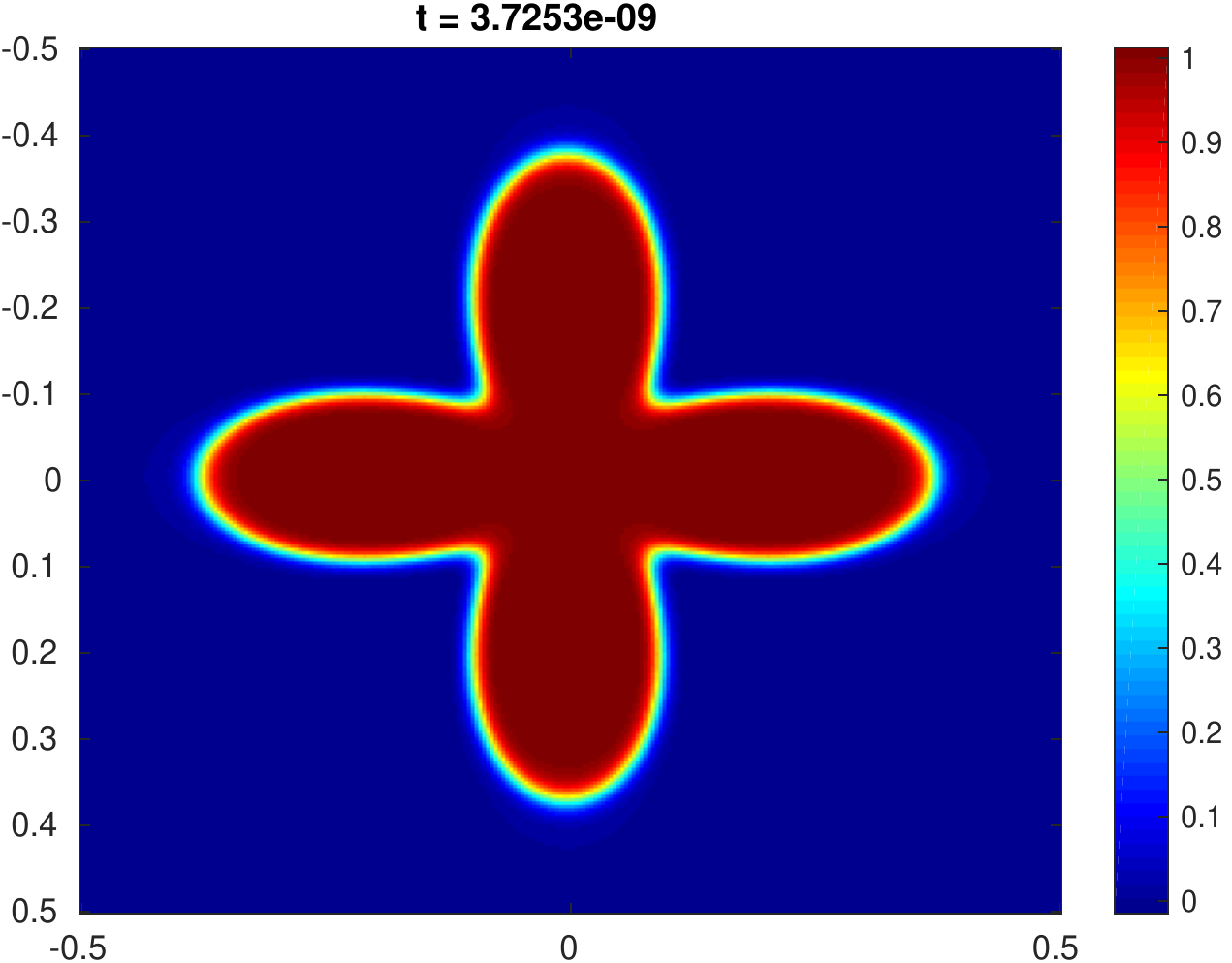}
	\includegraphics[width=3.5cm]{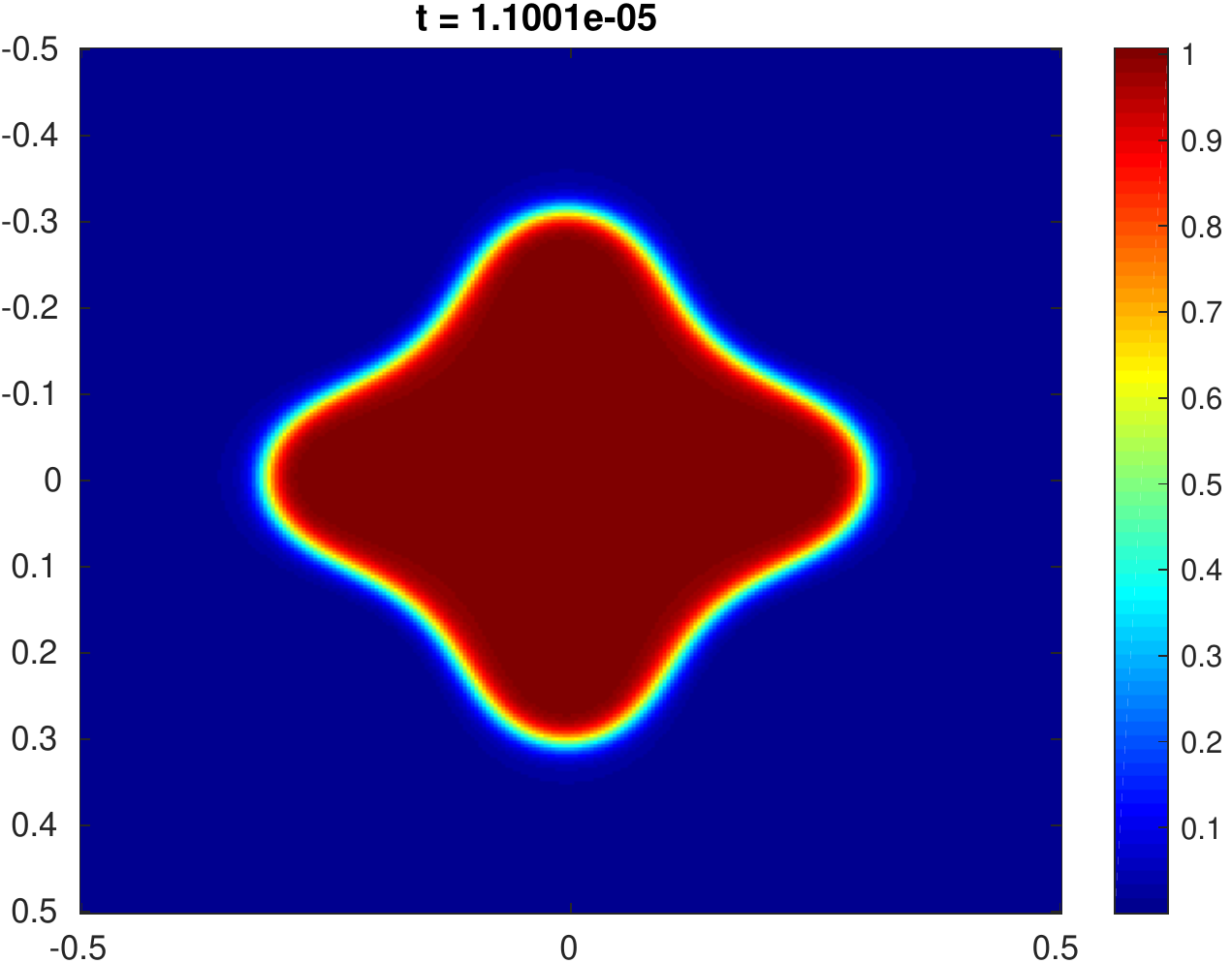}
	\includegraphics[width=3.5cm]{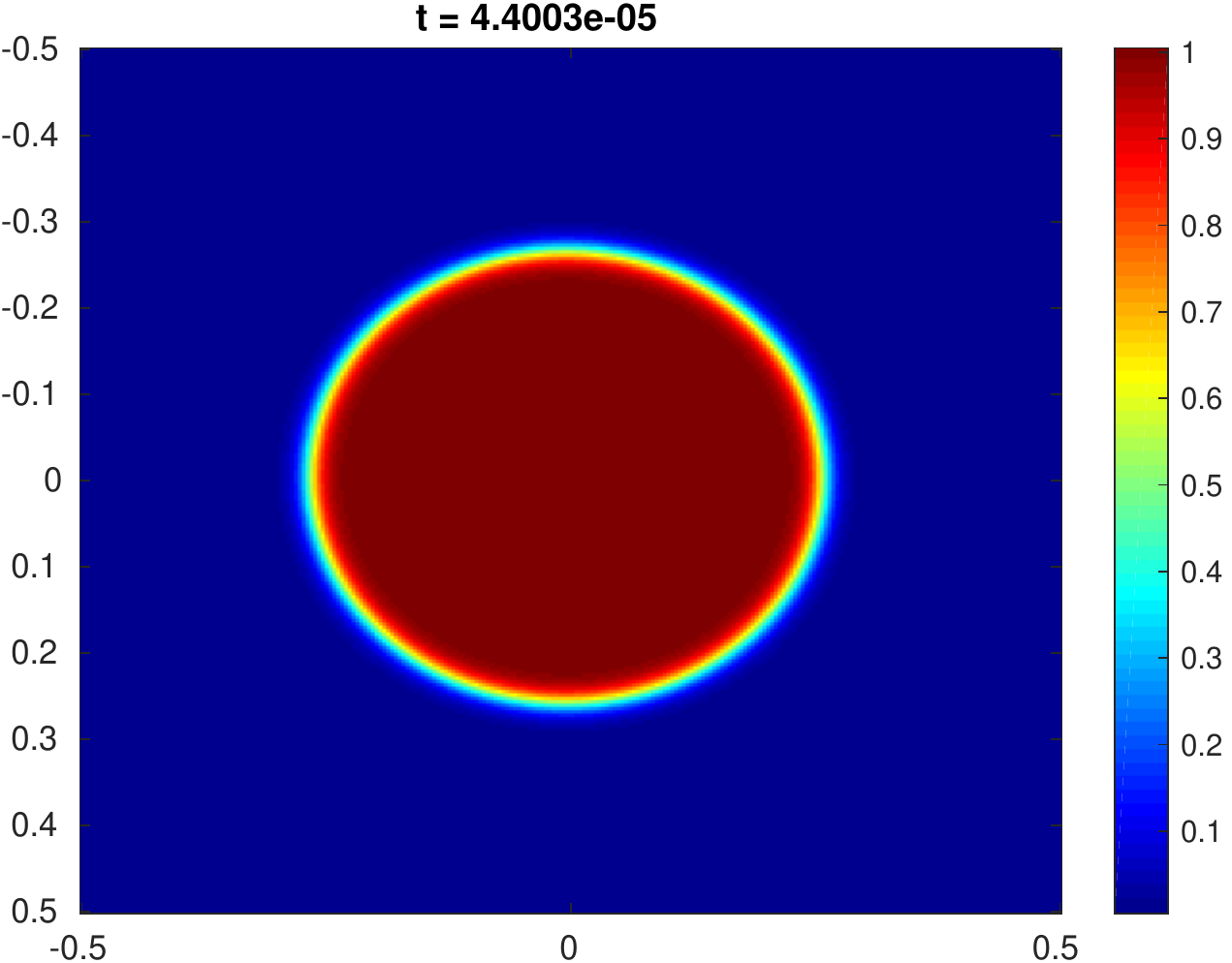}
	\includegraphics[width=3.5cm]{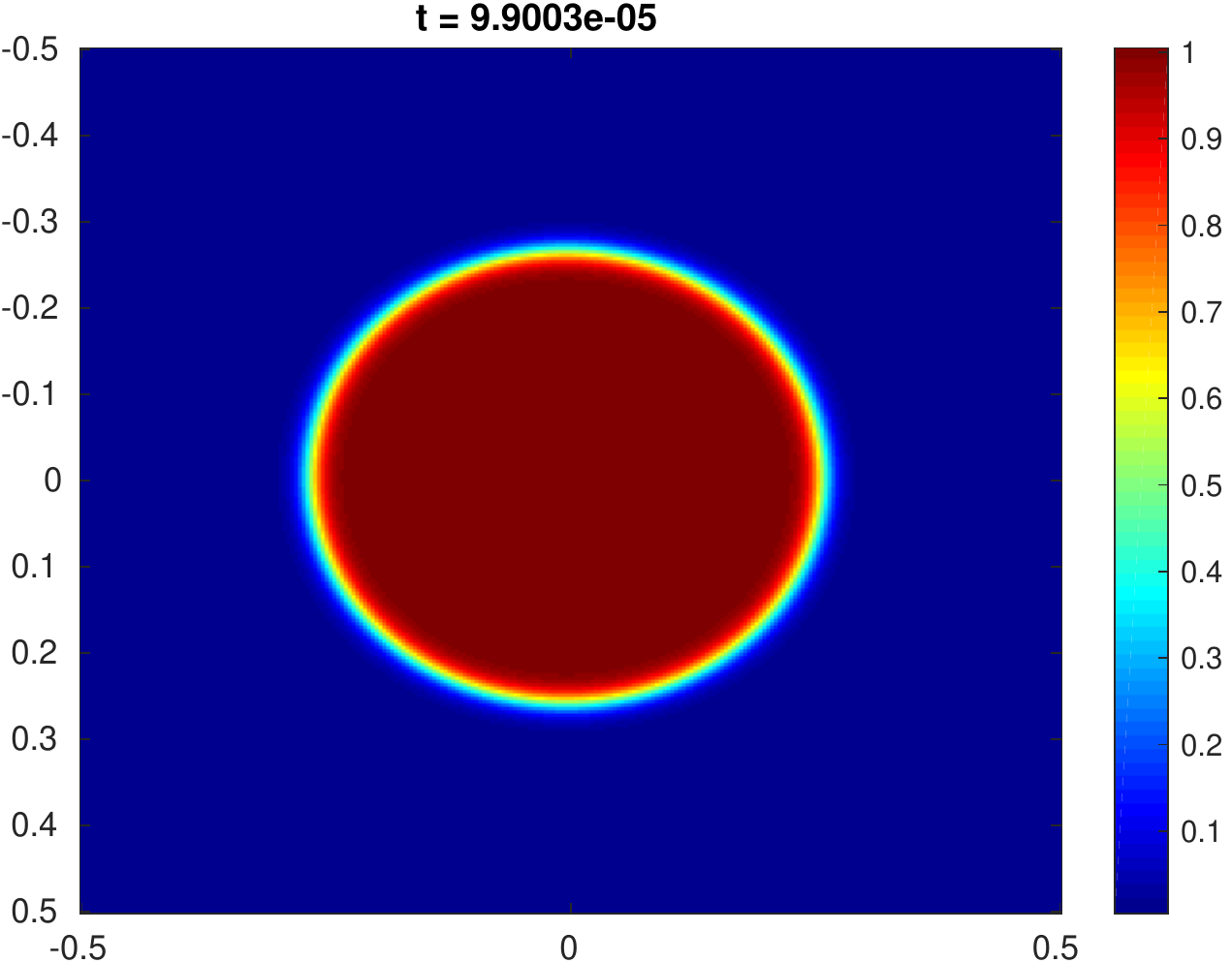} \\
	\includegraphics[width=3.5cm]{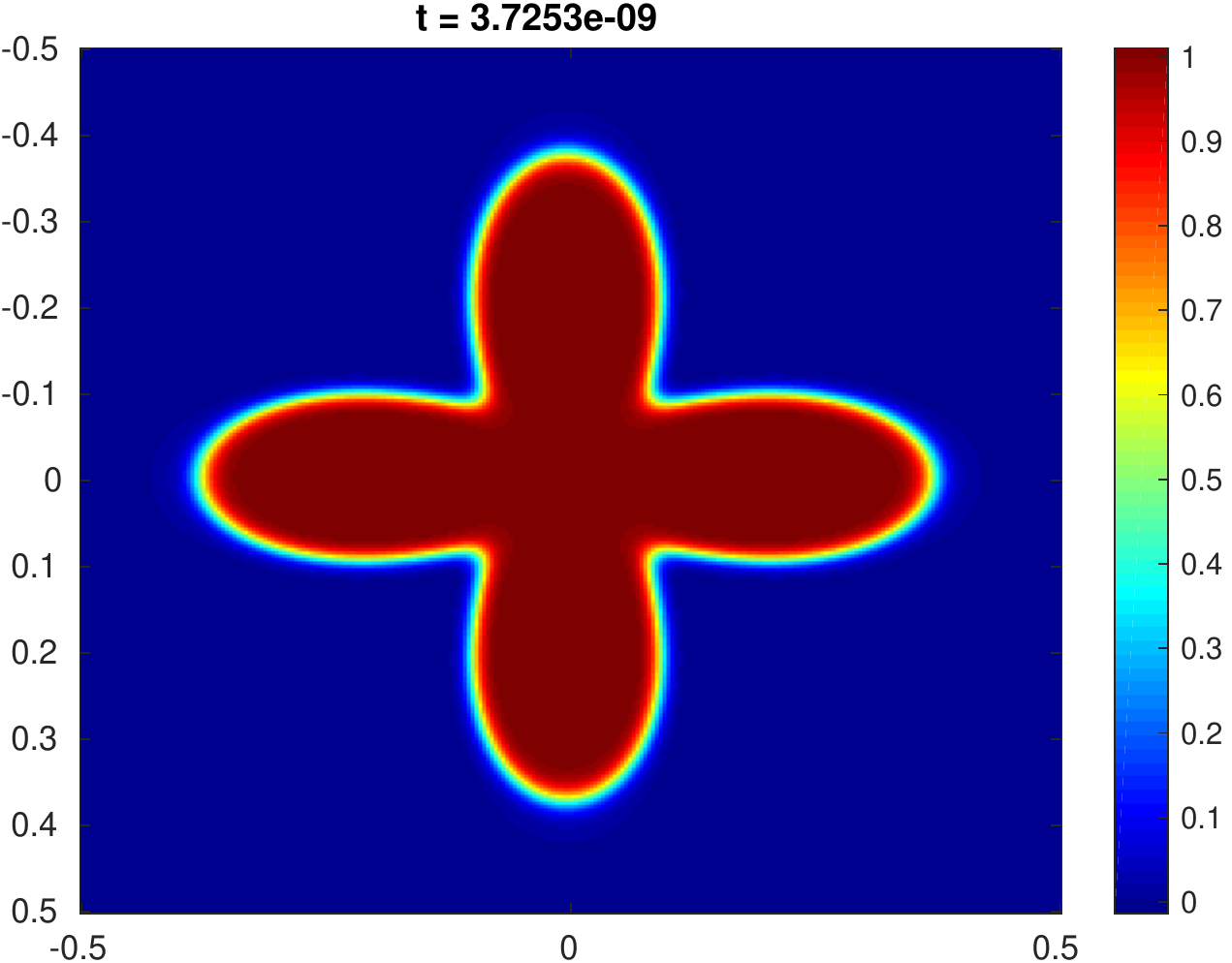}
	\includegraphics[width=3.5cm]{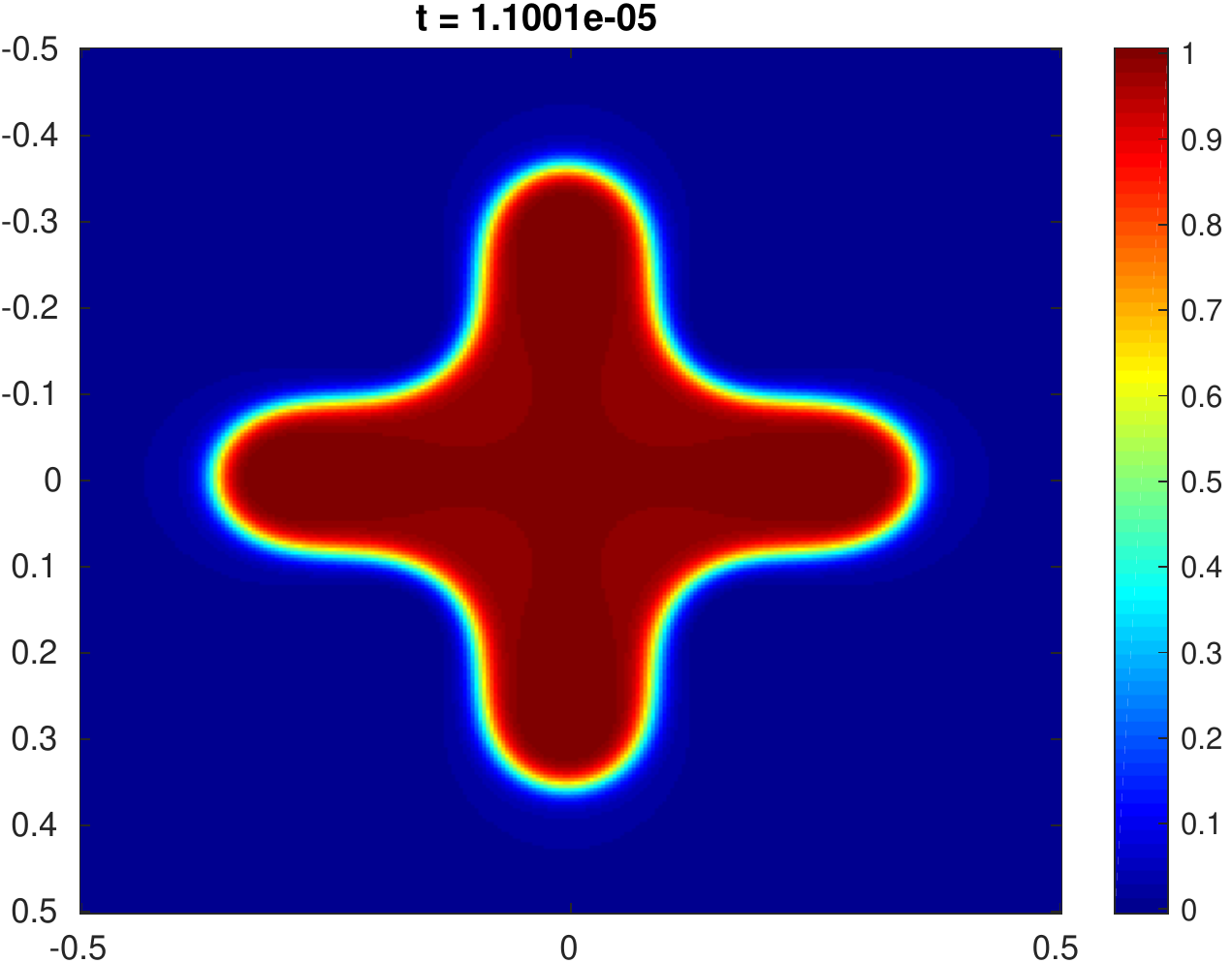}
	\includegraphics[width=3.5cm]{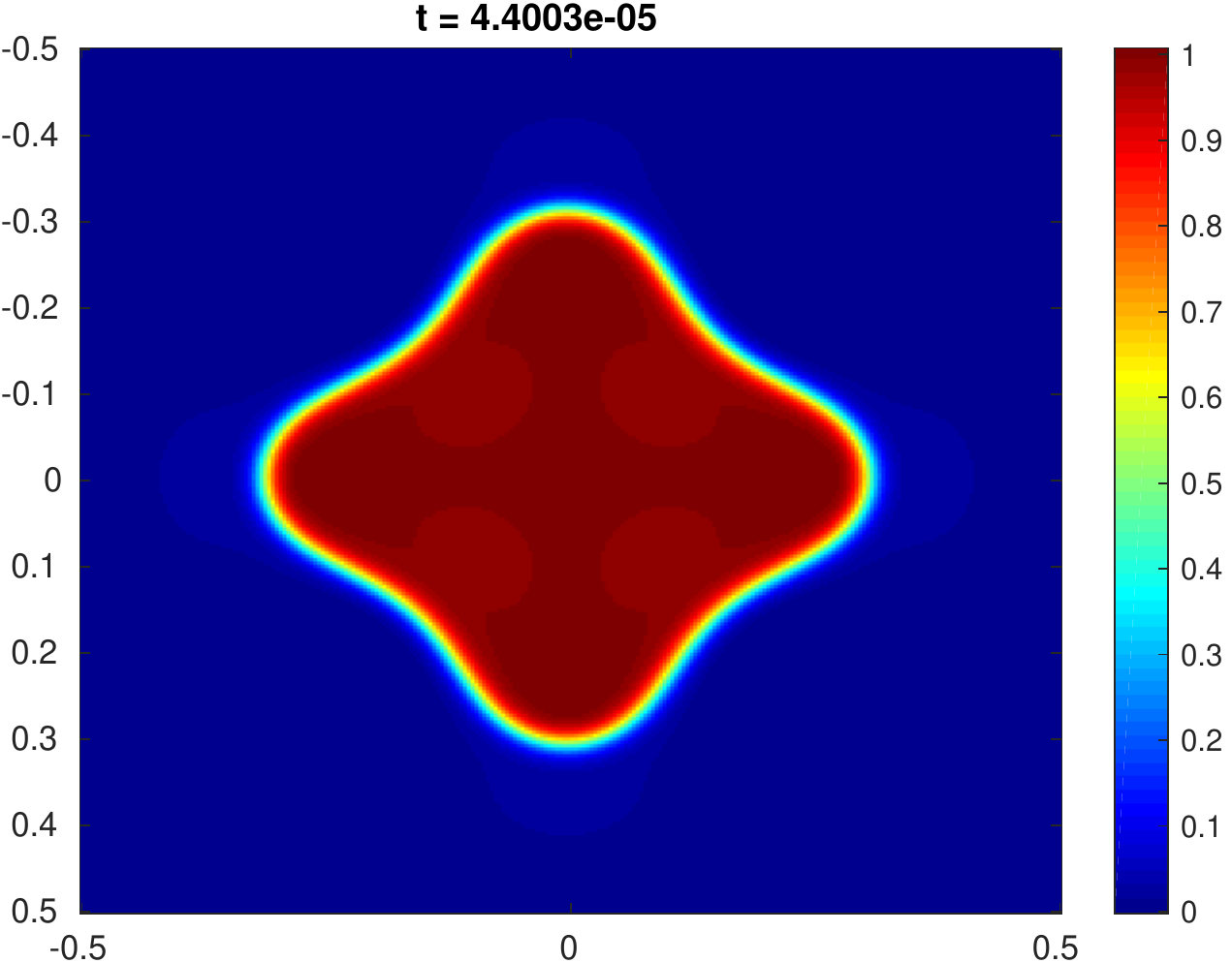}
	\includegraphics[width=3.5cm]{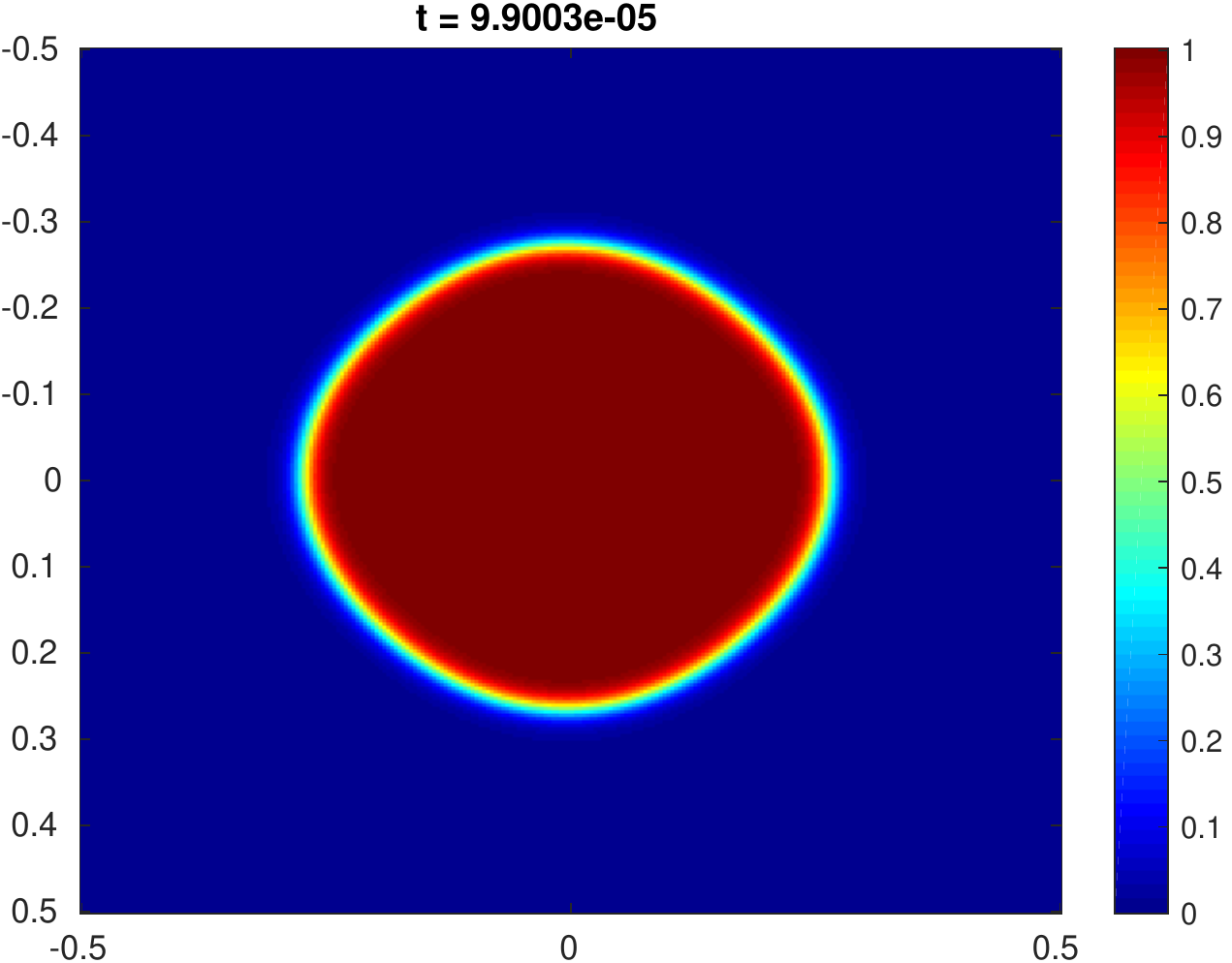} \\
	\includegraphics[width=3.5cm]{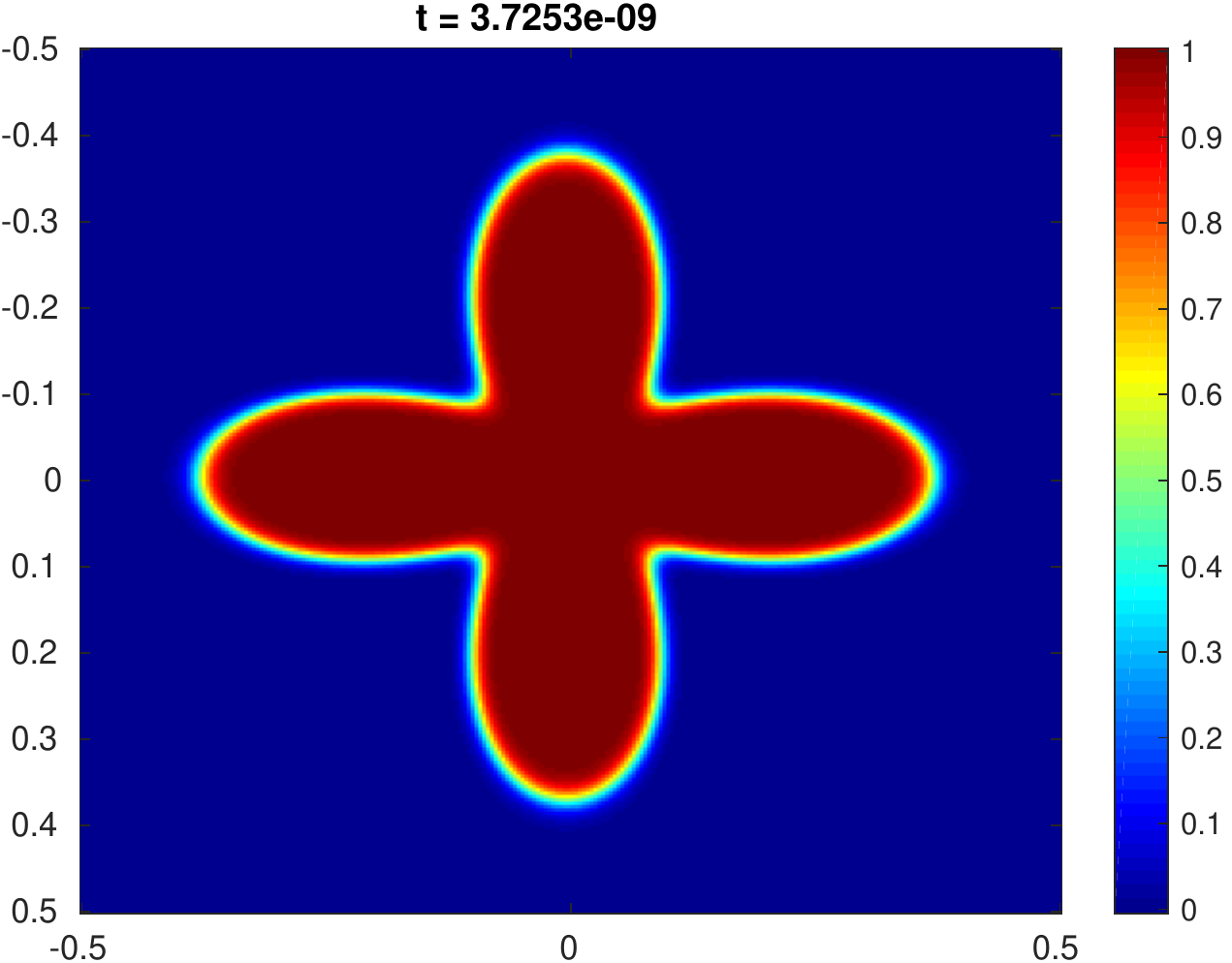}
	\includegraphics[width=3.5cm]{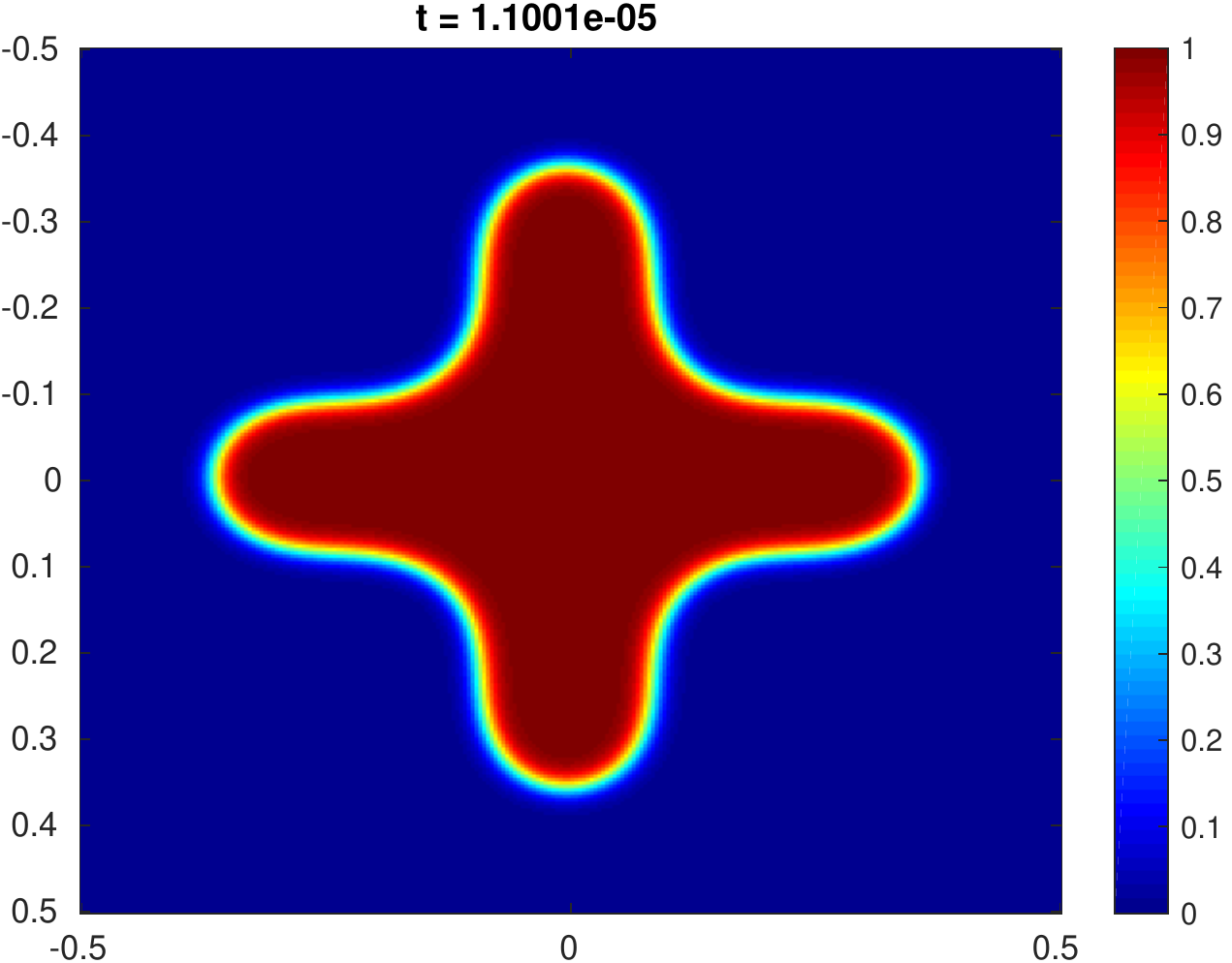}
	\includegraphics[width=3.5cm]{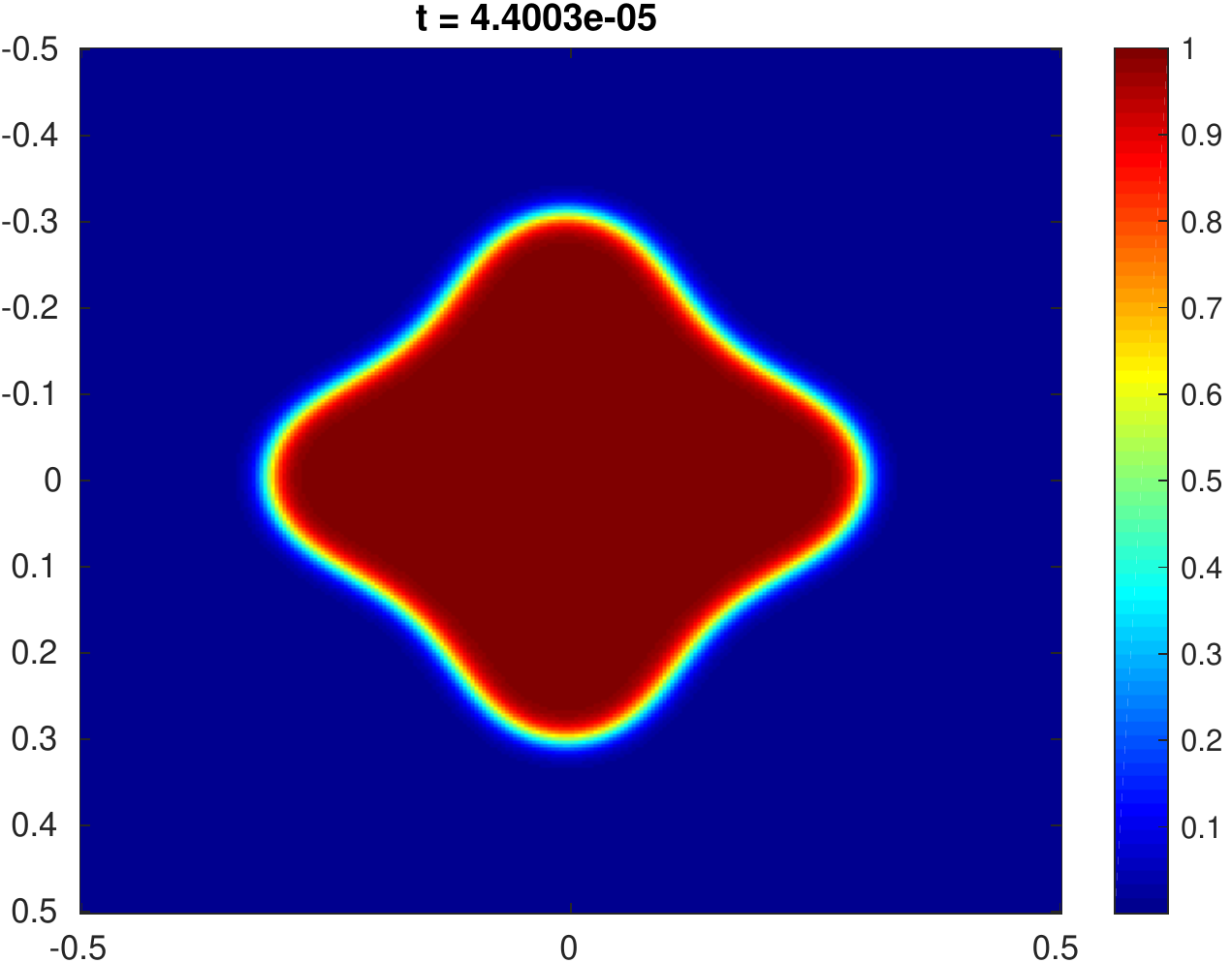}
	\includegraphics[width=3.5cm]{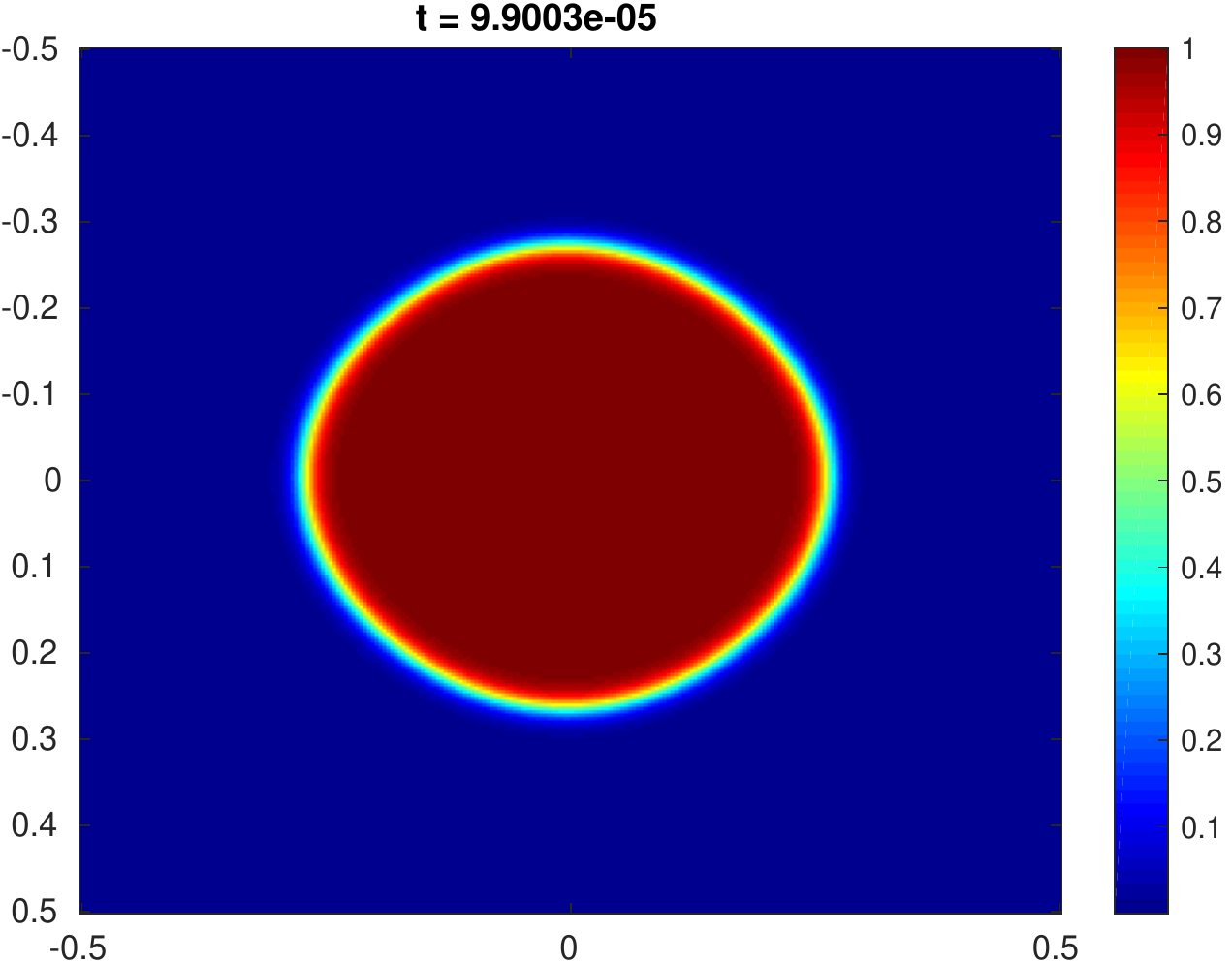} \\
\caption{First numerical comparison of the three different {\bf CH} models: Evolution of $u$ along the iterations;
First line using the {\bf C-CH} model, Second line, using the {\bf M-CH} model;
last line using the {\bf NMN-CH} model.}
\label{fig_test1}
\end{figure}

The first numerical example concerns the evolution of an initial connected set. For each Cahn--Hilliard model,
we plot on figure \eqref{fig_test1} the phase field function $u^{n}$ computed at different times $t$. 
Each experiment is performed using the same numerical parameters: 
$\delta_x = \frac{1}{2^8}$, $\epsilon = 2 \delta_x$, $\delta_t = \epsilon^4$,
$\alpha = 2/\epsilon^2$, $m=1$, and $\beta= 2/\epsilon^2$. The first, second and third lines on \eqref{fig_test1} correspond respectively to
the solution $u$ given by the {\bf C-CH} model, the {\bf  M-CH} model and the {\bf NMN-CH} model.  
The first remark is that, as expected, the {\bf C-CH} model, whose limit flow is 
the Hele-Shaw model \cite{pego1989front,alikakos1994convergence}) gives a slightly different flow compared to the other two models.
On the other hand, the numerical experiments obtained using the {\bf M-CH} model and the {\bf NMN-CH} model are very similar and 
should give a good approximation of the  surface diffusion flow.
In addition, for each model, the stationary flow limit appears to correspond to a ball of the same mass as that of the initial set. \\

\begin{figure}[htbp]
\centering
	\includegraphics[width=4.5cm]{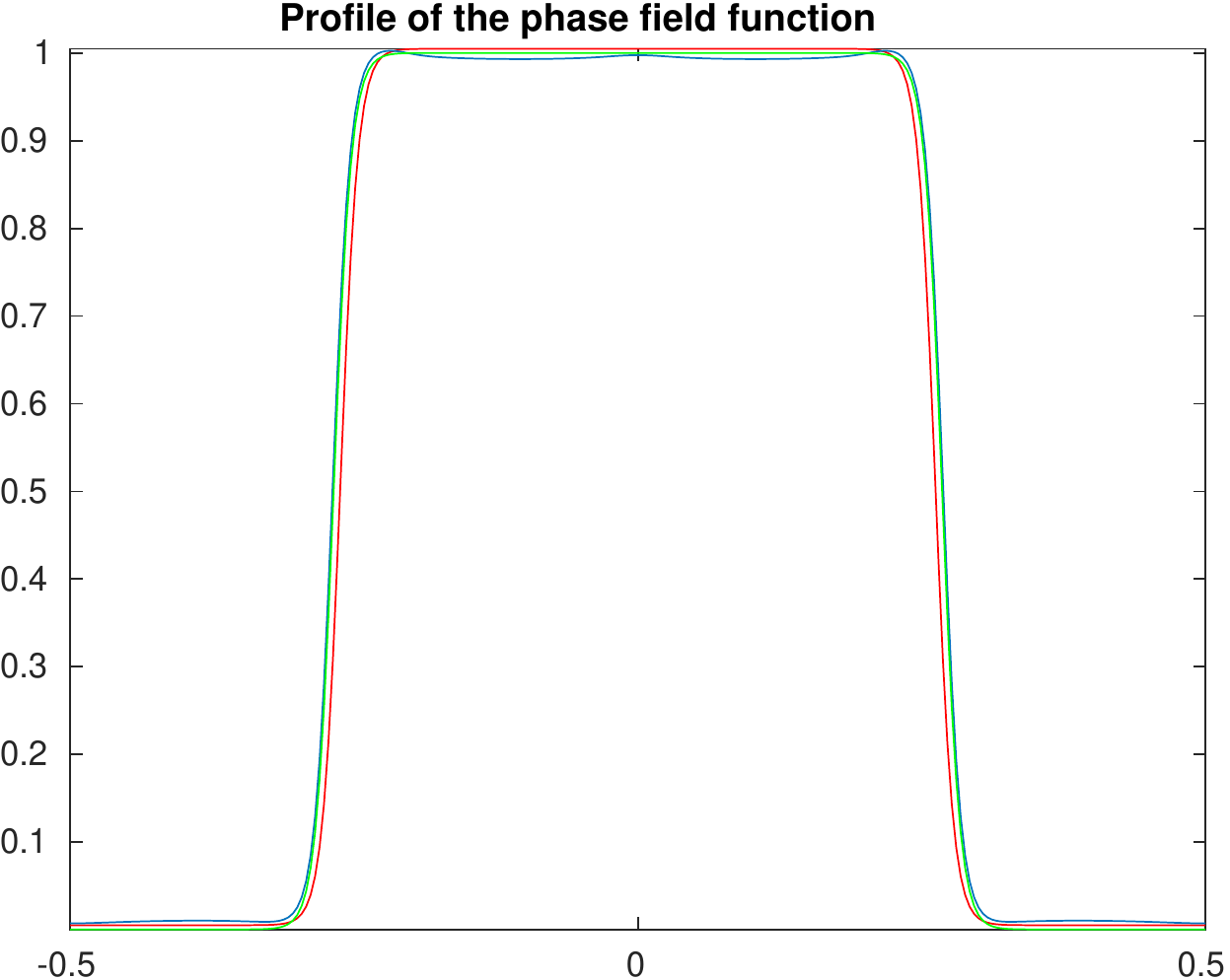}
	\includegraphics[width=4.5cm]{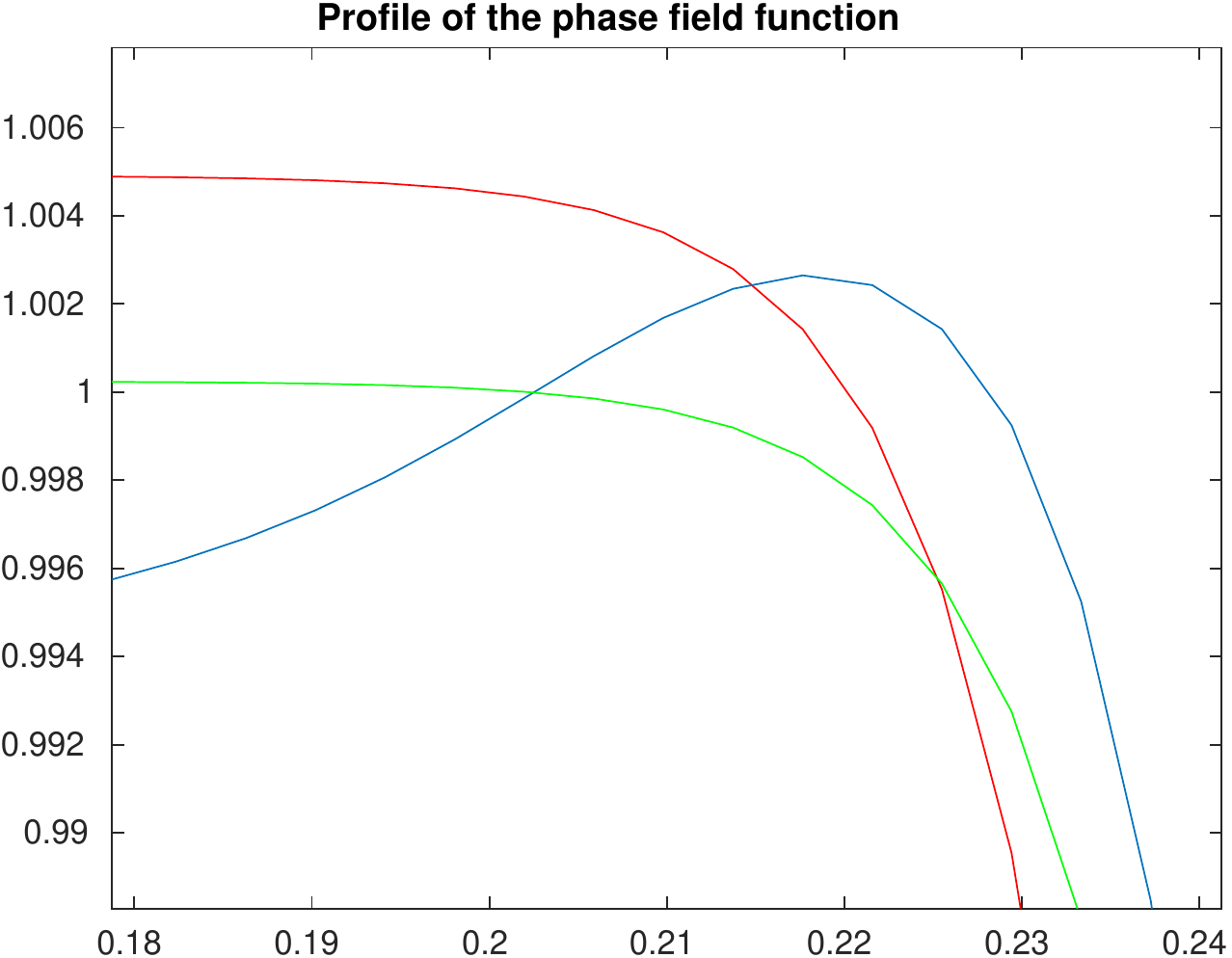}
	\includegraphics[width=4.5cm]{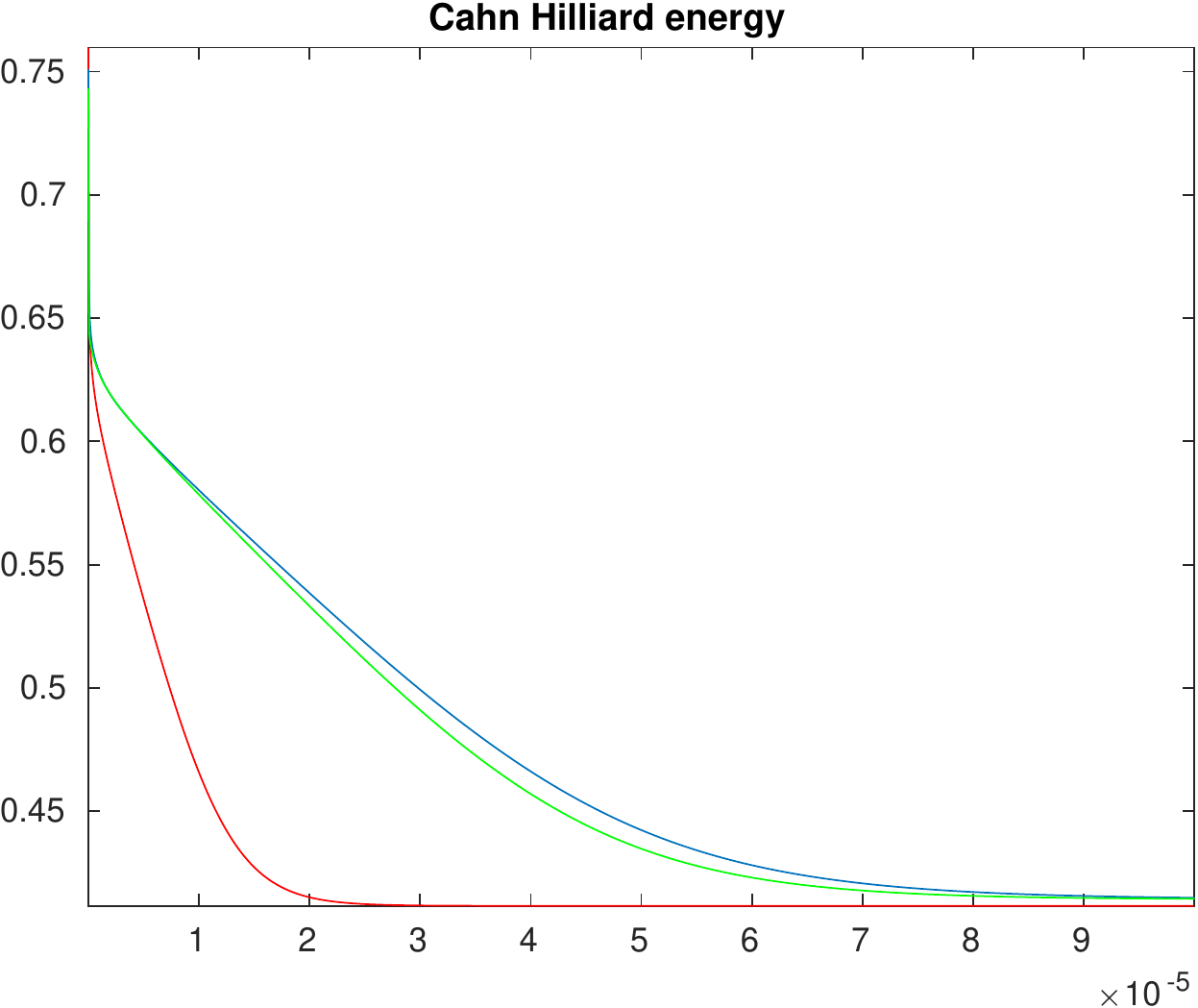}
\caption{Comparison of the three different models: profil and energy;  C-CH model in red,  M-CH model in blue, NMN-CH model in green; 
First figure: slice of $u$: $x_1 \mapsto u(x_1,0)$ ; Second figure: zoom on the slice of $u$;  last figure: evolution of the Cahn--Hilliard
energy along the flow.}
\label{fig_test1_profile}
\end{figure}

To illustrate  the asymptotic expansion performed in Section \ref{Sec:DeuxMob}, we plot on \eqref{fig_test1_profile} (first two pictures) 
the slice $x_1 \mapsto u(x_1,0)$  at the final time $T = 10^{-4}$. The profil associated to the {\bf C-CH} model  is plotted in red
and clearly indicates that the solution $u$ does not remain in the interval $[0,1]$  with an overshoot of order $O(\varepsilon)$.  
As for the {\bf M-CH} model (in blue), we can also observe a perturbation of order $O(\varepsilon)$ of the best profile  $q(z)$
and $u$ does not remain in $[0,1]$. 
In contrast, the profile obtained using the {\bf NMN} model (in green) seems to be very close to $q$ and remains in $[0,1]$ up to an error of order $O(\epsilon^2)$. Finally, we plot the evolution of the Cahn--Hilliard energy 
along the flow for each model on the last picture of \eqref{fig_test1_profile}. We can clearly observe a decrease of the energy in each case.\\

In conclusion, this first numerical experiment confirms the asymptotic expansion obtained in the previous section,
and highlights the interest of our {\bf NMN} model to approximate surface diffusion flows.

\subsubsection{Influence of the mobility: a local conservation of mass}

\begin{figure}[htbp]
\centering
	\includegraphics[width=3.5cm]{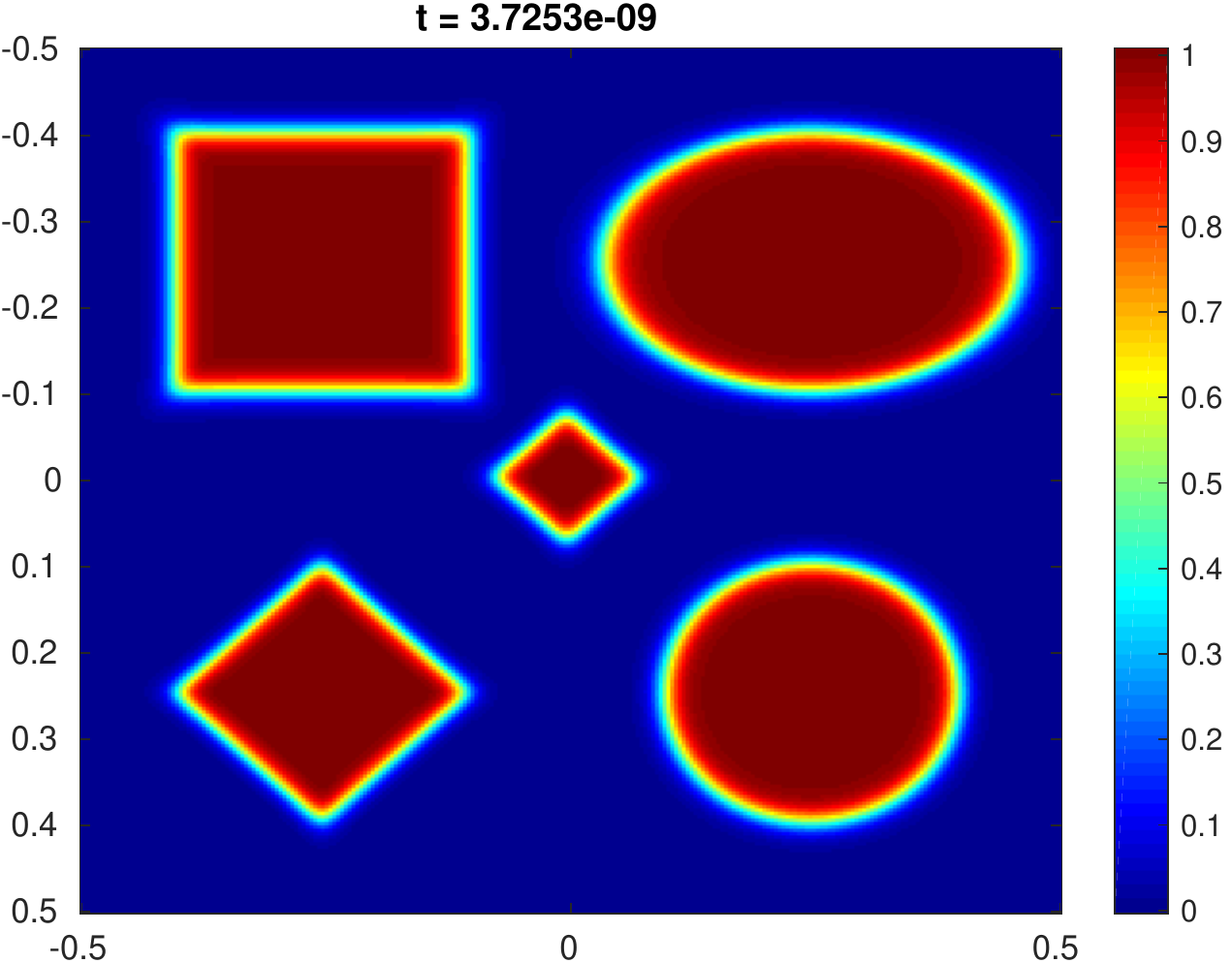}
	\includegraphics[width=3.5cm]{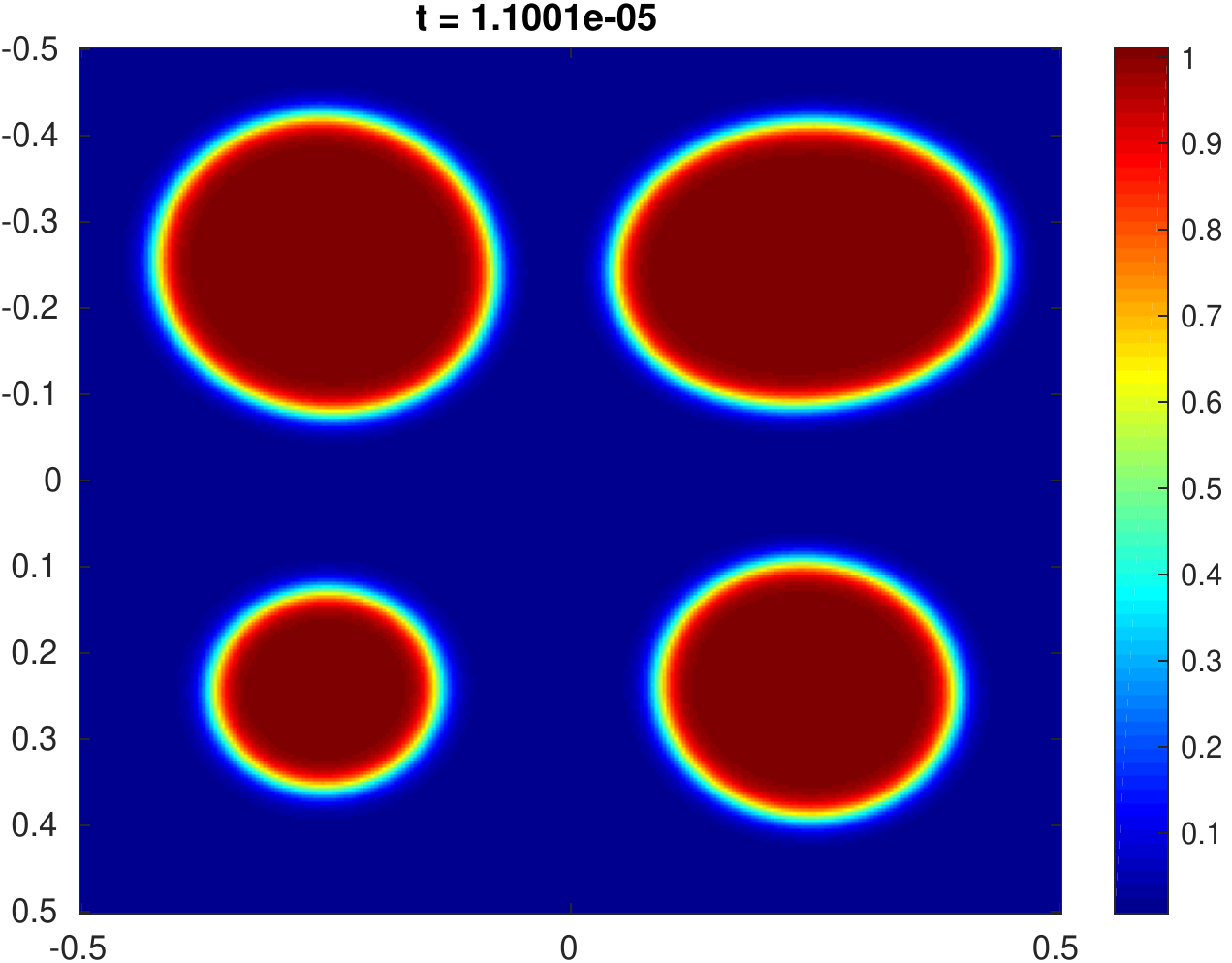}
	\includegraphics[width=3.5cm]{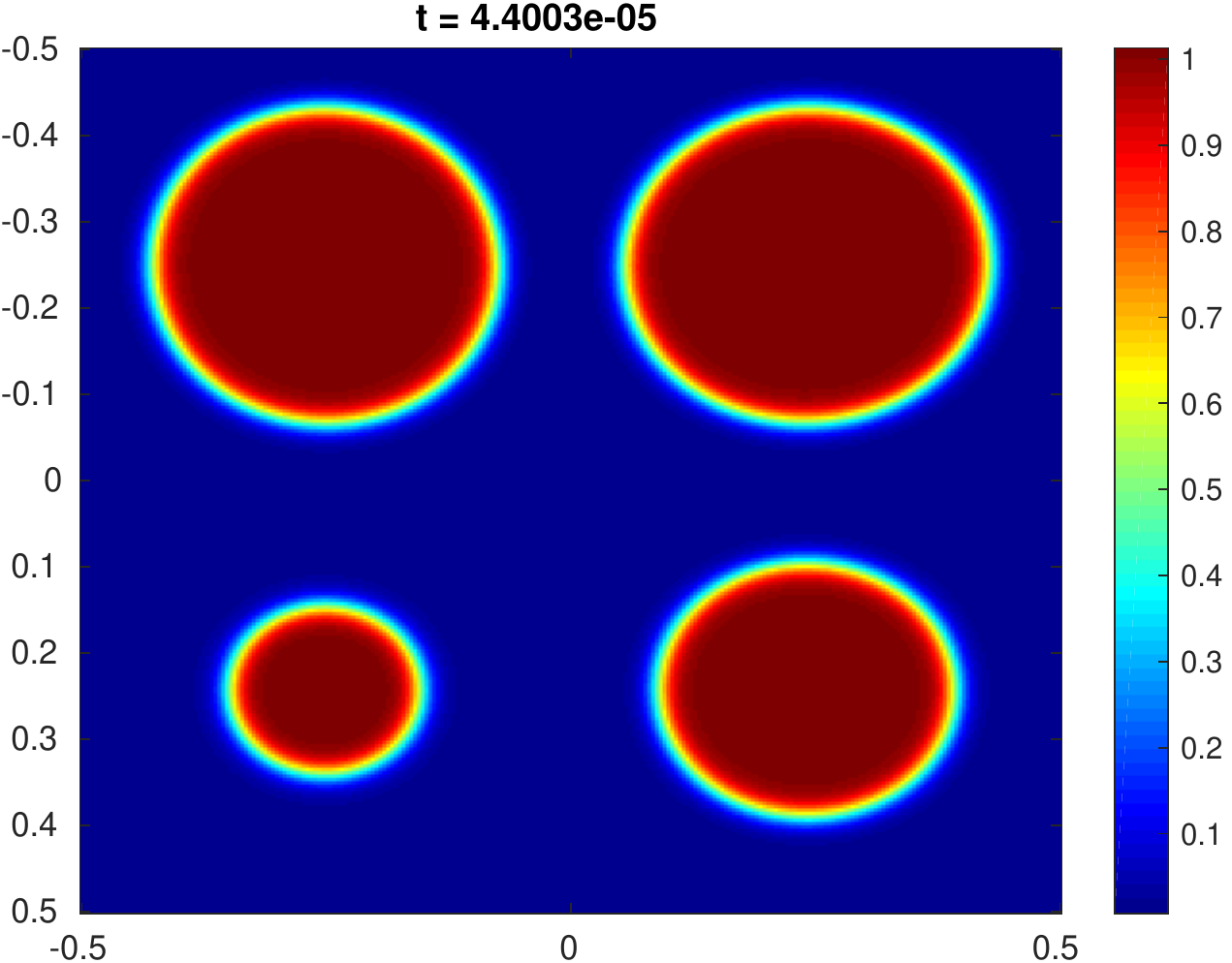}
	\includegraphics[width=3.5cm]{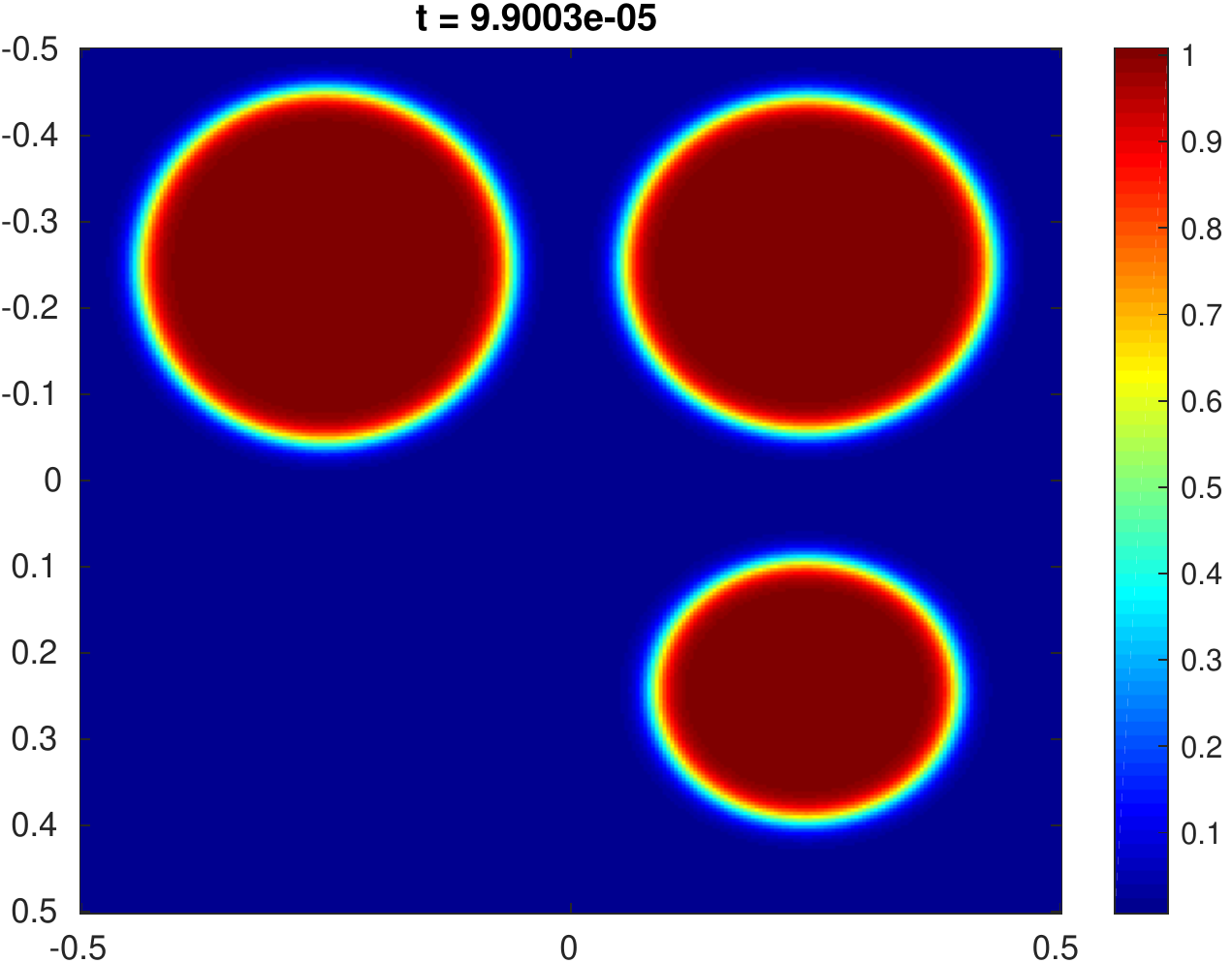} \\
	\includegraphics[width=3.5cm]{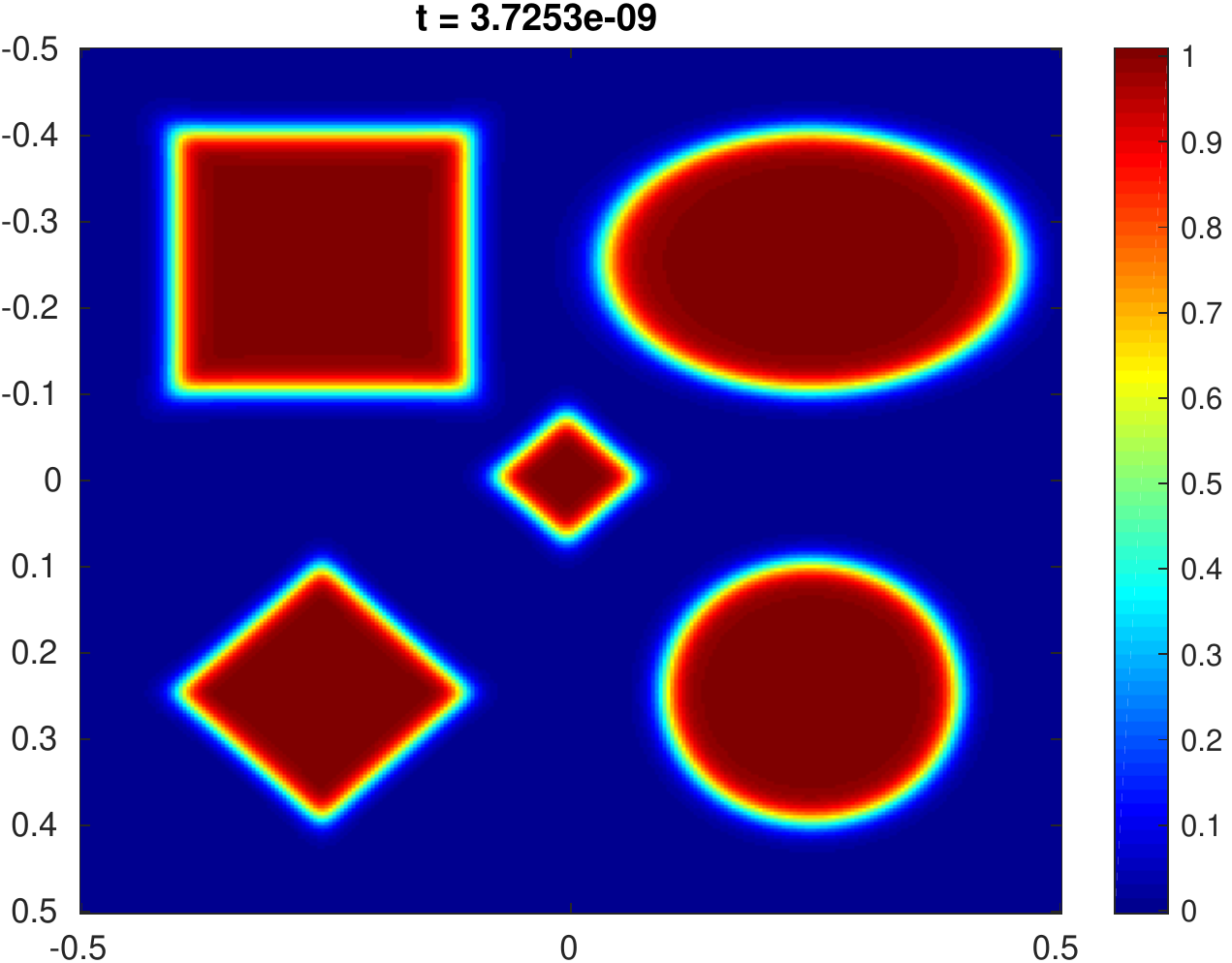}
	\includegraphics[width=3.5cm]{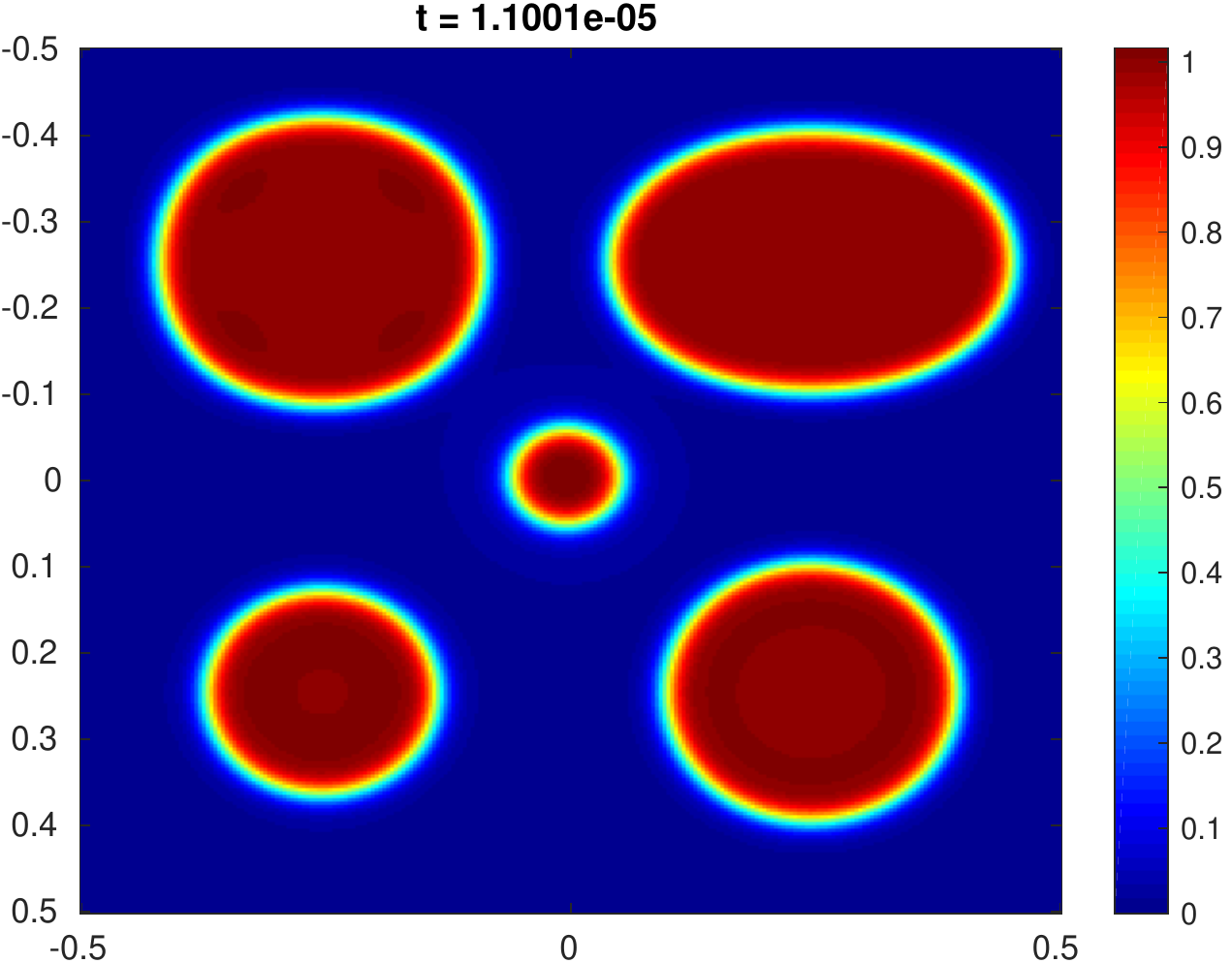}
	\includegraphics[width=3.5cm]{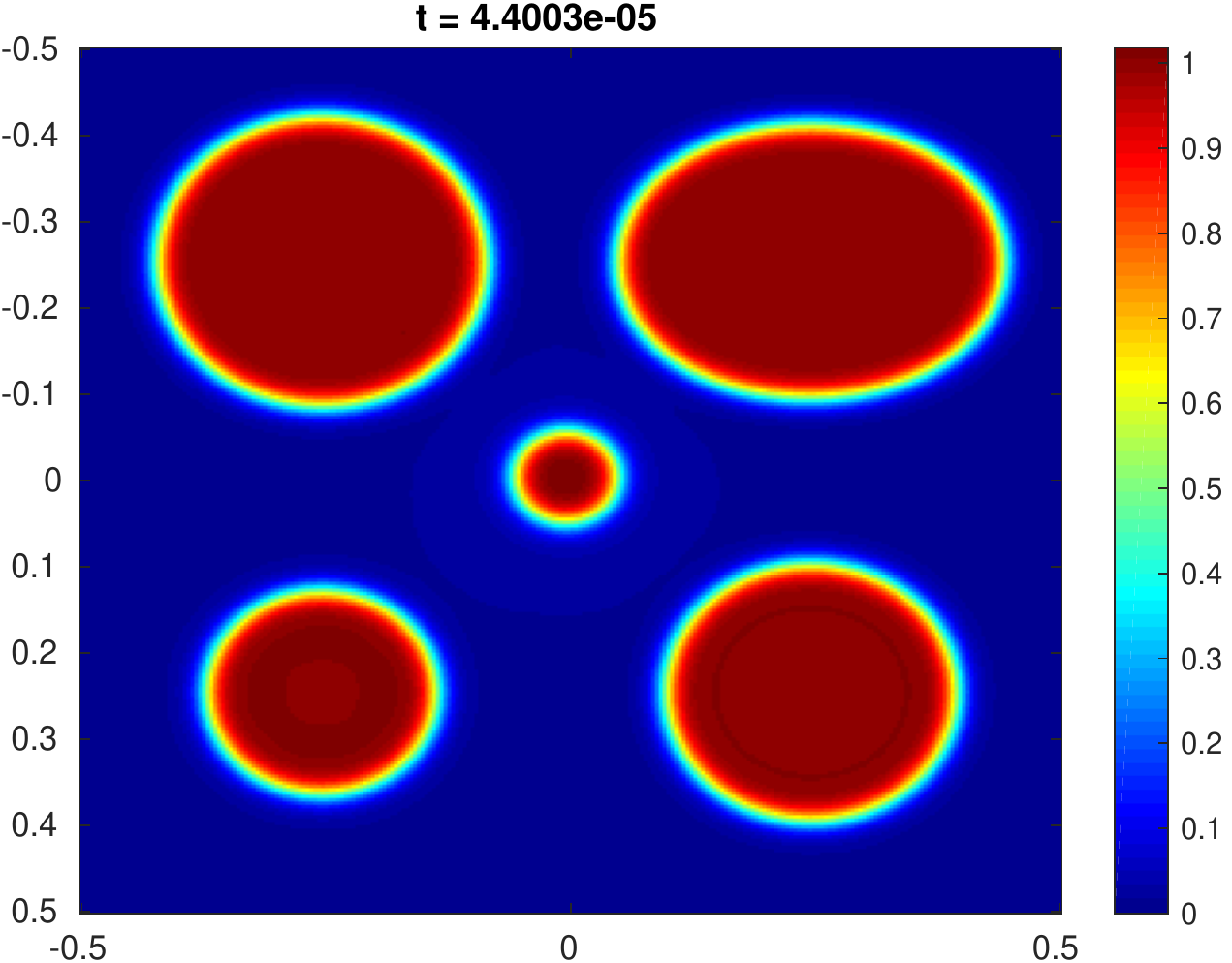}
	\includegraphics[width=3.5cm]{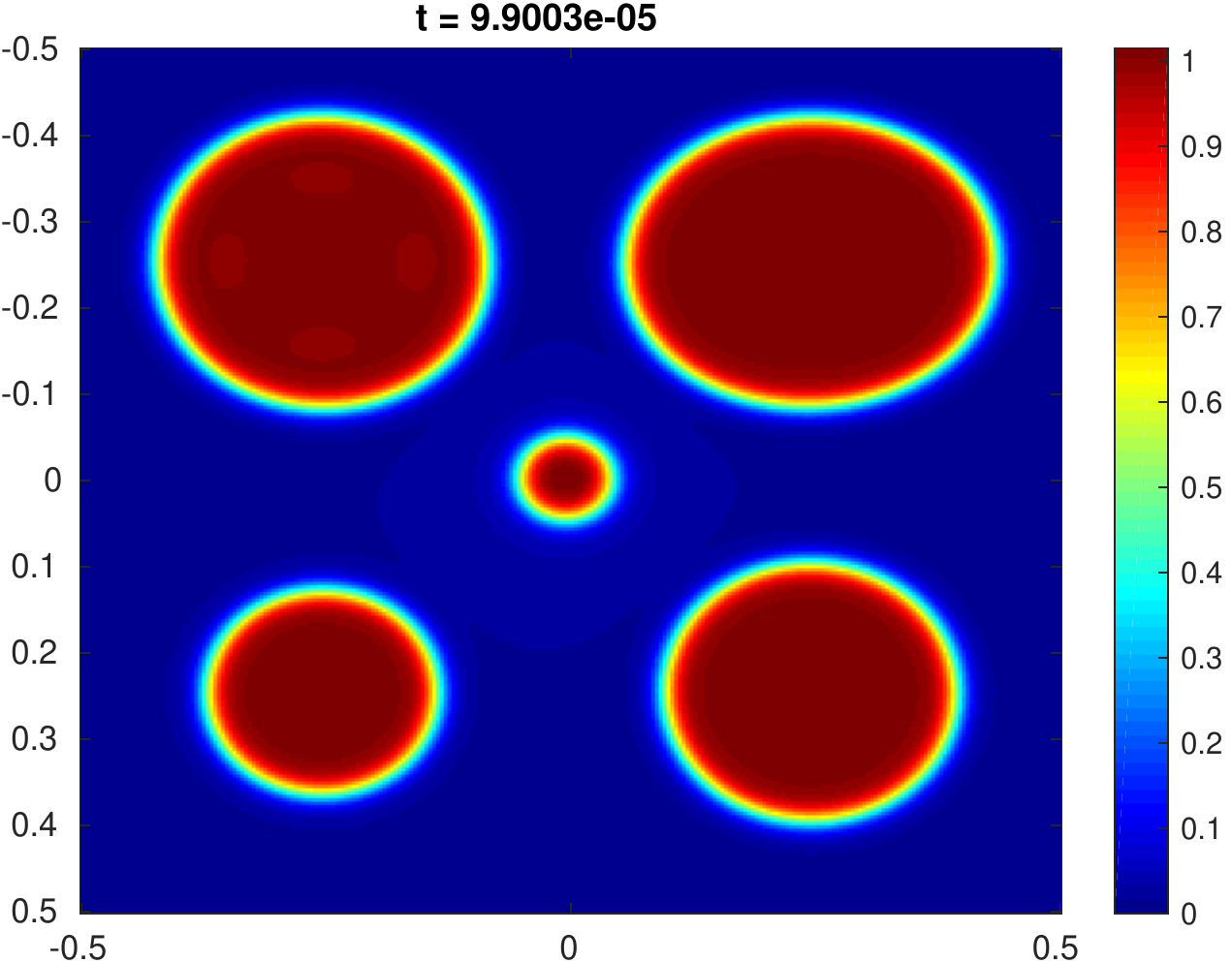} \\
	\includegraphics[width=3.5cm]{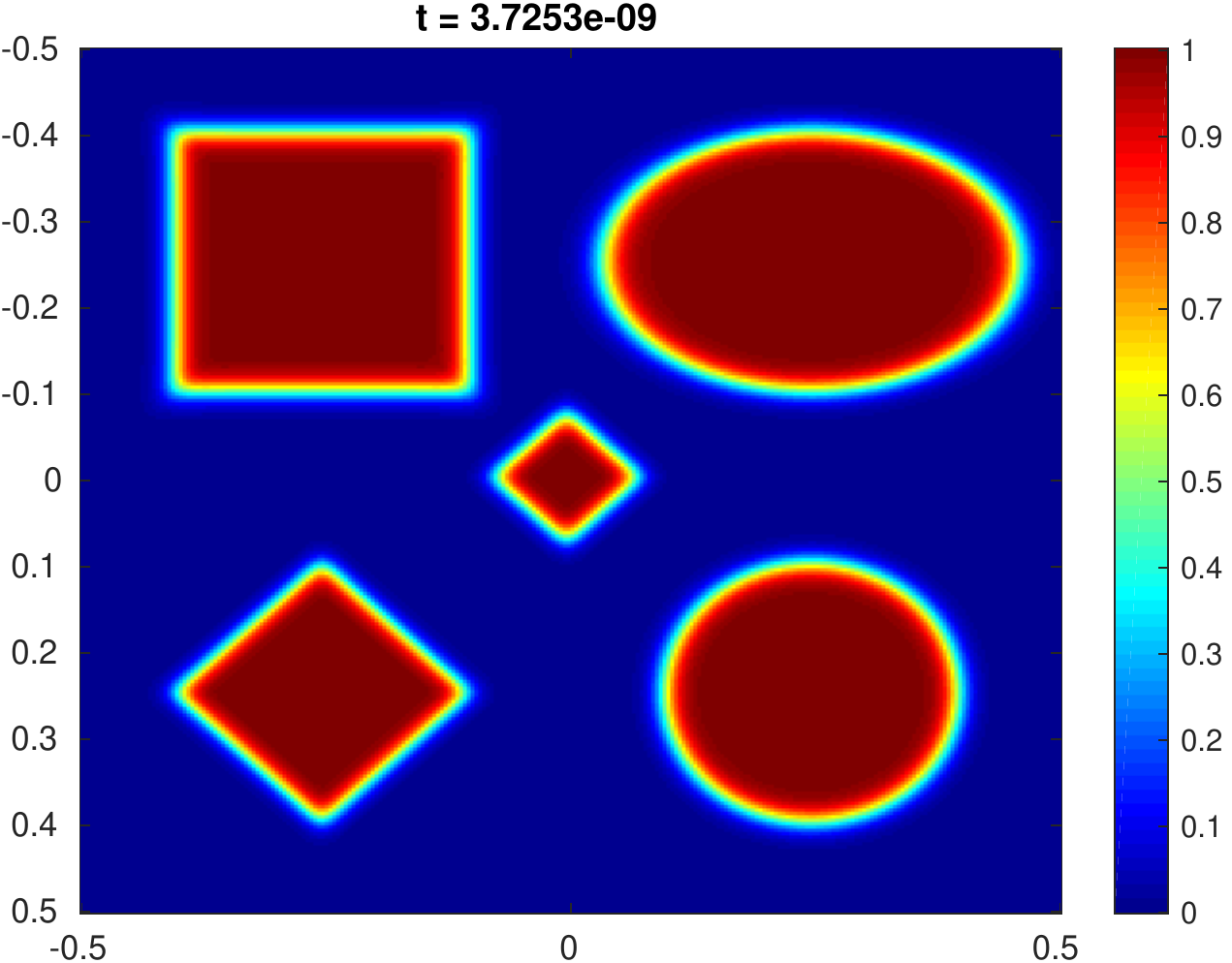}
	\includegraphics[width=3.5cm]{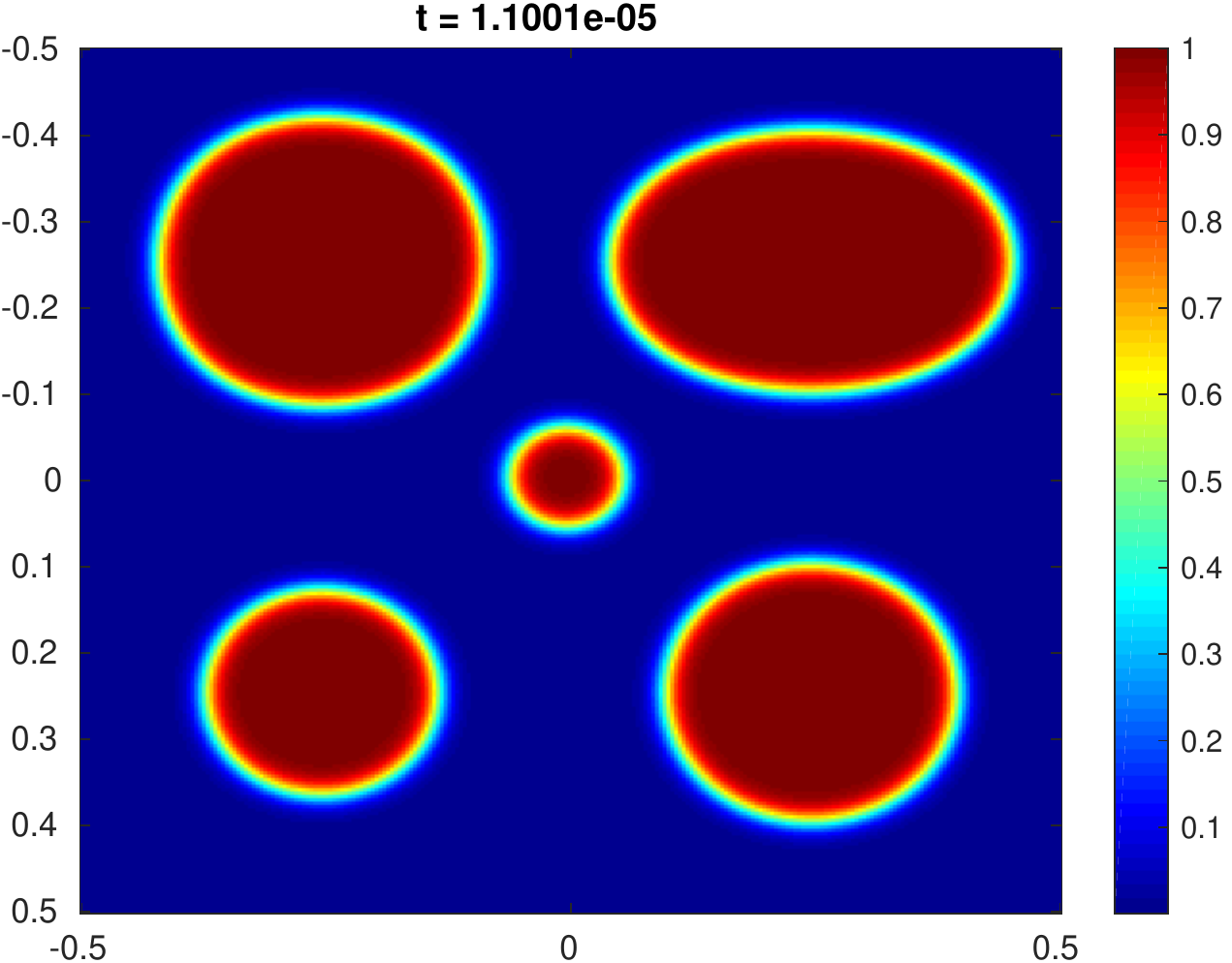}
	\includegraphics[width=3.5cm]{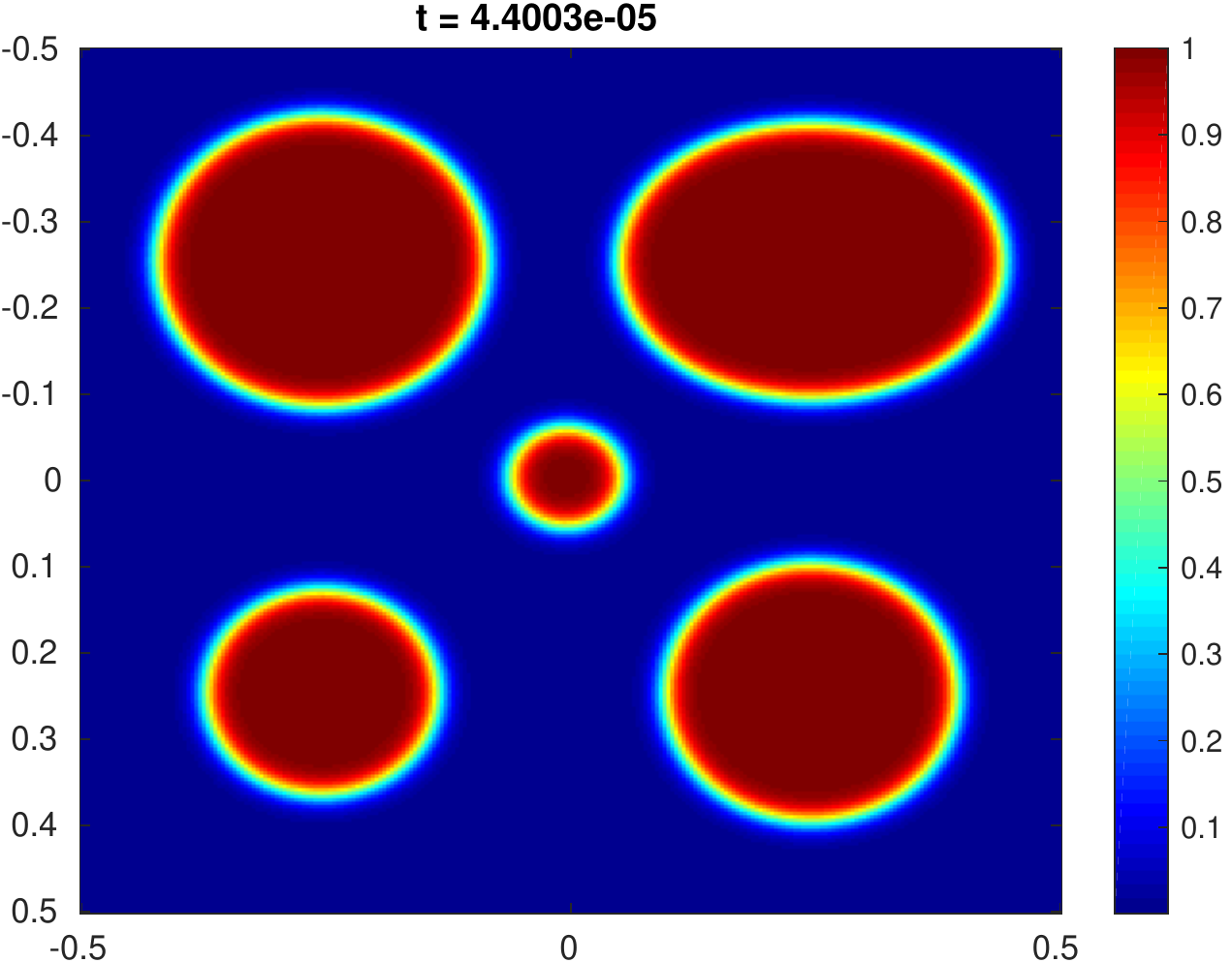}
	\includegraphics[width=3.5cm]{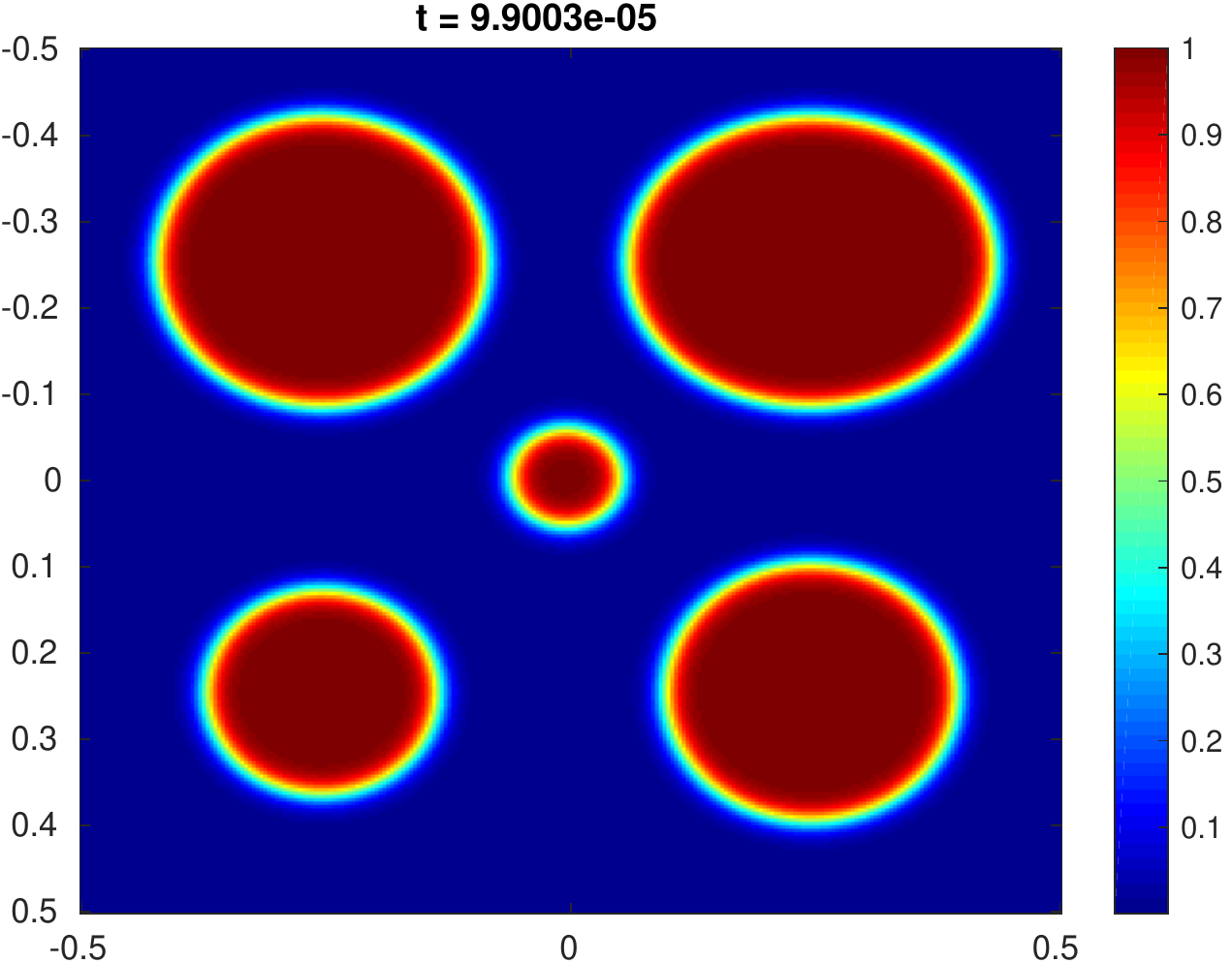} \\
\caption{Numerical comparison of the three different phase field models: Local conservation of the mass; 
Evolution of $u$ along the iterations; First line with the {\bf C-CH} model. Second line with the {\bf M-CH} model;
Last line with the {\bf NMN-CH} model.}
\label{fig_test2}
\end{figure}

The second numerical experiment is intended to show the advantage of adding mobility to the classical Cahn--Hilliard model
to preserve a local conservation of the mass. As previously, we use the same  numerical parameter in each case: 
$\delta_x = \frac{1}{2^8}$, $\epsilon = 2/N$, $\delta_t = \epsilon^4$,
$\alpha = 2/\epsilon^2$, $m=1$, and $\beta= 2/\epsilon^2$. Then we plot on figure \eqref{fig_test2} the phase field function $u$
obtained at  different times $t$ using the different phase field models (first line: {\bf C-CH} model, second line: {\bf M-CH} model, third line 
{\bf NMN-CH} model).
Here, the initial set is a disjoint union of five small sets. As expected, the evolutions obtained using the {\bf M-CH} and the 
{\bf NMN-CH} models show an independent evolution of each small disjoint set that converges to a ball of equivalent volume.
This last point is clearly not the case using the {\bf C-CH} model where the limit appears to be the union of three balls only. 
It suggests that the mass of the smaller set moves towards the larger set.  This emphasizes the interest of adding mobility in the Cahn--Hilliard model to get a local conservation of 
mass, which is particularly relevant for various physical applications, for example the simulation of dewetting phenomena.

\subsubsection{Numerical experiments with thin structures in dimension $3$}
\begin{figure}[htbp]
\centering
	\includegraphics[width=3.5cm]{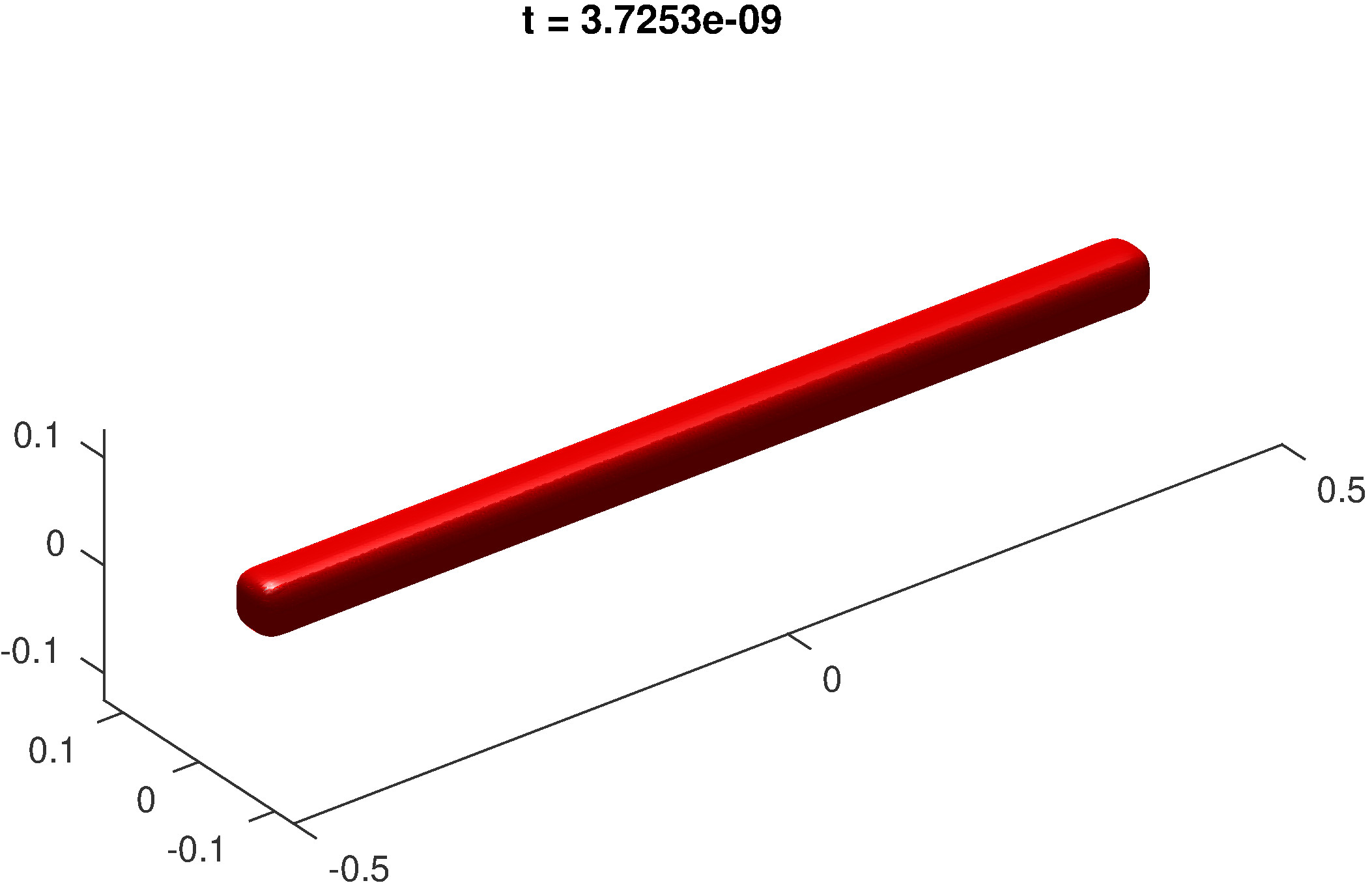}
	\includegraphics[width=3.5cm]{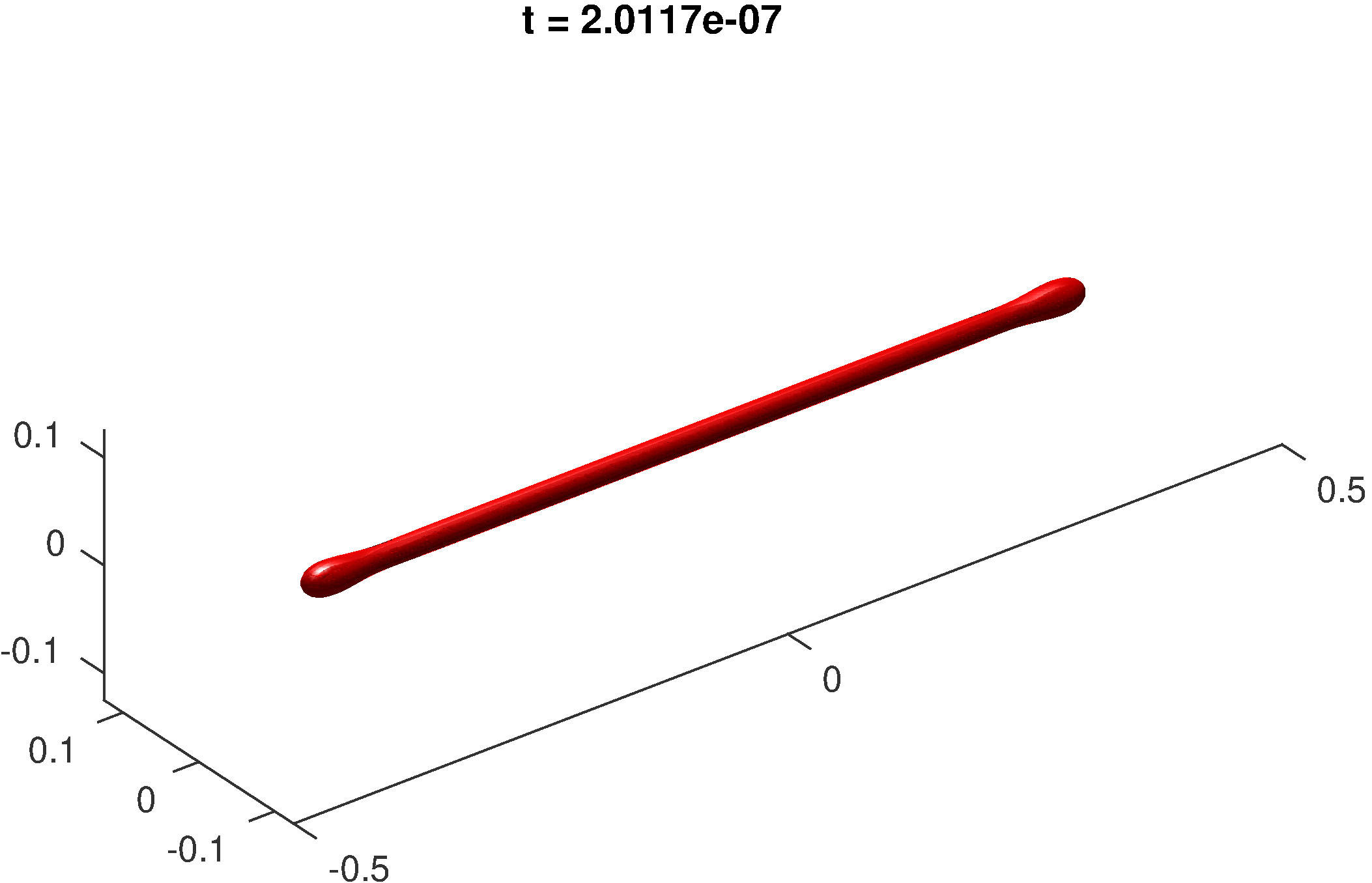}
	\includegraphics[width=3.5cm]{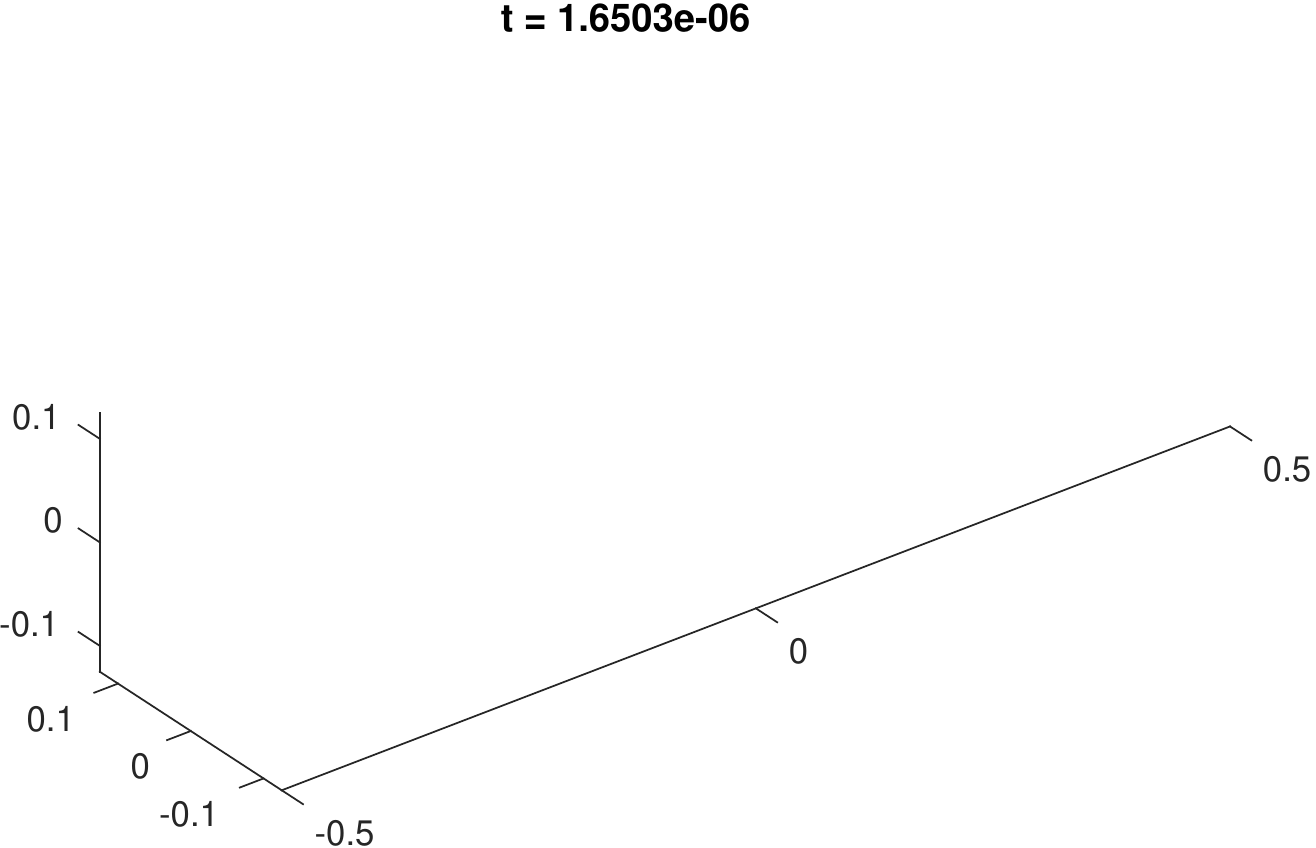}
	\includegraphics[width=3.5cm]{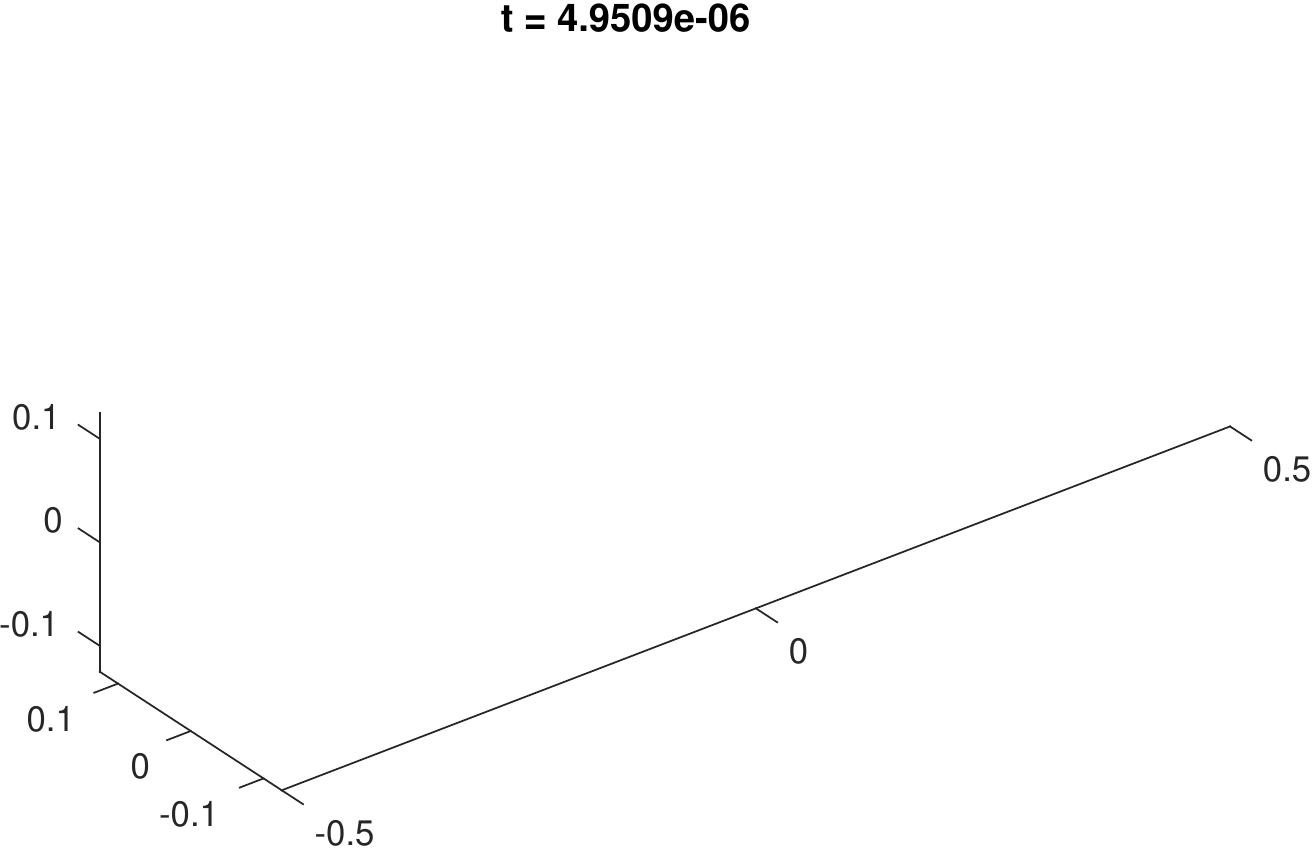} \\
    \includegraphics[width=3.5cm]{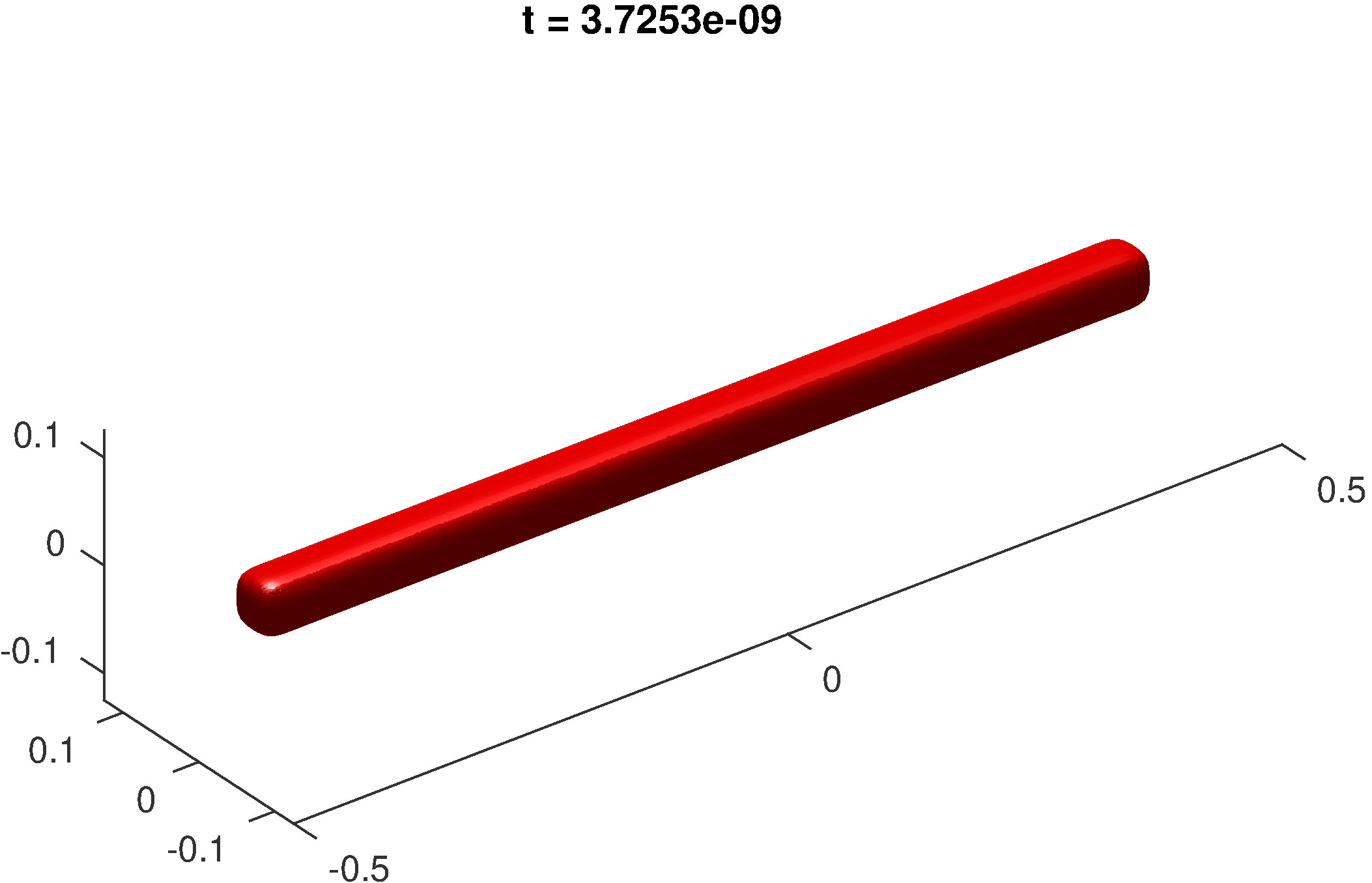}
	\includegraphics[width=3.5cm]{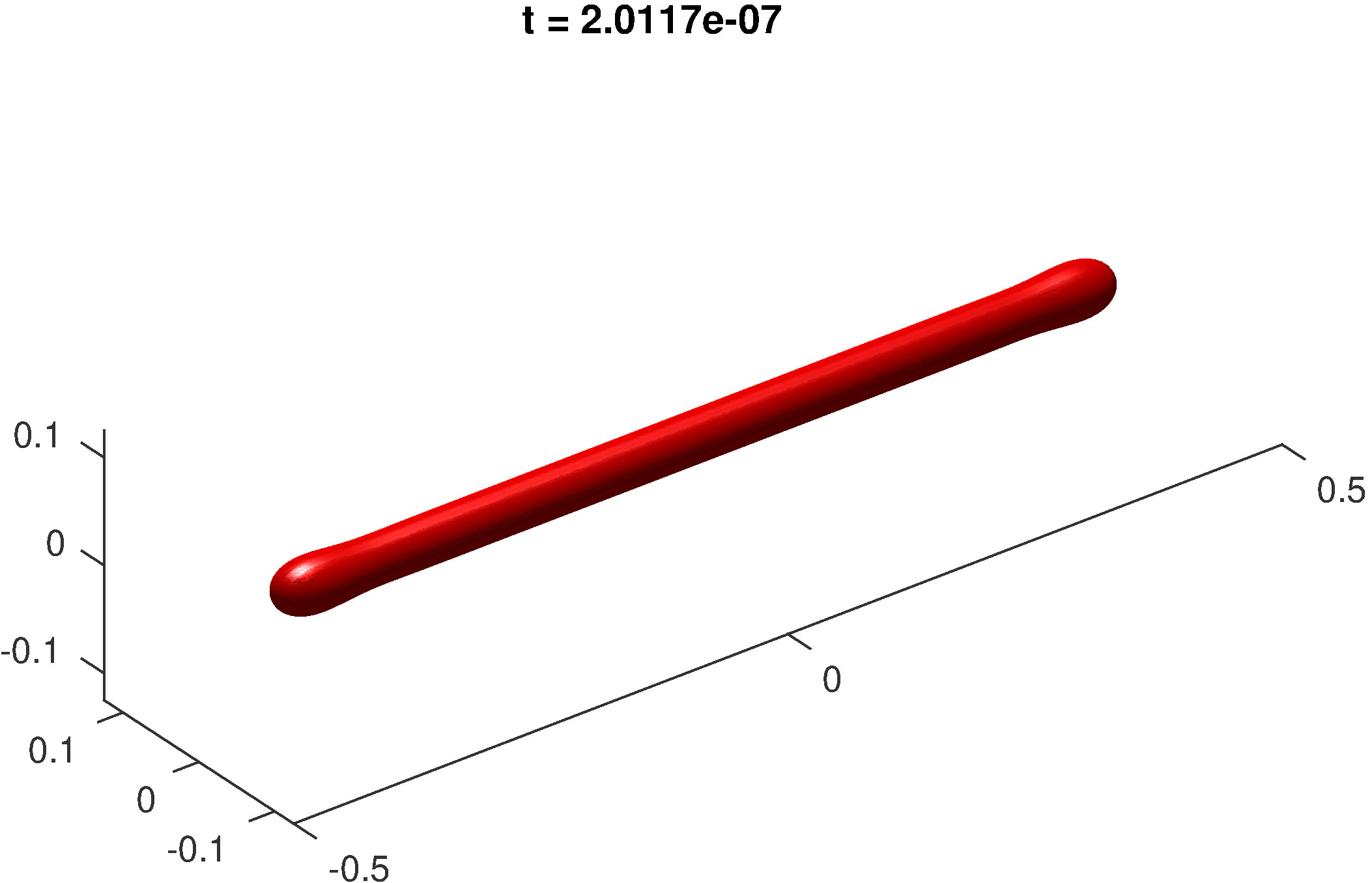}
	\includegraphics[width=3.5cm]{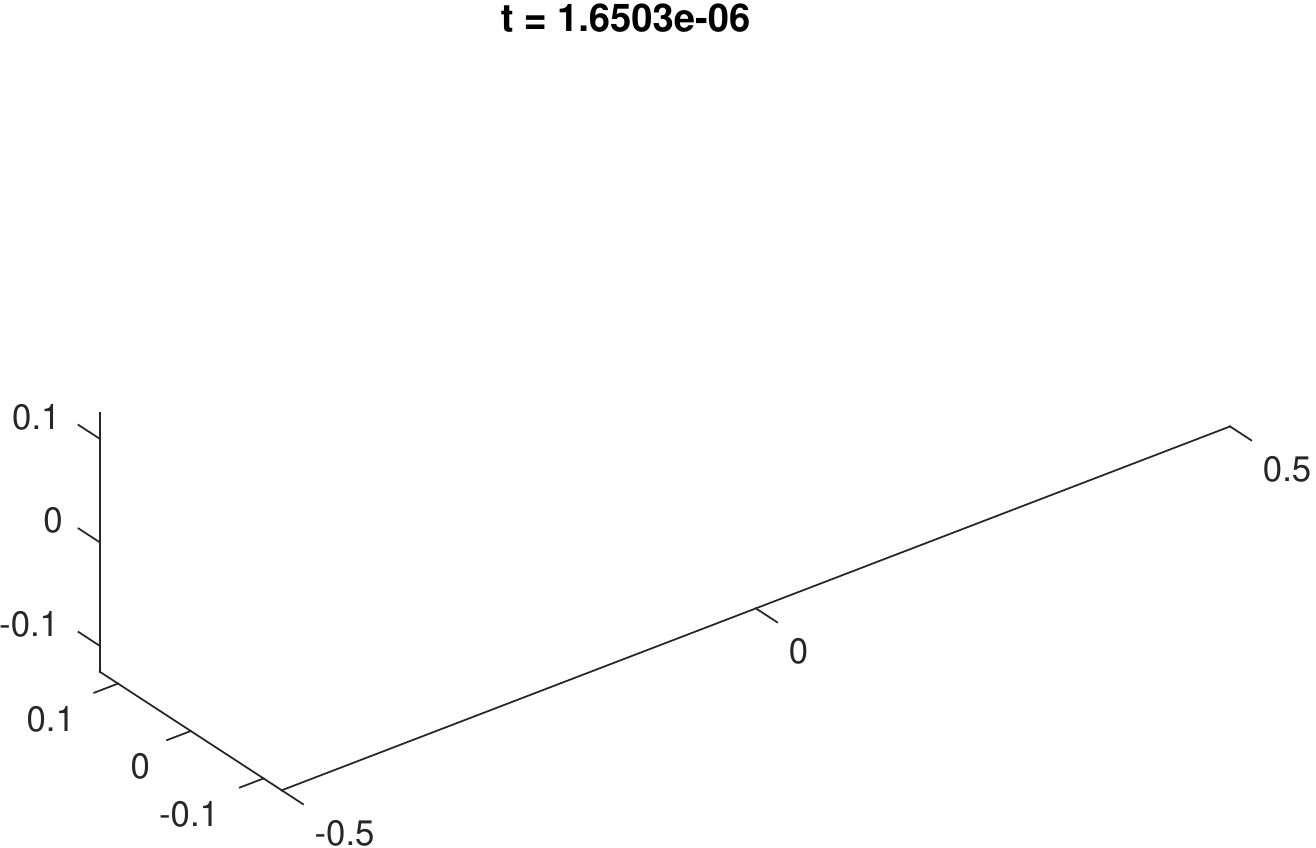}
	\includegraphics[width=3.5cm]{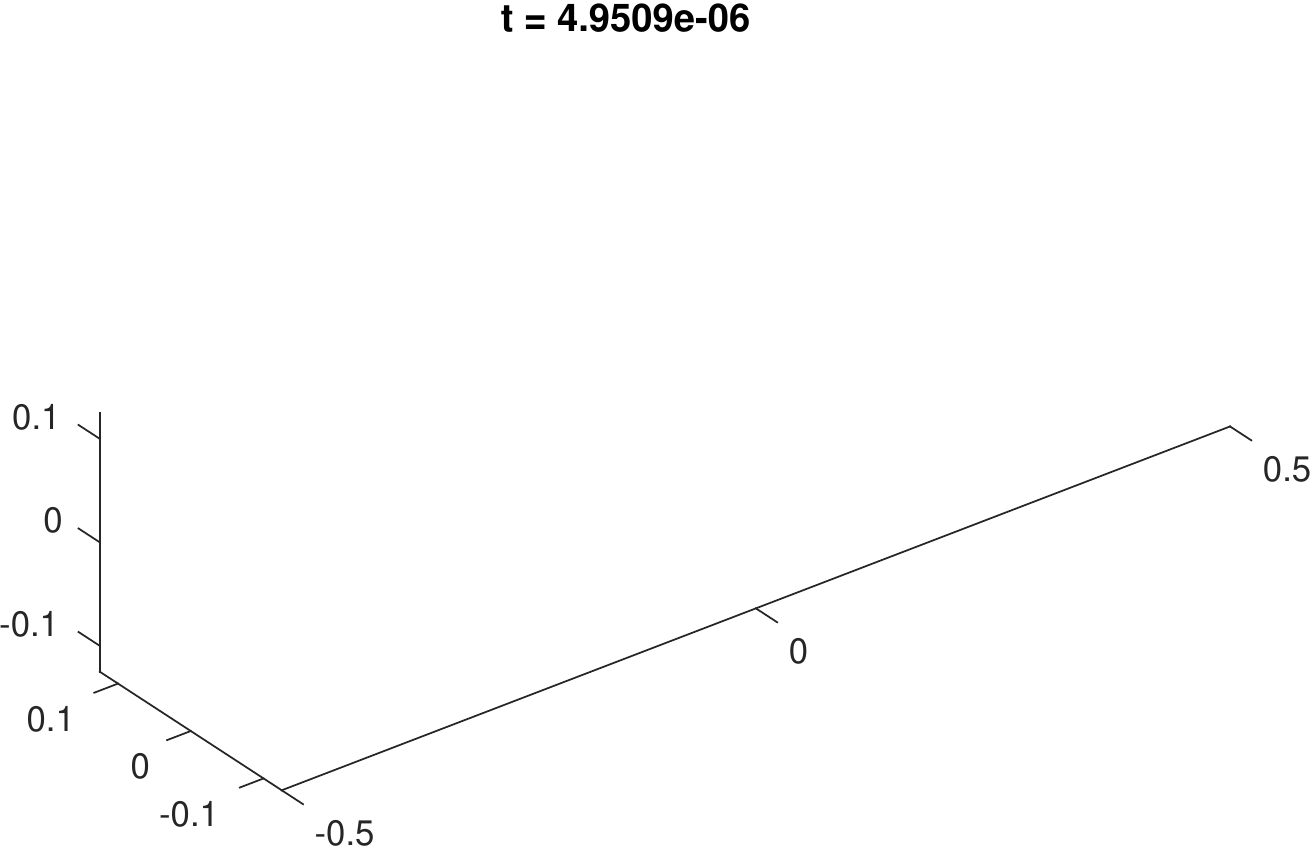} \\
    \includegraphics[width=3.5cm]{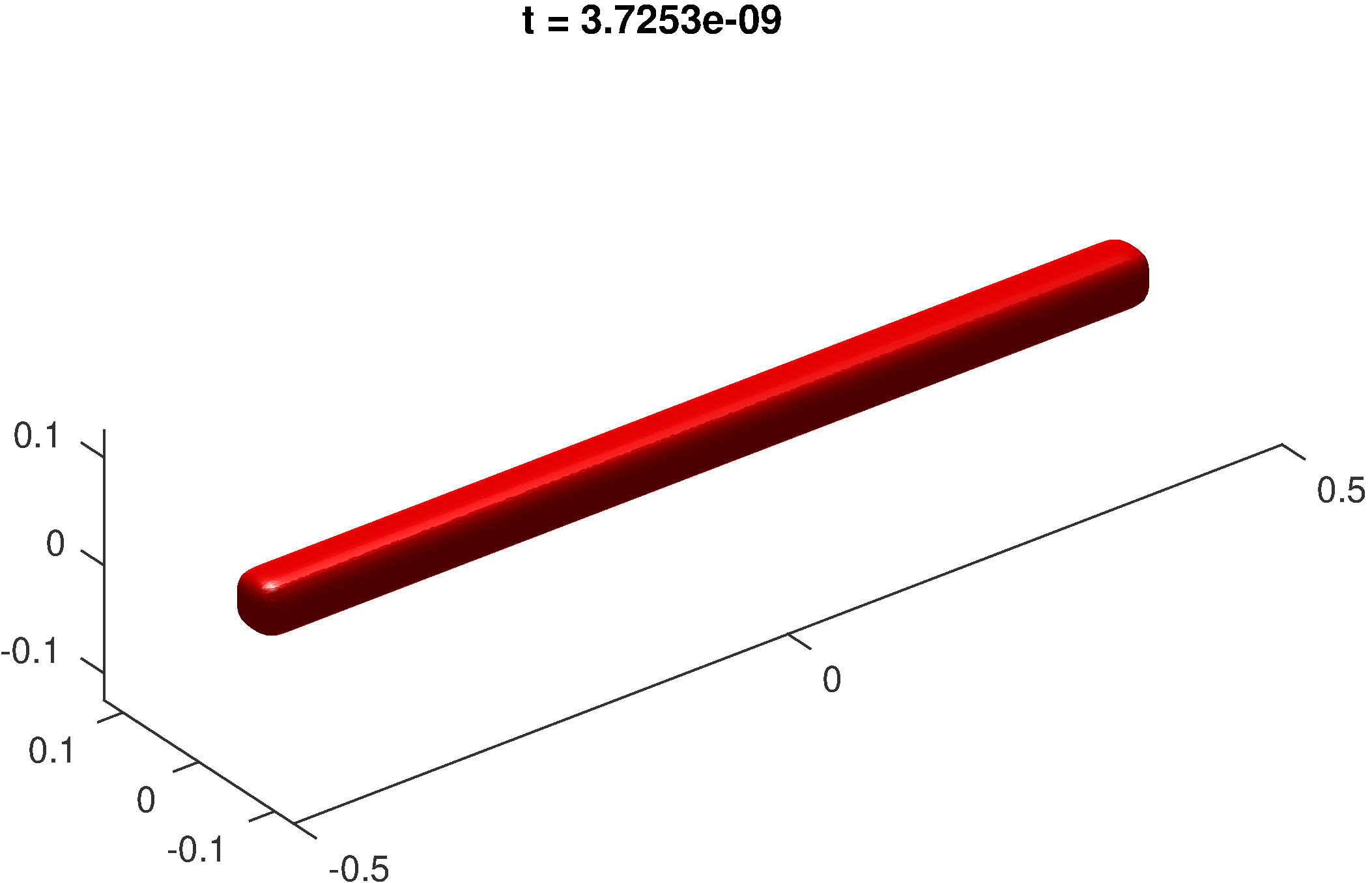}
	\includegraphics[width=3.5cm]{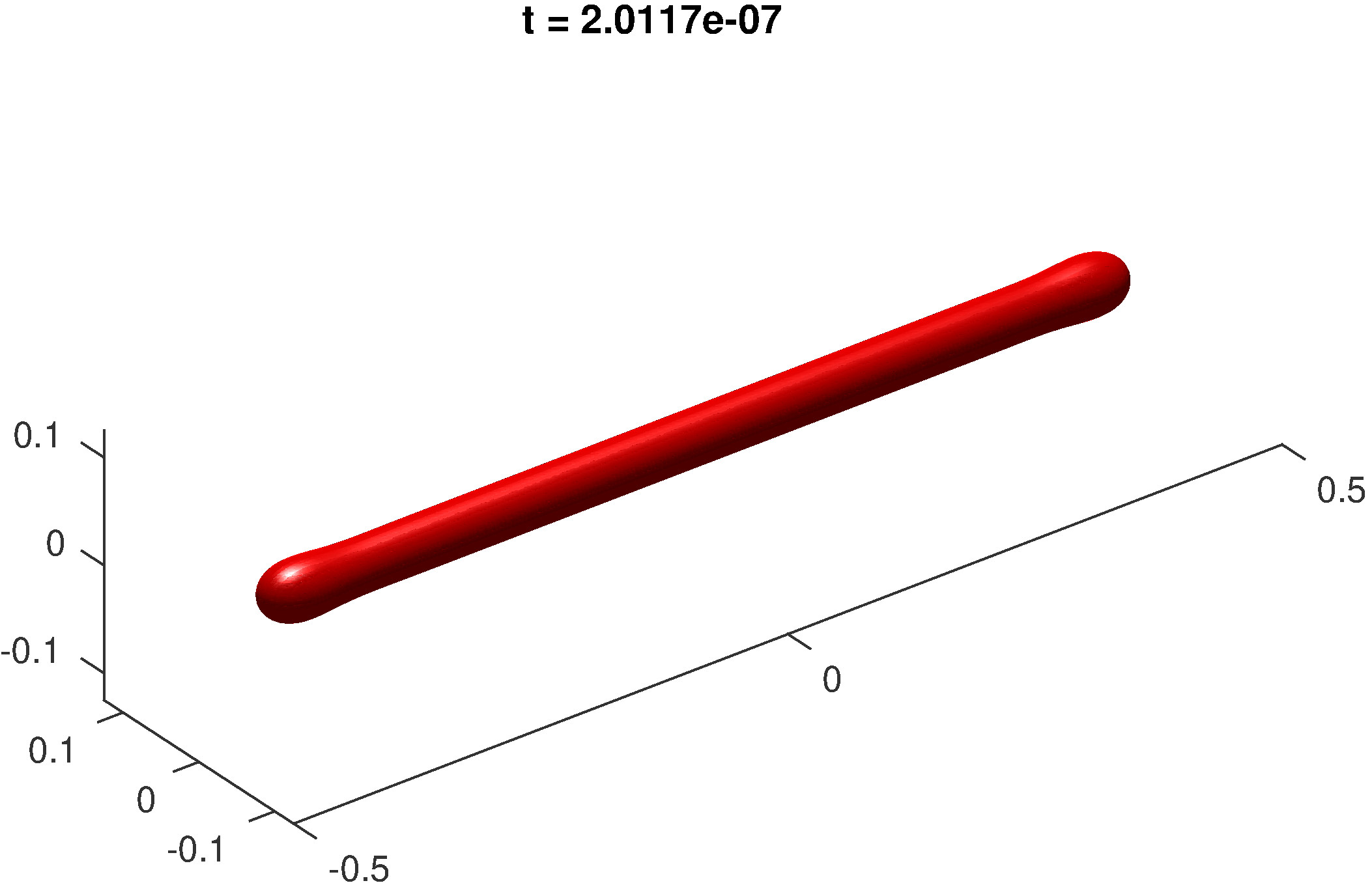}
	\includegraphics[width=3.5cm]{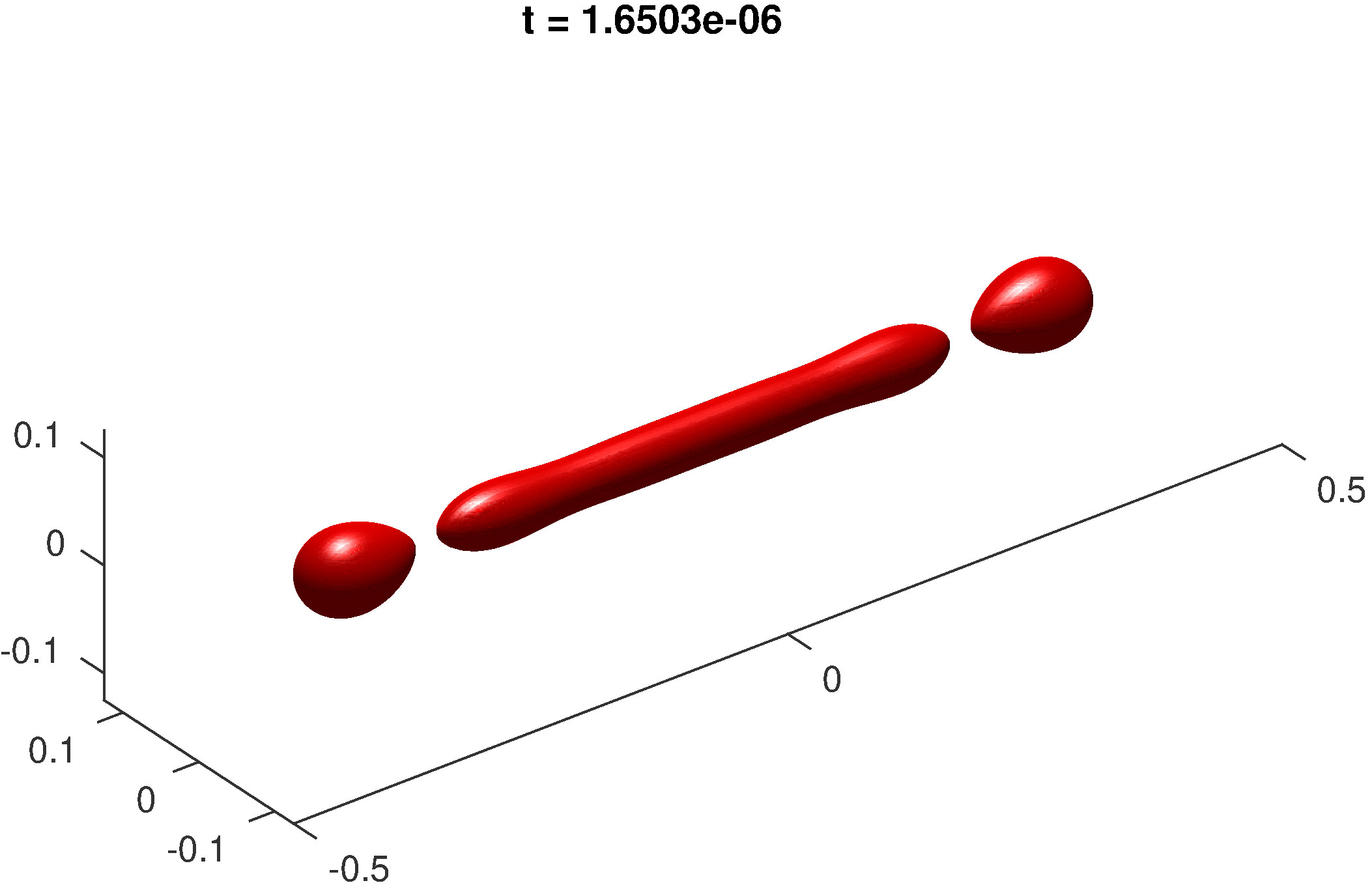}
	\includegraphics[width=3.5cm]{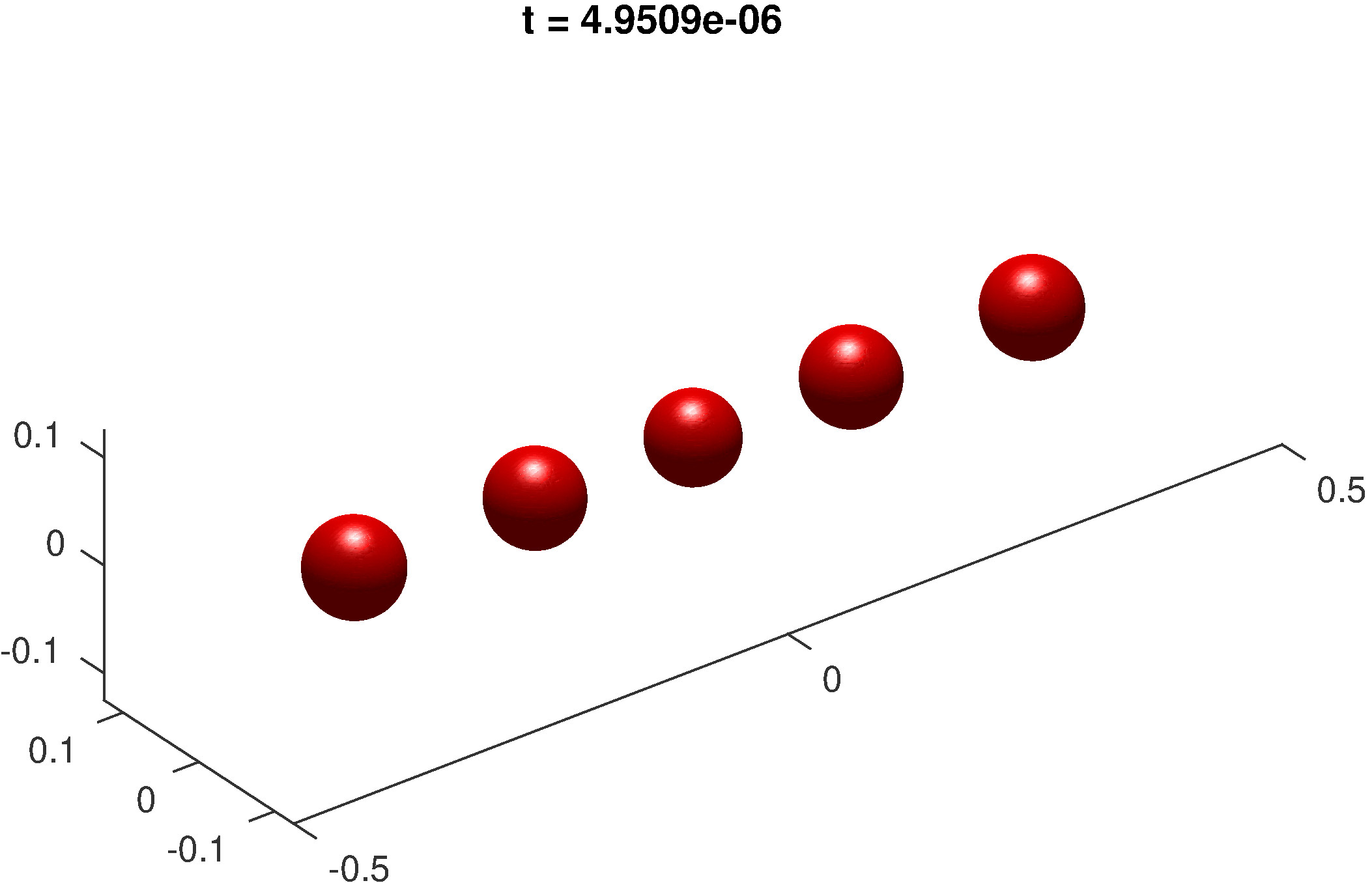} \\
	\caption{Comparison of the different models in the case of a thin structure in dimension $3$. First line corresponds to the \textbf{C-CH} model, second line to the \textbf{M-CH} model, and third line to the \textbf{NMN-CH} model.}
\label{fig_test5}
\end{figure}

We propose now a numerical experiment in dimension $3$ where the initial set is a thin tube. 
Our motivation here is to show the importance of having a model of order 2 in the phase field function $u$ in the complicated case of a thin structure evolution.
Similarly to the previous computations, the numerical parameters are given by $\delta_x = \frac{1}{2^8}$, $\epsilon = 2/N$, $\delta_t = \epsilon^4$,
$\alpha = 2/\epsilon^2$, $m=1$, and $\beta= 2/\epsilon^2$. We plot on each picture of \eqref{fig_test5} the $1/2$-level set of $u$ 
for different times $t$. The first, second and third line correspond, respectively, to the {\bf C-CH}, {\bf M-CH} and {\bf NNM-CH} models. 
We observe that the evolutionary set disappears using the {\bf C-CH} and {\bf M-CH} models whereas the {\bf NMN-CH} model seems to 
 have better volume conservation properties and the stationary set  is given as the sum of five small spheres. 
 
\begin{figure}[htbp]
\centering
	\includegraphics[width=5cm]{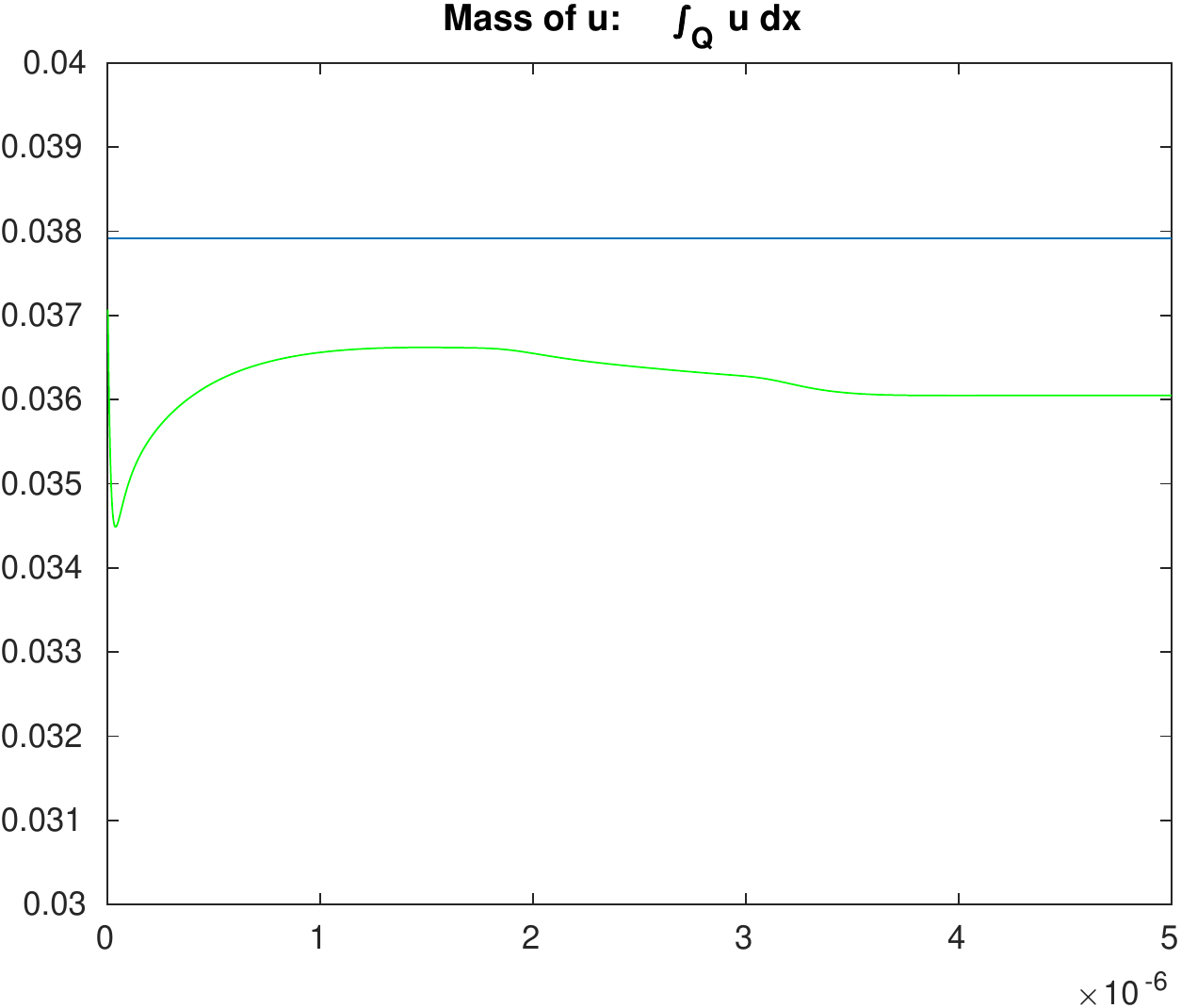}
	\caption{Comparison of the  different models in the case of a thin structure in dimension $3$. Evolution of the mass of $u$ given by $\int u dx$ along the iterations;
	Using the {\bf C-CH} model in red, the {\bf M-CH} model in blue, and the {\bf NMN-CH} model in green.}   
\label{fig_test5_vol}
\end{figure}

The results are surprising  at first glance as the mass of $u$ ($\int_Q u dx$) is well preserved using  the {\bf C-CH} and {\bf M-CH} models. 
So, to convince oneself that the problem arises from the phase field model order and not the numerical discretization,
we plot on figure \eqref{fig_test5_vol} the numerical evolution of the mass $t \mapsto \int_Q u dx$ along the flow for each model. We observe a very good conservation in the case of {\bf C-CH} and {\bf M-CH} models despite the disappearance of the structure. \\
Moreover, recall that we plot on figure \eqref{fig_test5} the $1/2$-level set of $u$:
$$\Omega_{\epsilon}(t) =  \left\{ x \in Q ;  u(x,t) \leq 1/2 \right\},$$
and that for a phase field model of order $1$ only, we have 
$$ Vol(\Omega_{\epsilon}(t))  =  \int_Q u(x,t) dx  + O(\epsilon).$$
This means that even if the mass of $u$ is conserved, we observe an error of order $O(\epsilon)$ on the volume of $\Omega_{\epsilon}$. The consequence is all the more dramatic in our example as the volume of the thin structure is of order $\epsilon^2$. In the end, the whole volume is lost because of this approximation error. 
Concerning the {\bf NMN-CH} model, we proved a volume approximation of order $2$,
$$ Vol(\Omega_{\epsilon}(t))  =  \int_Q u(x,t) dx  + O(\epsilon^2),$$
This explains the good numerical behavior of the {\bf NMN-CH} model in comparison with the other models.

In conclusion, this 3D numerical experiment showcases the inefficiency of models {\bf C-CH} and {\bf M-CH} to approximate the evolution
of a thin structure, where a much smaller $\epsilon$ is required. On the other hand, the second order {\bf NMN-CH}
phase field model 
seems to give a good approximation of surface diffusion even if the mass of $u$ is not perfectly conserved 
(Green plot on Figure \eqref{fig_test5_vol}). 

\subsubsection{Dewetting and surface diffusion of a thin plate  }

\begin{figure}[htbp]
\centering
	\includegraphics[width=3.5cm]{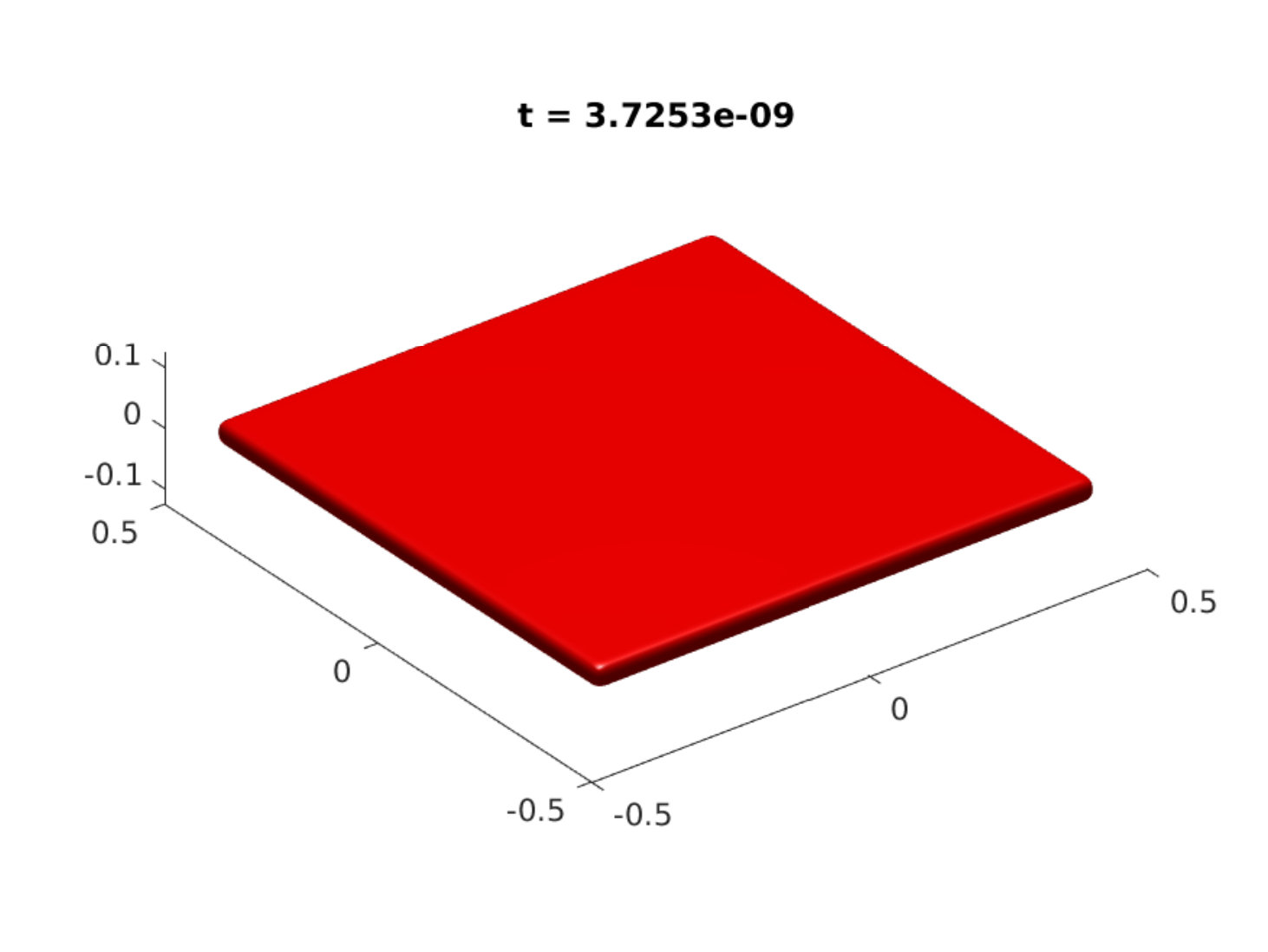}
	\includegraphics[width=3.5cm]{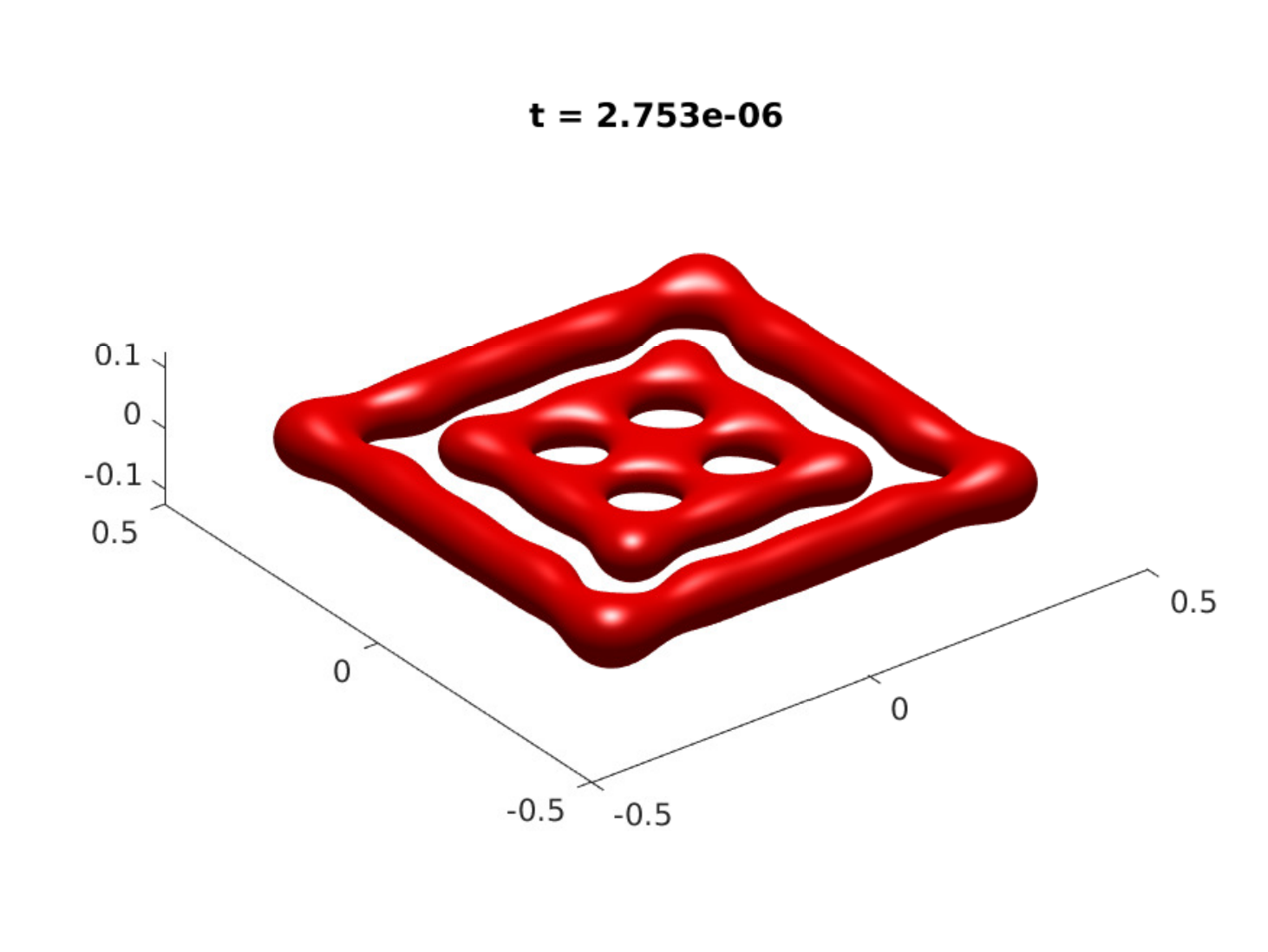} 
	\includegraphics[width=3.5cm]{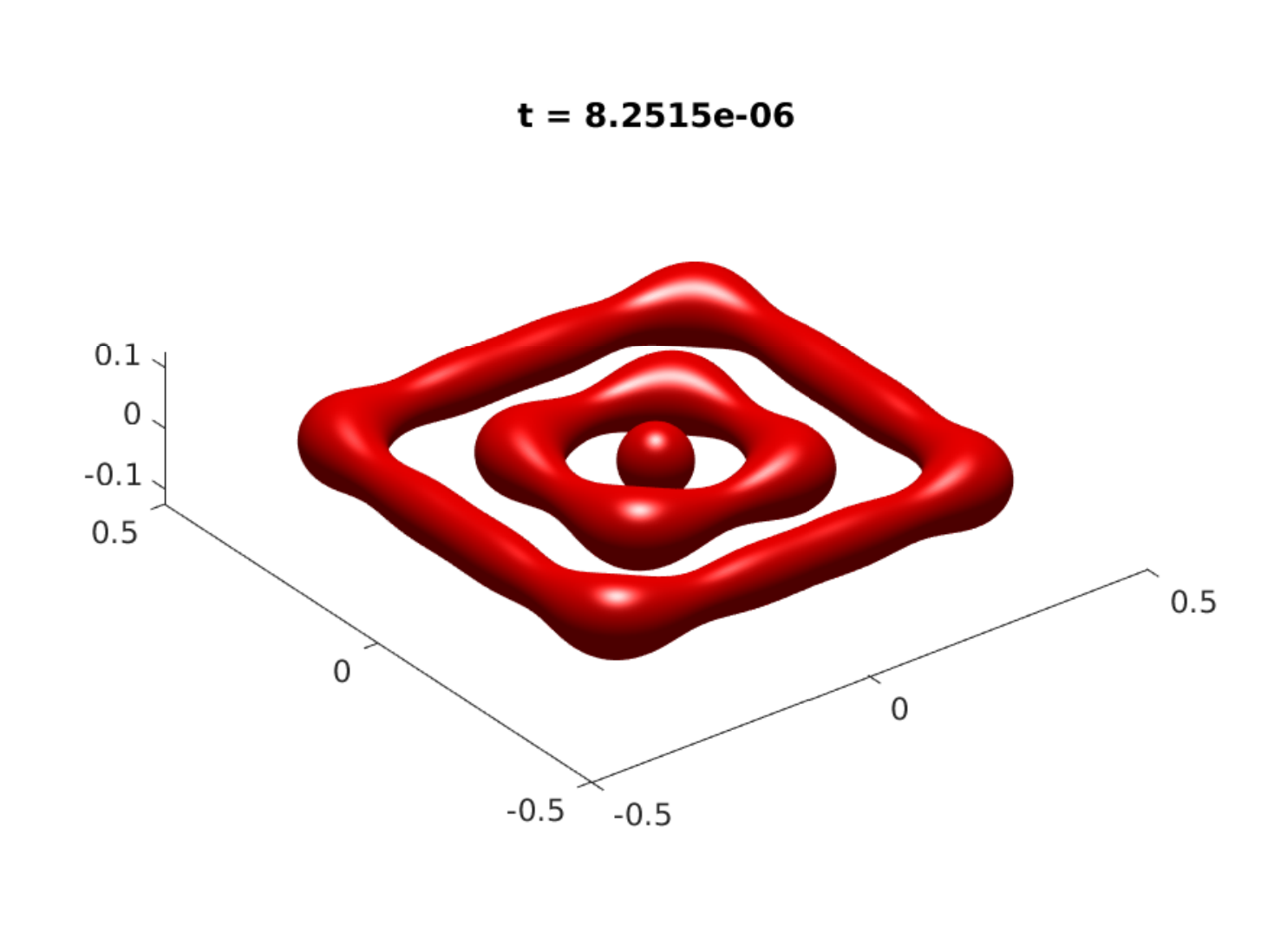}
	\includegraphics[width=3.5cm]{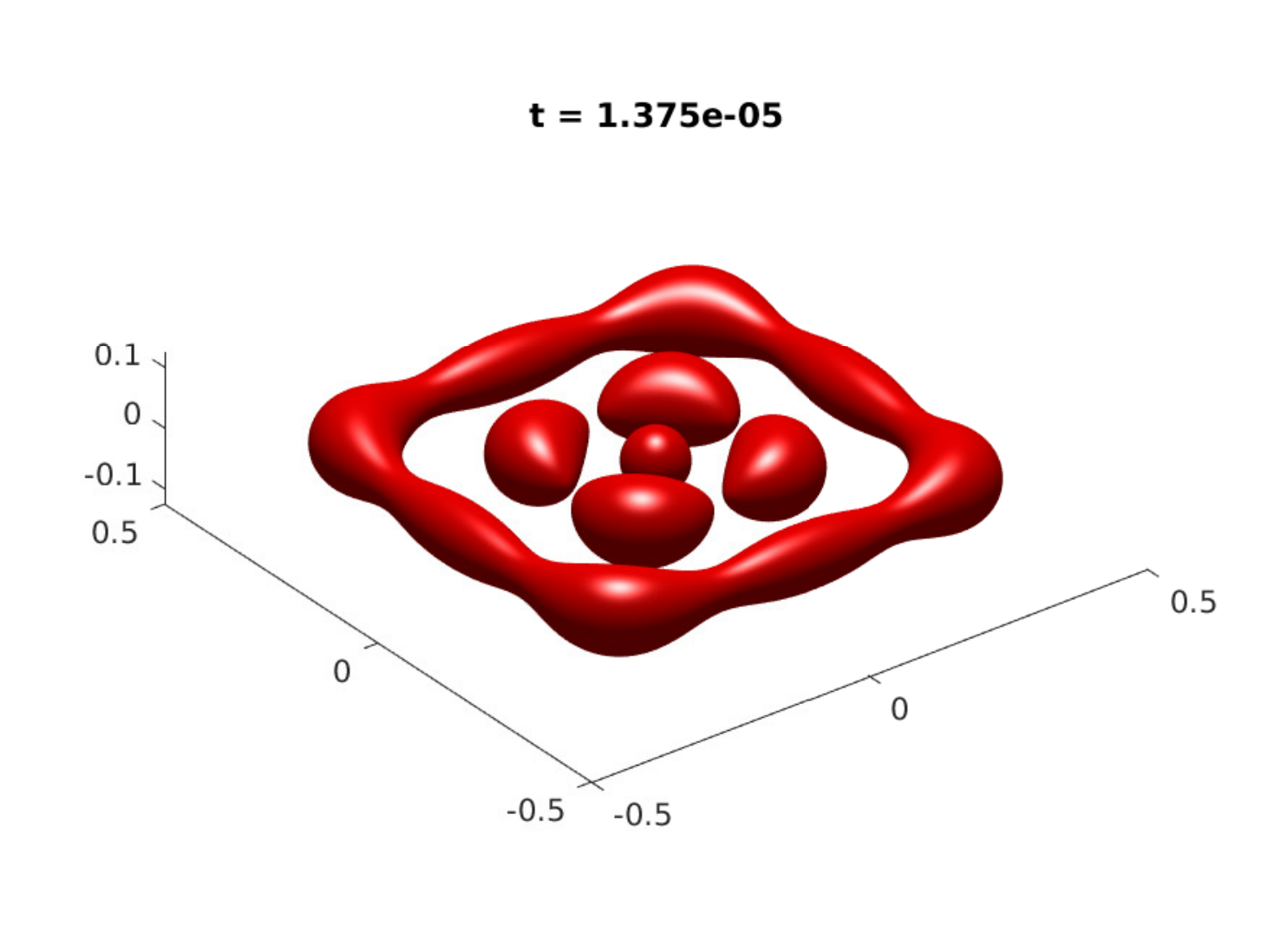} \\
	\includegraphics[width=3.5cm]{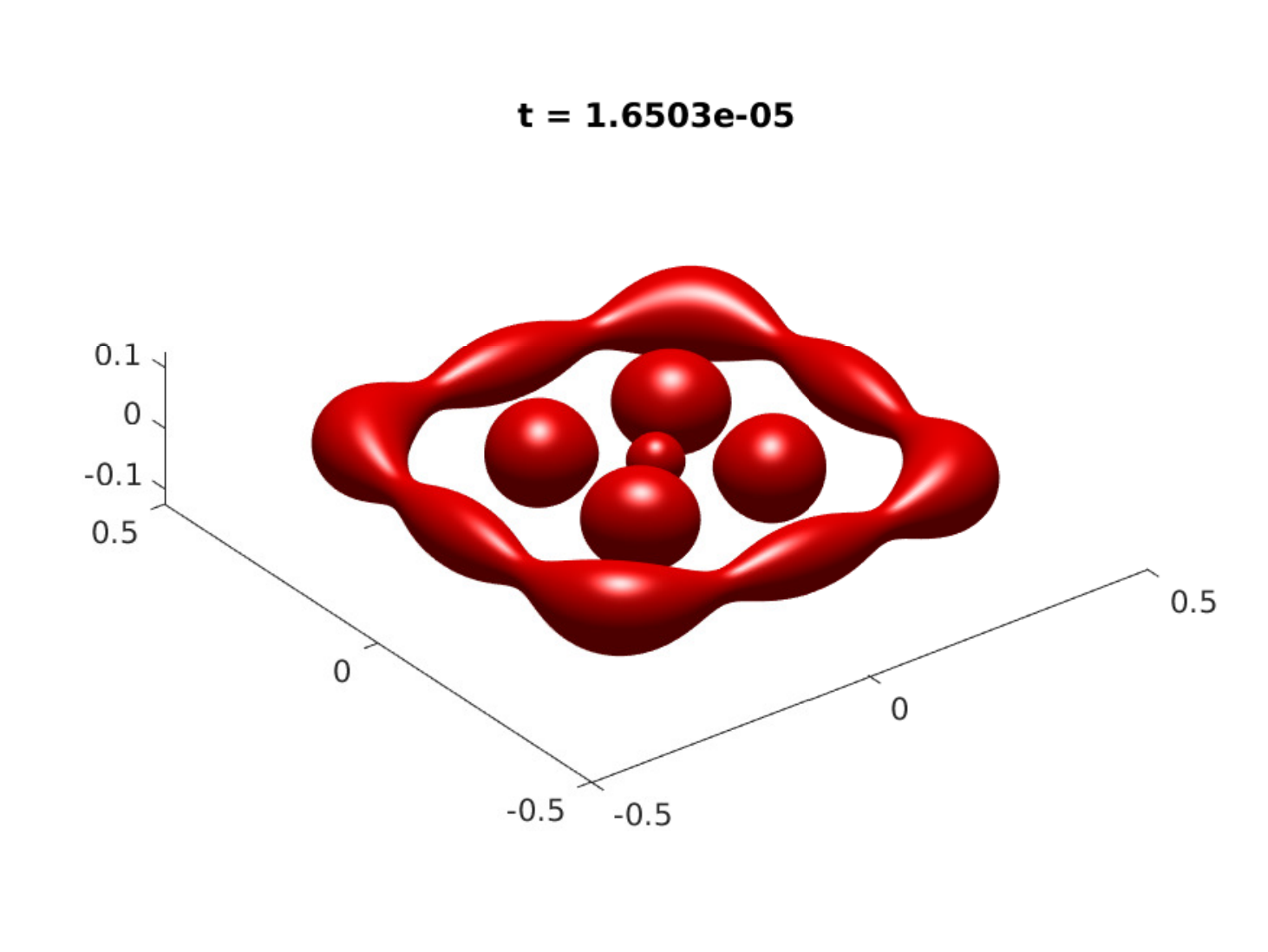}
    \includegraphics[width=3.5cm]{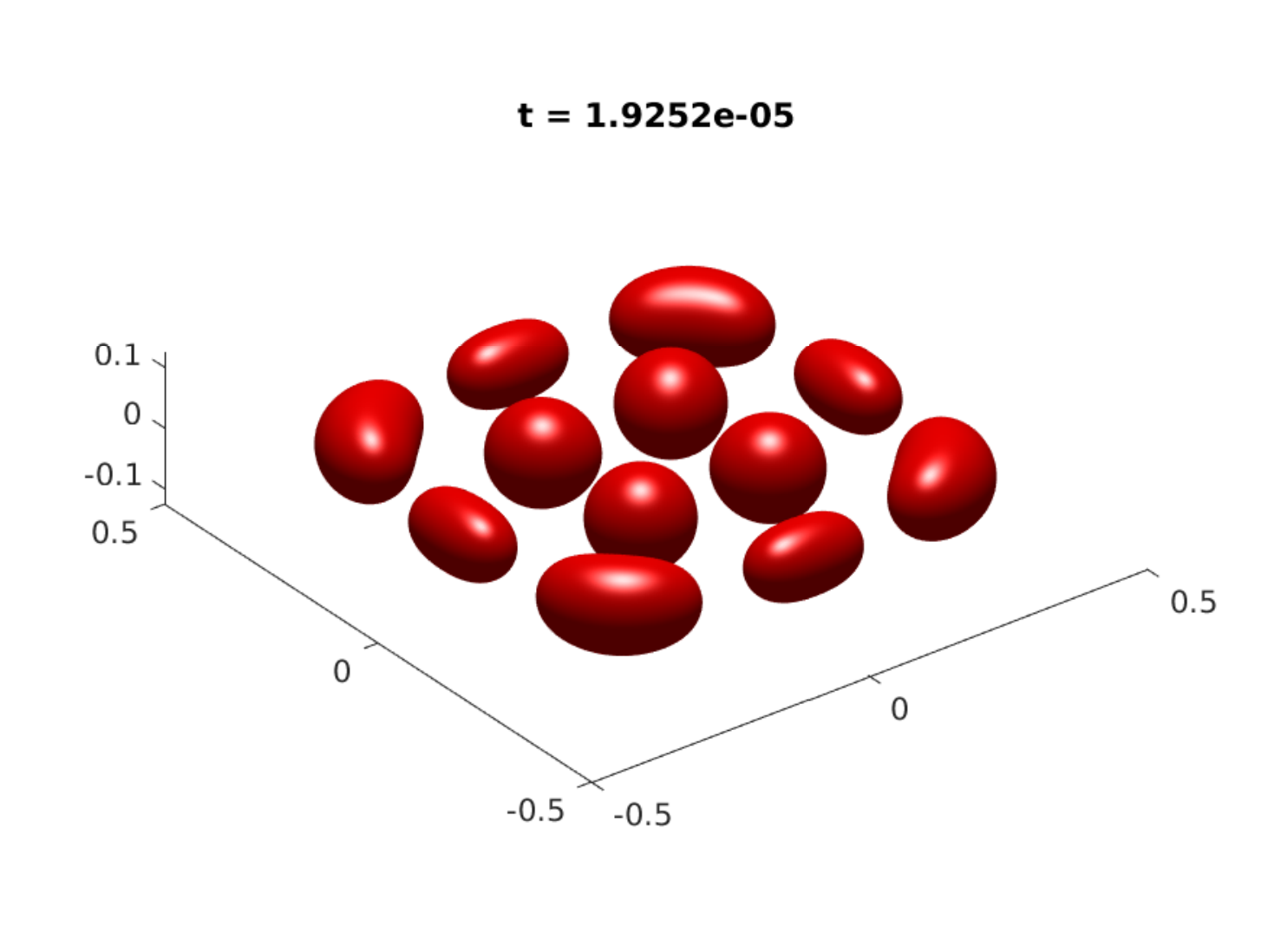}
    \includegraphics[width=3.5cm]{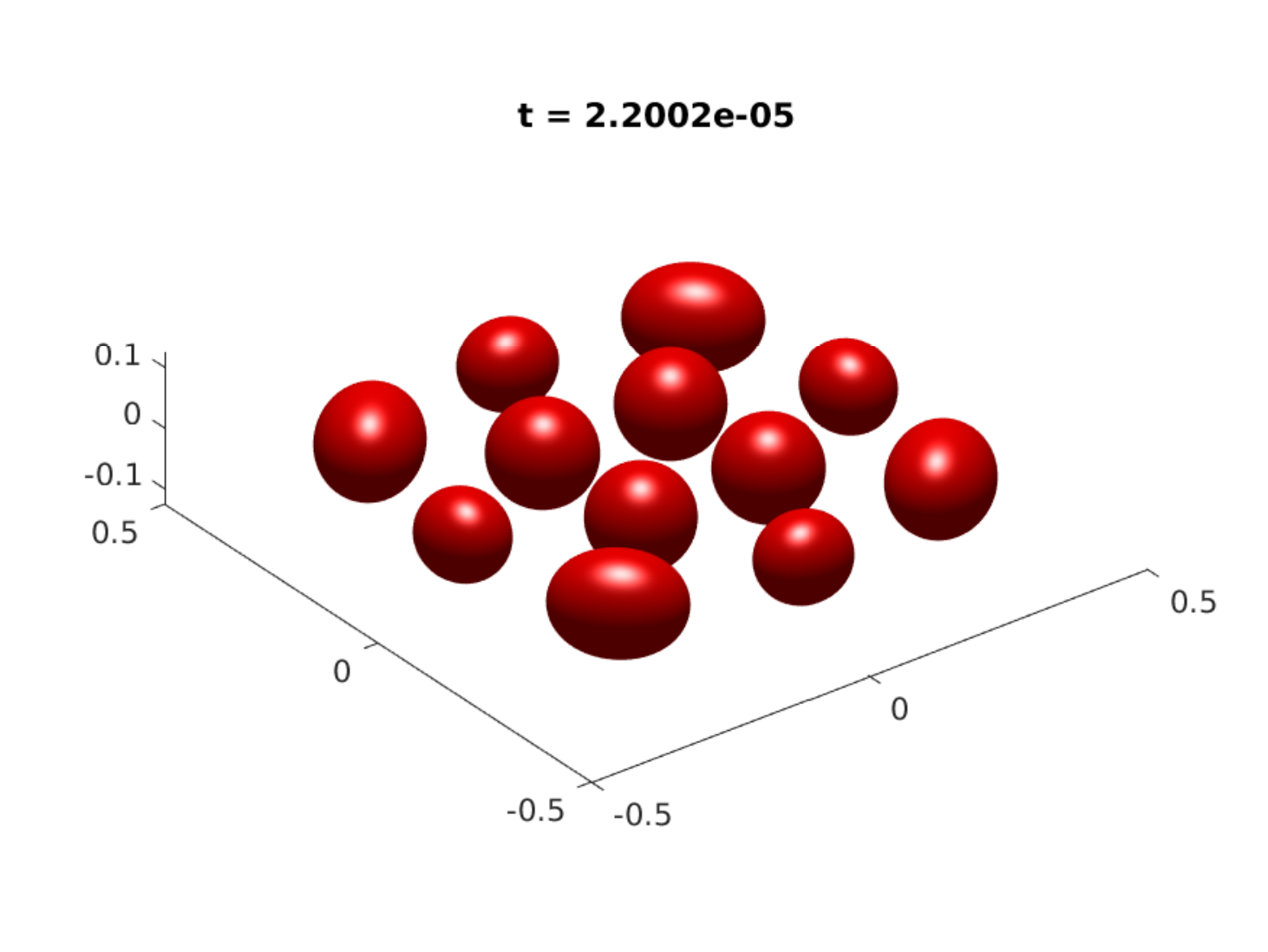} 
    \includegraphics[width=3.5cm]{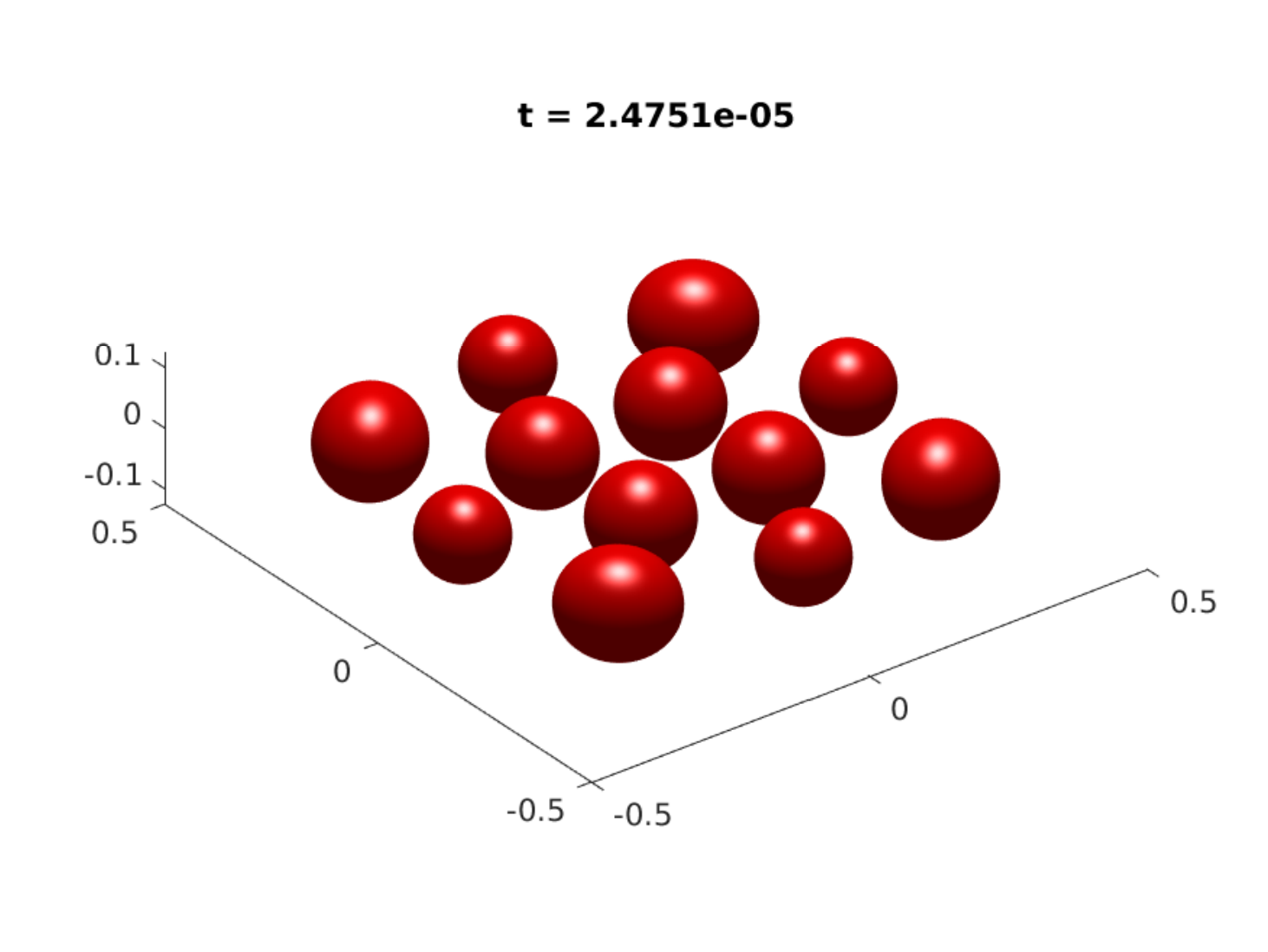} 
	
\caption{Example of dewetting in dimension 3 using the {\bf NMN-CH} model. Evolution of $u$ along the iterations.}
\label{fig_test4}
\end{figure}

The last numerical example is the evolution of a thin plate using the {\bf NMN-CH} model. As previously, 
the parameters are chosen as $\delta_x = \frac{1}{2^8}$, $\epsilon = 2/N$, $\delta_t = \epsilon^4$, $\alpha = 2/\epsilon^2$, $m=1$, and $\beta= 2/\epsilon^2$. 
We can observe on figure \eqref{fig_test4} an evolution similar to the one observed in real dewetting experiments\cite{MR3874087}.

\section*{Acknowledgment}
The authors thank Roland Denis for fruitful discussions. They acknowledge support from the French National Research Agency (ANR) under grants ANR-18-CE05-0017 (project BEEP) and  ANR-19-CE01-0009-01 (project MIMESIS-3D). Part of this work was also supported by the LABEX MILYON (ANR-10-LABX-0070) of Universit\'e de Lyon, within the program "Investissements d'Avenir" (ANR-11-IDEX- 0007) operated by the French National Research Agency (ANR).


\end{document}